\documentclass[11pt]{amsart}

\usepackage[margin=1in]{geometry}
\usepackage{amsmath,amsthm,amssymb}
\usepackage{color}
\usepackage{soul}
\usepackage{mathrsfs}
\usepackage{hyperref}
\usepackage{constants}
\newconstantfamily{abcon}{symbol=c}

\usepackage[capitalise]{cleveref}

\newtheorem{theorem}{Theorem}
\newtheorem{lemma}[theorem]{Lemma}
\newtheorem{proposition}[theorem]{Proposition}
\newtheorem{corollary}[theorem]{Corollary}

\newtheorem*{hypothesis}{Hypothesis}
\newtheorem{conjecture}[theorem]{Conjecture}
\theoremstyle{remark} \newtheorem*{remark}{Remark}

\newtheoremstyle{named}{}{}{\itshape}{}{\bfseries}{.}{.5em}{\thmnote{#3}}
\theoremstyle{named}
\newtheorem{namedtheorem}{}

\numberwithin{equation}{section}
\numberwithin{theorem}{section}

\usepackage{stackengine}

\newcommand{\kp}{\mathfrak{p}}
\newcommand{\Gal}{\mathrm{Gal}}
\newcommand{\re}{\mathrm{Re}}
\newcommand{\im}{\mathrm{Im}}
\newcommand{\kq}{\mathfrak{q}}
\newcommand{\Q}{\mathbb{Q}}
\newcommand{\R}{\mathbb{R}}
\newcommand{\kn}{\mathfrak{n}}
\newcommand{\kD}{\mathfrak{D}}
\newcommand{\N}{\mathrm{N}}

\renewcommand{\epsilon}{\varepsilon}
\renewcommand{\hat}{\widehat}
\definecolor{dartmouthgreen}{rgb}{0.0,0.412,0.243}
\definecolor{orange}{rgb}{1,0.5,0}

\renewcommand{\tilde}{\widetilde}

\title[An approximate form of  Artin's holomorphy conjecture and non-vanishing]{An approximate form of Artin's holomorphy conjecture and non-vanishing of Artin $L$-functions}

\author{Robert J. Lemke Oliver}
\address{Department of Mathematics, Tufts University, Medford, MA 02155}
\email{robert.lemke$\_{ }$oliver@tufts.edu}

\author{Jesse Thorner}
\address{Department of Mathematics, University of Illinois, Urbana, IL 61801}
\email{jesse.thorner@gmail.com}

\author{Asif Zaman}
\address{Department of Mathematics, University of Toronto, Toronto, Ontario, Canada M5S 2E4}
\email{zaman@math.toronto.edu}

\begin{document}

\begin{abstract}
Let $k$ be a number field and $G$ be a finite group.  Let $\mathfrak{F}_{k}^{G}(Q)$ be the family of number fields $K$ with absolute discriminant $D_K$ at most $Q$ such that $K/k$ is normal with Galois group isomorphic to $G$.  If $G$ is the symmetric group $S_n$ or any transitive group of prime degree, then  we unconditionally prove that for all $K\in\mathfrak{F}_k^G(Q)$ with at most $O_{\epsilon}(Q^{\epsilon})$ exceptions, the $L$-functions associated to the faithful Artin representations of $\Gal(K/k)$ have a region of holomorphy and non-vanishing commensurate with predictions by the Artin conjecture and the generalized Riemann hypothesis.  This result is a special case of a more general theorem.  As applications, we prove that:
\begin{enumerate}
	\item there exist infinitely many degree $n$ $S_n$-fields over $\Q$ whose class group is as large as the Artin conjecture and GRH imply, settling a question of Duke;
	\item for a prime $p$, the periodic torus orbits attached to the ideal classes of almost all totally real degree $p$ fields $F$ over $\Q$ equidistribute on $\mathrm{PGL}_p(\mathbb{Z})\backslash\mathrm{PGL}_p(\R)$ with respect to Haar measure; 
	\item for each $\ell\geq 2$, the $\ell$-torsion subgroups of the ideal class groups of almost all degree $p$ fields over $k$ (resp. almost all degree $n$ $S_n$-fields over $k$) are as small as GRH implies; and
	\item an effective variant of the Chebotarev density theorem holds for almost all fields in such families.
\end{enumerate}
\end{abstract}

\maketitle

\vspace{-0.5in}

\section{Introduction} \label{sec:intro}

Let $K/k$ be a normal extension of number fields with Galois group $G$.  Many natural arithmetic properties of $K$ are controlled by the Artin $L$-functions $L(s,\rho)$ attached to the irreducible complex representations $\rho$ of $G$.  The Artin conjecture asserts that $L(s,\rho)$ is entire for every nontrivial $\rho$, and the generalized Riemann hypothesis (GRH) asserts that $L(s,\rho)\neq 0$ for $\re(s)>\frac{1}{2}$.  Artin's conjecture is known for very few groups $G$, and GRH remains open.  In this paper, we substantially enlarge the region of holomorphy and non-vanishing for Artin $L$-functions in an average sense.

We begin with work of Aramata, Artin, and Brauer towards the Artin conjecture:  the quotient $\zeta_K(s)/\zeta_k(s)$ of Dedekind zeta functions is entire, and for each Artin representation $\rho$ of $\Gal(K/k)$, the $L$-function $L(s,\rho)$ is holomorphic and non-vanishing in any region where $\zeta_K(s)\neq 0$.  In many applications, existing zero-free regions are not strong enough to deduce the desired results.  Therefore, one might hope to average over number fields $K$ in a family $\mathfrak{F}$ and prove that apart from a small collection of $K\in\mathfrak{F}$, the ratios $\zeta_K(s)/\zeta_k(s)$ have a large zero-free region that can be used to deduce a large region of holomorphy and non-vanishing for each $L(s,\rho)$ associated with $K/k$.  We fix a finite nontrivial group $G$ and consider the family $\mathfrak{F}_k^G$ of normal extensions $K/k$ such that $\Gal(K/k)\cong G$.  When $G$ is abelian, each nontrivial representation is one-dimensional and can be realized as a Hecke character whose $L$-function is entire, so we focus on nonabelian groups $G$.

Two obstacles quickly emerge.  First, we do not yet know that the entire $L$-function $\zeta_K(s)/\zeta_k(s)$ factors as a product of $L$-functions associated to cuspidal automorphic representations defined over $k$, as the strong Artin conjecture asserts.  This hinders most methods of averaging over $K\in\mathfrak{F}_k^G$.  
The second obstacle, which we loosely term the {\it subfield problem}, arises from the work of Aramata, Artin, and Brauer---if $F \subseteq K$ is a subextension of $k$, then a zero of $\zeta_F(s)$ is a zero of $\zeta_K(s)$.  Moreover, the same field $F$ may arise as a subfield of many different $K$ in the family $\mathfrak{F}_{k}^G$.  For example, this is the case for the family of $S_n$-extensions of $k$, where different fields may share a common quadratic subfield; see \eqref{eqn:structural} below.  Consequently, any problematic zero of a single $\zeta_F(s)$ could propagate to many different $\zeta_K(s)$, so our ability to study the zeros of the ratios $\zeta_K(s)/\zeta_k(s)$ as we average over $K\in\mathfrak{F}_{k}^G$ is both technically limited by one's ability to control the frequency with which fields $K \in \mathfrak{F}_{k}^G$ intersect and structurally limited by the fact that such intersections do occur.

Pierce, Turnage--Butterbaugh, and Wood \cite{PTW} were the first to codify the subfield problem in the context of modern arithmetic statistics by relating it to other, more understood problems.  Additionally, when $k=\Q$, they use unconditional field counting results for certain nonabelian groups $G$, including the symmetric groups $S_n$ for $n\in\{3,4,5\}$ and the dihedral groups $D_{p}$ with $p$ an odd prime, to control the subfield problem by restricting to $K\in\mathfrak{F}_{\Q}^G$ satisfying certain ramification conditions.  Apart from $G=S_5$, the strong Artin conjecture is known for each of these groups, so the quotient $\zeta_K(s)/\zeta_\mathbb{Q}(s)$ factors in terms of automorphic $L$-functions.  Using a zero density estimate for automorphic $L$-functions due to Kowalski and Michel \cite{KM}, they then proved that for almost all of the $K\in\mathfrak{F}_{\Q}^G$ satisfying certain ramification conditions, the ratio $\zeta_K(s)/\zeta_{\Q}(s)$ and each Artin $L$-function in the associated factorization enjoy a very wide zero-free region.

Avoiding unproven hypotheses on the automorphy of the factors of $\zeta_K(s)$, Thorner and Zaman \cite{ThornerZaman} proved for {\it any} finite group $G$ for which the subfield problem may be controlled that for almost all fields $K\in\mathfrak{F}_{\Q}^G$, $\zeta_K(s)/\zeta_\mathbb{Q}(s)$ has a much larger zero-free region than was known previously, commensurate with what GRH implies.  This work also further clarified the subfield problem, distilling it into a question about so-called ``intersection multiplicities'' (see \eqref{eqn:intersection_multiplicity} below).  
As described earlier, the classical work of Aramata, Artin, and Brauer then shows that each irreducible Artin $L$-function must be both non-vanishing and holomorphic in this region, independent of whether it is known to be automorphic.
Thus, the work in \cite{ThornerZaman} avoids assumptions of automorphy, which it does by proving a new large sieve for the family $\mathfrak{F}_{\Q}^G$, but it does not address the technical and structural limitations presented by the subfield problem.  In particular, their work avoids these limitations when $G$ is simple (where the subfield problem trivially disappears) but it cannot avoid them in most other families, including the natural situation when $G = S_n$. 

We offer a new approach to producing large regions of holomorphy and non-vanishing for Artin $L$-functions associated to almost all fields in a family of Galois extensions that simultaneously addresses the subfield problem and the absence of automorphy results.  The novelty lies in reducing these obstacles to a group theoretic computation that is tractable for groups like $S_n$, but which in full generality apparently requires the complete classification of finite simple groups.  Our approach has two independent components.  The first component changes the averaging process in \cite{ThornerZaman} in a way that avoids the subfield problem (\cref{subsec:component_1}).  This new average, however, no longer produces a result that is amenable to the classical results of Aramata, Artin, and Brauer.  The second component, therefore, is a new result in character theory that handles this complication by expressing the characters of certain Artin representations in terms of inductions of one-dimensional characters with restricted components, which neither Artin nor Brauer induction can address (\cref{subsec:component_2}).  We view these two components together as an approximate form of Artin's conjecture that also provides a strong zero-free region for almost all $K\in\mathfrak{F}_{k}^{G}$ for any number field $k$.

Before we describe our method, we give a representative example of what our approach can prove for faithful Artin representations associated to the fields $K\in \mathfrak{F}_k^{G}$, where $G$ is the symmetric group $S_n$ or a transitive group of prime degree.  There is little to no progress towards the Artin conjecture when the degree of $G$ is at least 5 or towards basic counting problems in arithmetic statistics when the degree of $G$ is at least 6.  Despite these setbacks, we prove: 
\begin{theorem}\label{thm:artin-pnt-intro}
Let $k$ be a number field.  Let $G$ be the symmetric group $S_n$ for some $n\geq 2$ or a transitive subgroup of $S_p$ for some prime $p$. Let $Q\geq 1$.  For all $\epsilon>0$, there exists an  effectively computable constant $\Cl[abcon]{main}=\Cr{main}(|G|,[k:\Q],\epsilon)>0$ such that for all except $O_{|G|,[k:\Q],\epsilon}(Q^\epsilon)$ normal extensions $K/k$ with $\mathrm{Gal}(K/k) \simeq G$ and absolute discriminant $D_K$ at most $Q$, each irreducible faithful Artin representation $\rho$ of $\mathrm{Gal}(K/k)$ satisfies 
	\[
	\Big|\sum_{ \mathrm{N}_{k/\mathbb{Q}} \mathfrak{p} \leq x} \mathop{\mathrm{tr}} \rho(\mathrm{Frob_\mathfrak{p}})\Big|\ll_{|G|,[K:\Q],\epsilon} x \exp(-\Cr{main} \sqrt{\log x})
	\]
	for all $ x \geq (\log D_K)^{81|G|/\epsilon}$.  The sum is over prime ideals of $k$ with absolute norm at most $x$.
\end{theorem}

\begin{remark}
Under GRH, the exponent $81|G|/\epsilon$ may be replaced by $2+\epsilon$ for all $K$.
\end{remark}

\begin{remark}
The constant $\Cr{main}$ is the same as in \cref{lem:error_term_data_bound} below.	
\end{remark}

\begin{remark}
We use the notation $f\ll_{\nu} g$ and $f=O_{\nu}(g)$ to denote the existence of an effectively computable constant $c_{\nu}>0$, depending at most on $\nu$, such that $|f|\leq c_{\nu}|g|$ in the range indicated.  We use the notation $f\asymp_{\nu}g$ to indicate that $f\ll_{\nu}g$ and $g\ll_{\nu}f$.
\end{remark}

\section{Applications}
\label{sec:Applications}

\cref{thm:artin-pnt-intro} and the ideas leading to it produce several desirable arithmetic results.  We will now sample a few of them.  In \cref{sec:approximate-artin}, we will discuss the ideas leading to \cref{thm:artin-pnt-intro}.

\subsection{The extremal order of class numbers}
\label{ssec:extremal}

Let $\mathcal{K}_n$ be the family of totally real number fields $F$ with $[F:\Q]=n$ whose normal closure over $\Q$ has the full symmetric group $S_n$ as its Galois group.  Using lower bounds on the regulator of $F$ due to Remak \cite{Remak}, Duke \cite{Duke03} proved under GRH and the Artin conjecture that if $F\in\mathcal{K}_n$, then
\[
|\mathrm{Cl}(F)|
		\ll_{n} \frac{D_F^{1/2} (\log\log D_F)^{n-1}}{(\log D_F)^{n-1}},
\]
where $D_F$ is the absolute discriminant of $F$.  Furthermore, still under the assumption of GRH and the Artin conjecture, Duke showed that this upper bound is sharp, in that
\begin{equation}
\label{eqn:duke}
\begin{aligned}
\textup{there exist $F\in\mathcal{K}_n$ with $D_F$ arbitrarily large and }|\mathrm{Cl}(F)|
		\asymp_{n} \frac{D_F^{1/2} (\log\log D_F)^{n-1}}{(\log D_F)^{n-1}}.
\end{aligned}
\end{equation}
The conclusion \eqref{eqn:duke} is proved without recourse to unproven hypotheses when $n\in\{2,3,4\}$ \cite{Cho,Daileda,MontgomeryWeinberger}.  Cho \cite{Cho} proved that \eqref{eqn:duke} holds when $n\geq 5$ using only the strong Artin conjecture, removing the reliance of Duke's argument on GRH.  We use Theorem \ref{thm:artin-pnt-intro} to prove the following unconditional result that removes the hypotheses of Artin's conjecture and GRH as well as the requirement that $F$ be totally real.

\begin{theorem} \label{thm:extremal-class-number}
For any fixed integers $r_1,r_2\geq 0$ with $n:=r_1+2r_2\geq 2$, there are number fields $F$ of signature $(r_1,r_2)$ with arbitrarily large discriminant $D_F$ whose normal closure has the full symmetric group $S_{n}$ as its Galois group, for which
\begin{equation}
\label{eqn:maximal_class}
	|\mathrm{Cl}(F)|
		\asymp_{r_1,r_2} \frac{D_F^{1/2} (\log\log D_F)^{r_1+2r_2-1}}{(\log D_F)^{r_1+r_2-1}}.
\end{equation}
\end{theorem}
\begin{remark}
Fix $0<\tau<1/(n^2-n)$.  Our proof shows that there exists a constant $\Cl[abcon]{large_class_lower_bound}=\Cr{large_class_lower_bound}(n,\tau)>0$ such that if $Q\geq \Cr{large_class_lower_bound}$, then there are at least $Q^{\tau}$ number fields $F$ that satisfy the conclusion of \cref{thm:extremal-class-number}.  See \cref{thm:extremal-class-number-body}.
\end{remark}

\subsection{Distribution of periodic torus orbits and subconvexity}

Let $F/\mathbb{Q}$ be a totally real field of degree $n$ with ring of integers $\mathcal{O}_F$.  Then $F$ may be naturally embedded into $\mathbb{R}^n$ by the product of its real embeddings, and in this embedding, the integers 
$\mathcal{O}_F$ form a full rank lattice.  More generally, a subset $\Lambda \subseteq F$ is a lattice if it is a free $\mathbb{Z}$-submodule of rank $n$.  The $F$-equivalence class of $\Lambda$, or the $F$-homothety class, is the set of lattices $\Lambda^\prime \subseteq F$ for which $\Lambda^\prime = \alpha \Lambda$ for some $\alpha \in F^\times$.  If we let $\mathcal{O} = \mathcal{O}_\Lambda := \{ \alpha \in F : \alpha \Lambda \subseteq \Lambda \}$, then $\mathcal{O}$ is an order in $\mathcal{O}_F$, and two equivalent lattices have the same associated order $\mathcal{O}$.  Moreover, there is a representative of the class of $\Lambda$ that is an ideal $\mathfrak{a}$ in $\mathcal{O}$, and the set of such ideal representatives constitute the ideal class of $\mathfrak{a}$ in $\mathcal{O}$.  Thus, equivalence classes of lattices in $F$ are naturally identified with ideal classes in orders $\mathcal{O} \subseteq \mathcal{O}_F$.

The space of lattices in $\mathbb{R}^n$ is naturally identified with $\mathrm{GL}_n(\mathbb{Z}) \backslash \mathrm{GL}_n(\mathbb{R})$.  By considering its action via multiplication, $F^\times$ embeds into the maximal split torus $H_n \subseteq \mathrm{GL}_n(\mathbb{R})$ consisting of diagonal matrices.  Thus, $F$-equivalent lattices give rise to elements of the manifold $\mathrm{PGL}_n(\mathbb{Z}) \backslash \mathrm{PGL}_n(\mathbb{R})$ equivalent under the action of $H_n$.  By the above discussion, we may think of an ideal class in an order $\mathcal{O} \subseteq \mathcal{O}_F$ as parametrizing a full $H_n$-orbit in $\mathrm{PGL}_n(\mathbb{Z}) \backslash \mathrm{PGL}_n(\mathbb{R})$.  In a pair of papers \cite{ELMV-Duke,ELMV-Annals}, Einsiedler, Lindenstrauss, Michel, and Venkatesh showed that every closed $H_n$-orbit on $\mathrm{PGL}_n(\mathbb{Z}) \backslash \mathrm{PGL}_n(\mathbb{R})$ arises as a periodic torus orbit in this manner.  When $n$ is prime, they connected the equidistribution of these torus orbits as the discriminant $\mathrm{disc}(\mathcal{O}) \to \infty$ to the problem of proving a discriminant-aspect subconvexity bound for $\zeta_F(s)$ of the form
\begin{equation}
	\label{eqn:subconvexity_dedekindzeta}
	|\zeta_F(\tfrac{1}{2}+it)|\ll_{[F:\Q]}D_F^{\frac{1}{4}-\theta}(1+|t|)^A,
\end{equation}
where $\theta\in(0,\frac{1}{4})$ and $A>0$ are constants that depend at most on $[F:\Q]$.  In the case $n=2$, this leads to a reinterpretation of Duke's theorem \cite{Duke88} on equidistribution of geodesics on the modular curve associated to real quadratic fields.  When  $n=3$, there was enough progress toward subconvexity that Einsieder, Lindenstrauss, Michel, and Venkatesh could prove that the analogous equidistribution result holds on $\mathrm{PGL}_3(\mathbb{Z}) \backslash \mathrm{PGL}_3(\mathbb{R})$.

Despite tremendous progress on proving subconvexity bounds for various families of automorphic $L$-functions, the bound \eqref{eqn:subconvexity_dedekindzeta} is only known when $F$ is a normal extension over a fixed base field $k$ with either an abelian or generalized dihedral Galois group, as well as the case when $F$ is an arbitrary cubic extension of a fixed field $k$.   In these cases, the Dedekind zeta function factors as a product of standard $L$-functions associated to cuspidal automorphic representations of $\mathrm{GL}_1(\mathbb{A}_k)$ or $\mathrm{GL}_2(\mathbb{A}_k)$, and the bound \eqref{eqn:subconvexity_dedekindzeta} follows from work of Michel and Venkatesh \cite{MV} (see also \cite{BHM,Burgess,DFI}).

We produce many new number fields $F$ that are extensions of $k$ satisfying \eqref{eqn:subconvexity_dedekindzeta} even if we do not yet know whether $\zeta_F(s)$ factors into a product of $L$-functions that are automorphic over $k$.  To state our result, we introduce some notation.  Let $k$ be a number field, $p$ be a prime, $n\geq 2$ be an integer, and $Q\geq 1$.  We define the families
\begin{equation}
\label{eqn:families}
\begin{aligned}
\mathscr{F}_{k}^{p}:=\{F\colon [F:k]=p\},\qquad &\mathscr{F}_{k}^{p}(Q):=\{F\in\mathscr{F}_k^p\colon D_F\leq Q\}\\
\mathscr{F}_{k}^{n,S_n}:=\{F\colon [F:k]=n,~\mathrm{Gal}(\tilde{F}/k)\cong S_n\},\quad &\mathscr{F}_{k}^{n,S_n}(Q):=\{F\in\mathscr{F}_k^{n,S_n}\colon D_F\leq Q\},
\end{aligned}
\end{equation}
where $\tilde{F}$ is the Galois closure of $F$ over $k$.  The ideas leading to \cref{thm:artin-pnt-intro} (see \cref{subsec:holomorphy_non-vanishing}) enable us to prove the following unconditional result.

\begin{theorem}
\label{thm:subconvexity-nice-families}
Let $t\in\R$, $Q\geq 1$, and $k$ be a number field.  Let $p$ be prime and $n\geq 2$.
\begin{enumerate}
\item Let $\epsilon>0$.  For all except $O_{p,[k:\mathbb{Q}],\epsilon}(Q^\epsilon)$ of the fields $F \in \mathscr{F}_k^{p}(Q)$, we have
\[
|\zeta_F(\tfrac{1}{2}+it)| \ll_{p,[k:\Q]} D_k^{O(\epsilon)}D_F^{\frac{1}{4}(1 - \frac{\epsilon}{10^{10}(p!)^2})}(1+|t|)^{O(p[k:\Q])}.
\]
\item Let $\epsilon>0$.  For all except $O_{n,[k:\mathbb{Q}],\epsilon}(Q^\epsilon)$ of the fields $F \in \mathscr{F}_k^{n,S_n}(Q)$, we have
\[
|\zeta_F(\tfrac{1}{2}+it)| \ll_{n,[k:\Q]} D_k^{O(\epsilon)}D_F^{\frac{1}{4}(1 - \frac{\epsilon}{10^{10}(n!)^2})}(1+|t|)^{O(n[k:\Q])}.
\]
\end{enumerate}
\end{theorem}
\begin{remark}
Ellenberg and Venkatesh \cite[Theorem 1.1]{EV_counting} proved that there exist effectively computable constants $\Cl[abcon]{EV_lower_n}=\Cr{EV_lower_n}(n,k)>0$ and $\Cl[abcon]{EV_lower_p}=\Cr{EV_lower_p}(p,k)>0$ such that if $Q\geq 1$, then
\begin{equation}
\label{eqn:lower_bounds}
|\mathscr{F}_k^{n,S_n}(Q)|\geq \Cr{EV_lower_n}Q^{\frac{1}{2}},\qquad |\mathscr{F}_k^p(Q)|\geq \Cr{EV_lower_p}Q^{\frac{1}{2}}.
\end{equation}
This ensures that \cref{thm:subconvexity-nice-families} is not vacuous.
\end{remark}

Taking $k=\Q$ in part 1 of \cref{thm:subconvexity-nice-families}, we deduce our next result.

\begin{theorem} \label{thm:torus-orbits}
Let $p \geq 5$ be prime and let $\mathscr{F}_\Q^{p,+} \subseteq \mathscr{F}_{\Q}^p$ be the set of totally real degree $p$ extensions of $\mathbb{Q}$.  For any $\epsilon>0$, there exists a set $\mathscr{E}_{\epsilon}^p\subseteq\mathscr{F}_{\Q}^{p,+}$ such that
\begin{enumerate}
	\item $|\{F\colon F\in \mathscr{E}_{\epsilon}^p\cap \mathscr{F}_{\Q}^{p,+}(Q)\}|\ll_{p,\epsilon}Q^{\epsilon}$ for all $Q\geq 1$, and
	\item if $(\mathcal{O}_j)_{j=1}^{\infty}$ is a sequence of orders in $\{\mathcal{O}\colon \text{there exists $F \in \mathscr{F}^{p,+}_{\mathbb{Q}}-\mathscr{E}_{\epsilon}^p$ such that $\mathcal{O} \subseteq \mathcal{O}_F$}\}$ with $\lim_{j\to\infty}\mathrm{disc}(\mathcal{O}_j)=\infty$, then as $j \to \infty$, the union of $H_p$-orbits associated to the ideal classes of $\mathcal{O}_j$ described above becomes equidistributed with respect to Haar measure on $\mathrm{PGL}_p(\mathbb{Z}) \backslash \mathrm{PGL}_p(\mathbb{R})$.
\end{enumerate}
In particular, if $(F_j)_{j=1}^{\infty}$ is a sequence of fields in $\mathscr{F}_{\Q}^{p,+}-\mathscr{E}_{\epsilon}^p$ ordered by discriminant, then the measures $\mu_{F_j}$ on $\mathrm{PGL}_p(\mathbb{Z}) \backslash \mathrm{PGL}_p(\mathbb{R})$ associated to the $H_p$-orbits of the ideal classes of $\mathcal{O}_{F_j}$ converge to Haar measure on $\mathrm{PGL}_p(\mathbb{Z}) \backslash \mathrm{PGL}_p(\mathbb{R})$ in the weak-* limit as $j\to\infty$.
\end{theorem}

\begin{remark}
It follows from minor modifications to the work of Ellenberg and Venkatesh \cite{EV_counting} that for all primes $p$, there exists a constant $\Cl[abcon]{pconst}=\Cr{pconst}(p)>0$ such that for all $Q\geq 1$, we have $|\mathscr{F}_\Q^{p,+}(Q)| \geq\Cr{pconst} Q^{\frac{1}{2}}$.  Consequently, this result is not vacuous.

\end{remark}

\subsection{$\ell$-torsion in class groups}

Let $F/k$ be an extension of number fields, and let $\mathrm{Cl}(F)$ denote the class group of $F$.  For any integer $\ell \geq 2$ and any $\epsilon>0$, it is expected that the $\ell$-torsion subgroup $\mathrm{Cl}(F)[\ell]$ satisfies $|\mathrm{Cl}(F)[\ell]| \ll_{\ell,[F:\Q],\epsilon} D_F^\epsilon$ \cite{Duke98}.  This is known only for prime $\ell$ when the normal closure of $F$ has a Galois group that is an $\ell$-group \cite{KlunersWang}.  The trivial bound $|\mathrm{Cl}(F)[\ell]| \leq |\mathrm{Cl}(F)| \ll_{[F:\Q],\epsilon} D_F^{1/2+\epsilon}$ follows from Minkowski's bound.  Ellenberg and Venkatesh \cite[Proposition 3.1]{EV} showed that GRH implies for all $\epsilon>0$ the improvement
\begin{equation}
\label{eqn:EV_GRH}
|\mathrm{Cl}(F)[\ell]|\ll_{[F:\Q],\ell,\epsilon}D_F^{\frac{1}{2}-\frac{1}{2\ell([F:k]-1)}+\epsilon}.
\end{equation}
The fields $F$ for which there unconditionally exists a constant $\delta>0$ (depending at most on $\ell$ and $[F:k]$) such that $|\mathrm{Cl}(F)[\ell]|\ll_{[F:\Q],\ell,\epsilon}D_F^{1/2-\delta}$ are scarce \cite{BSTTTZ,EV,MR2220098,MR2254390,Wang-Nilpotent,Wang}.   The key to the $\ell$-torsion bounds in \cite{EV,Wang-Nilpotent,Wang} is a lemma of Ellenberg and Venkatesh \cite[Lemma 2.3]{EV} that exploits non-inert primes of norm at most $D_F^{1/2\ell([F:k]-1)}$.

The works of An \cite{Chen_An}; Ellenberg, Pierce, and Wood \cite{EPW}; Pierce, Turnage-Butterbaugh, and Wood \cite{PTW}; and Thorner and Zaman \cite{ThornerZaman} consider the problem of proving  that \eqref{eqn:EV_GRH} holds for  $k=\Q$  and all number fields $F$ with $D_F\leq Q$ in certain families, provided that an exceptional set of relative density zero is omitted.  Each of these results uses the work of Ellenberg and Venkatesh \cite[Lemma 2.3]{EV} to reduce the problem to the study of small primes that split completely in the fields under consideration.  Table 1 below summarizes the current progress that makes no recourse to unproven hypotheses (with $n\geq 2$ denoting an integer and $p$ denoting a prime).
{\small\begin{table}[h!]
\caption{On-average $\ell$-torsion results for extensions of $\Q$ from \cite{Chen_An,EPW,PTW,ThornerZaman}}
\label{table}
\begin{tabular}{|l|l|l|l|l|l|}
\hline Source & Galois structure & Restrictions & Family size & Exceptional set size\\
\hline \cite{PTW} & degree $n$ $\mathbb{Z}/n\mathbb{Z}$-fields & on tamely ramified primes & $\sim c_n Q^{\frac{1}{n-1}}$ & $Q^{\epsilon}$ \\
\hline \cite{Chen_An} & degree $4$ $D_{4}$-fields & none & $\sim b_4 Q$  & $Q^{\frac{1}{2}+\epsilon}$ \\
\hline \cite{PTW} & degree $p$ $D_{p}$-fields, $p\geq 3$ & on tamely ramified primes & $\gg_p Q^{\frac{2}{p-1}}$  & $Q^{\frac{1}{p-1}+\epsilon}$ \\
\hline \cite{EPW} & degree 3 $S_3$-fields & none & $\sim c_3 Q$ & $Q^{1-\frac{1}{4\ell}+\epsilon}$\\
\hline \cite{PTW} & degree 3 $S_3$-fields & squarefree discriminants & $\sim d_3 Q$ & $Q^{\frac{1}{3}+\epsilon}$\\
\hline \cite{EPW} & degree 4 $S_4$-fields & $\ell\geq 8$ &  $\sim c_4 Q$ & $Q^{1-\frac{1}{6\ell}+\epsilon}$ \\
\hline \cite{PTW} & degree 4 $S_4$-fields & squarefree discriminants & $\sim d_4 Q$ & $Q^{\frac{1}{2}+\epsilon}$\\
\hline \cite{EPW} & degree 5 $S_5$-fields &  $\ell\geq 25$ & $\sim c_5 Q$  & $Q^{1-\frac{1}{8\ell}+\epsilon}$\\
\hline \cite{ThornerZaman} & degree $n$ $A_n$-fields, $n\geq 5$ & none & $\gg_n Q^{\frac{1}{30}}$ & $Q^{\epsilon}$ \\
\hline
\end{tabular}
\end{table}}%

\noindent Using \cref{thm:artin-pnt-intro} in concert with \cite[Lemma 2.3]{EV} we obtain the following result.

\begin{theorem}\label{thm:ell-torsion}
Let $\ell \geq 2$ be an integer, $k$ be a number field, and $Q\geq 1$.  Let $p$ be prime, $n\geq 2$ be an integer, and the families $\mathscr{F}_k^p(Q)$ and $\mathscr{F}_k^{n,S_n}(Q)$ be as in \eqref{eqn:families}.
\begin{enumerate}
	\item Let $\epsilon>0$ and $\eta>0$.  For all except $O_{p,[k:\Q],\epsilon}(Q^{\epsilon})$ fields $F\in\mathscr{F}_k^p(Q)$, there holds
\[
	|\mathrm{Cl}(F)[\ell]|\ll_{p,k,\ell,\epsilon,\eta}D_F^{\frac{1}{2}-\frac{1}{2\ell([F:k]-1)}+\eta}.
\]
	\item Let $\epsilon>0$ and $\eta>0$.  For all except $O_{n,[k:\Q],\epsilon}(Q^{\epsilon})$ fields $F\in\mathscr{F}_k^{n,S_n}(Q)$, there holds
\[
	|\mathrm{Cl}(F)[\ell]|\ll_{n,k,\ell,\epsilon,\eta}D_F^{\frac{1}{2}-\frac{1}{2\ell([F:k]-1)}+\eta}.
\]
\end{enumerate}
\end{theorem}
\begin{remark}
Unlike the work in \cite{Chen_An,EPW,PTW,ThornerZaman}, Theorem \ref{thm:ell-torsion} crucially relies on the fact that the lemma of Ellenberg and Venkatesh exploits non-inert primes, not just primes that split completely.
\end{remark}
\begin{remark}
The lower bounds in \eqref{eqn:lower_bounds} ensure that \cref{thm:ell-torsion} is not vacuous.
\end{remark}

\cref{thm:ell-torsion} is new for degree $p$ fields over a base field $k\neq\Q$, and it is new when $k=\Q$ and $p\geq 7$.  When $p=3$ or $5$, \cref{thm:ell-torsion} greatly reduces the sizes of the exceptional sets in \cite{EPW,PTW}.  \cref{thm:ell-torsion} is new for degree $n$ $S_n$-fields over any given base field $k\neq\Q$ when $n\geq 4$ (see An \cite{2020arXiv200414410A} for $n=3$), and it is new when $k=\Q$ for $n\geq 6$.  When $n\leq 5$, \cref{thm:ell-torsion} greatly reduces the sizes of the exceptional sets in \cite{EPW,PTW}. 

We also prove a mutual refinement of \cref{thm:extremal-class-number} and \cref{thm:ell-torsion} wherein we produce an infinitude of number fields $F$ with a given signature $(r_1,r_2)$ whose Galois closure over $\Q$ has Galois group $S_n$ and whose class group satisfies both \eqref{eqn:EV_GRH} for any fixed integer $\ell\geq 2$ and \eqref{eqn:maximal_class}.  This gives the first examples of number fields of high degree whose class groups have a nontrivial upper bound on the $\ell$-torsion subgroup that provably does not hold for the full class group.  See \cref{thm:extremal-class-number-body} below.  Additionally, when we choose $r_2=0$ so that such $F$ are totally real, we can show that each of the aforementioned fields have a point of exact order $\ell_0$ in their class group, where $\ell_0\geq 2$ is an arbitrary fixed integer.  When $\ell=\ell_0$, this gives the first examples of number fields of high degree with a nontrivial upper bound on the $\ell$-torsion subgroup when the $\ell$-torsion subgroup itself is nontrivial. See \cref{thm:extremal-bilu-luca} below.

\subsection{An effective Chebotarev density theorem for fibers} \label{ssec:fibered-chebotarev}

Let $x\geq 1$.  For a normal extension $K/k$, let $\mathcal{C} \subseteq G \simeq \mathrm{Gal}(K/k)$ be a conjugacy class, and define
\[
	\pi_\mathcal{C}(x; K/k)
		:= \#\{ \mathfrak{p} \subseteq \mathcal{O}_k \text{ prime}\colon  \mathrm{N}_{k/\Q}\mathfrak{p} \leq x, \mathrm{Frob}_\mathfrak{p} \in \mathcal{C}\},
\]
where $\mathcal{O}_k$ denotes the ring of integers of $k$.  The Chebotarev density theorem asserts that
\begin{equation}
\label{eqn:CDT_asymp}
\pi_\mathcal{C}(x; K/k) \sim \frac{|\mathcal{C}|}{|G|} \pi_k(x)\qquad\text{as $x\to\infty$,}
\end{equation}
where $\pi_k(x)$ is the prime ideal counting function of $k$.  When $\zeta_K(s)$ has no Landau--Siegel zero, Thorner and Zaman \cite[Corollary 1.2]{TZ3} proved a stronger result, namely
\[
\pi_{\mathcal{C}}(x; K/k)\sim\frac{|\mathcal{C}|}{|G|}\pi_k(x)\qquad\textup{as }\frac{\log x}{\log([K:\Q]^{[K:\Q]}D_K)}\to\infty.
\]
A similar result holds when a Landau--Siegel zero exists.  This improves previous work of Lagarias and Odlyzko \cite{LO} and V. K. Murty \cite{VKM}.

Let $n\geq 5$ be an integer, let $G$ be the full symmetric group $S_n$, and define
\[
\mathfrak{F}_k^{S_n}(Q):=\{\textup{$K$: $K/k$ is normal, Gal$(K/k)\cong S_n$, $D_K\leq Q$}\}.
\]
Using Theorem \ref{thm:artin-pnt-intro}, we obtain an effective variant of the Chebtarev density theorem that holds in a much wider range for almost all $K \in \mathfrak{F}_{k}^{S_n}(Q)$.  However, since \cref{thm:artin-pnt-intro} imposes the restriction that the representations $\rho$ be {\it faithful}, we do not obtain an equidistribution result in the sense of \eqref{eqn:CDT_asymp}.  Instead, we show that as $\kp$ varies, the conjugacy class of $\mathrm{Frob}_\mathfrak{p}$ is equidistributed in each fiber of the projection $S_n \to S_n/A_n$.

\begin{theorem}\label{thm:fibered-chebotarev-Sn}
	Let  $n \geq 5$ be an integer and $k$ be a number field.  Let $\mathcal{C} \subseteq S_n$ be a conjugacy class.  Given $K\in\mathfrak{F}^{S_n}_k(Q)$, let $k(\sqrt{\Delta_K})/k$ be the unique quadratic extension contained in $K/k$.  Let $Q\geq 1$.  For all $\epsilon>0$, there exists a constant $\Cr{main}=\Cr{main}(|S_n|,[k:\Q],\epsilon)>0$ such that for all except $O_{n,[k:\Q],\epsilon}(Q^{\epsilon})$ fields $K \in \mathfrak{F}^{S_n}_k(Q)$, one has that for any $x \geq (\log D_K)^{81n!/\epsilon}$, there holds
\[
	\pi_\mathcal{C}(x; K/k)
		= \frac{2|\mathcal{C}|}{n!} \pi_{\mathrm{sgn}(\mathcal{C})}(x; k(\sqrt{\Delta_K})/k) + O_{n,[k:\Q],\epsilon}(x \exp(-\Cr{main}\sqrt{\log x})).
\]
\end{theorem}

We highlight two immediate corollaries of Theorem \ref{thm:fibered-chebotarev-Sn}.  First, the Chebotarev density theorem implies that $\pi_{\mathrm{sgn}(\mathcal{C})}(x;k(\sqrt{\Delta_K})/k) \sim  \frac{1}{2}\pi_k(x)$, so Theorem \ref{thm:fibered-chebotarev-Sn} is consistent with \eqref{eqn:CDT_asymp}.  Indeed, by fixing the quadratic subfield $k(\sqrt{\Delta_K})$, we find the following.

\begin{corollary}
\label{cor:CDT1}
	Let $n\geq 5$ be an integer, $k$ be a number field, and $\Delta \in k$ be a non-square element. Let $\epsilon > 0$ and $Q\geq 1$.  For all except $O_{n,[k:\Q],\epsilon}(Q^\epsilon)$ of the fields $K\in\mathfrak{F}_k^{S_n}(Q)$ whose quadratic subfield is $k(\sqrt{\Delta})$, there holds for any $x \geq (\log D_K)^{81n!/\epsilon}$,
	\[
	\pi_\mathcal{C}(x; K/k)
		= \frac{|\mathcal{C}|}{n!} \mathrm{Li}(x) + O_{n,k,\Delta,\epsilon}(x \exp(-\Cr{main}\sqrt{\log x})).
	\]
\end{corollary}
\begin{remark}
For all non-square $\Delta\in k$ and all integers $n\geq 5$, it follows from \cite[Theorem 1.3]{LLOT} that there exist effectively computable constants $\Cl[abcon]{fixed_resolvent}=\Cr{fixed_resolvent}(n,k,\Delta)>0$ and $\Cl[abcon]{fixed_resolvent2}=\Cr{fixed_resolvent2}(n)>0$ such that the number of $K\in\mathfrak{F}_k^{S_n}(Q)$ that contain the quadratic subfield $k(\Delta)/k$ is at least $\Cr{fixed_resolvent}Q^{\Cr{fixed_resolvent2}}$.  This ensures that \cref{cor:CDT1} is not vacuous.
\end{remark}

In the full family $\mathfrak{F}_k^{S_n}(Q)$, where the quadratic resolvent is not assumed fixed, Theorem \ref{thm:fibered-chebotarev-Sn} does not directly permit access to primes whose Frobenius element lies in a single conjugacy class $\mathcal{C}$.  However, in many applications of effective Chebotarev density theorems (e.g., to bounding $\ell$-torsion subgroups of the class group), it is desirable to produce primes whose Frobenius element lies in one of several conjugacy classes.  Theorem \ref{thm:fibered-chebotarev-Sn} provides access to such primes, provided that not all of the desired classes have the same sign.  As a particularly simple instance of this, we have:

\begin{corollary}\label{cor:cdt-opposite-parity}
	Let  $n \geq 5$ be an integer and $k$ be a number field.  Let $\mathcal{C}, \mathcal{C}^\prime \subseteq S_n$ be conjugacy classes of opposite parity.  Let $Q\geq 1$ and $\epsilon>0$.  For all except $O_{n,[k:\Q],\epsilon}(Q^{\epsilon})$ fields $K \in \mathfrak{F}^{S_n}_k(Q)$, one has that for any $x \geq (\log D_K)^{81n!/\epsilon}$, there holds
		\[
			\frac{n!}{|\mathcal{C}|} \pi_\mathcal{C}(x;K/k) + \frac{n!}{|\mathcal{C}^\prime|} \pi_\mathcal{C^\prime}(x;K/k)
				= 2 \pi_k(x) + O_{n,[k:\Q],\epsilon}(x \exp(-\Cr{main}\sqrt{\log x})).
		\]
\end{corollary}

Finally, we note that while we have stated Theorem \ref{thm:fibered-chebotarev-Sn} for the family of $S_n$ extensions, an analogous result will hold for the family of fields $K$ whose Galois group $G$ is a given transitive group of prime degree.  Such groups have a unique minimal normal subgroup $N$, and the role of the quadratic subfield $k(\sqrt{\Delta_K})$ will instead be played by the subfield of $K$ fixed by $N$.  An analogue of Corollary \ref{cor:cdt-opposite-parity} will hold for any set of conjugacy classes that surjects onto the set of conjugacy classes of the quotient $G/N$.

\section*{Organization}

In \cref{sec:approximate-artin}, we state and give context for our main technical results (\cref{thm:zero-density-intro,thm:zfr-transfer,thm:main_result} and \cref{cor:approximate-dedekind}) after summarizing the work in \cite{PTW,ThornerZaman}.

In \cref{sec:artin-prelim}, we recall the definition and basic properties of Artin $L$-functions.

In \cref{sec:transfer}, we prove Theorem \ref{thm:zfr-transfer}, which ensures that a zero-free region for the quotient $\zeta_K(s)/\zeta_{K^N}(s)$ ``transfers'' to Artin $L$-functions not coming from $K^N$.

\cref{sec:zero-density} and \ref{sec:zero-free} contain the main analytic results of this paper, including the proofs of \cref{thm:zero-density-intro,thm:main_result}.

In \cref{sec:degree-p}, we discuss how our results apply to the family of degree $p$ extensions, regardless of Galois structure.

In \cref{sec:chebotarev-proof}, we prove the fibered Chebotarev density theorem in \cref{thm:fibered-chebotarev}.

In \cref{sec:class-groups}, we prove the applications to class groups in \cref{thm:extremal-class-number,thm:ell-torsion}.

In \cref{sec:subconvexity}, we prove the applications to subconvexity and the equidistribution of periodic torus orbits in \cref{thm:subconvexity-nice-families,thm:torus-orbits}.

In \cref{ssec:multiplicity}, we give heuristics for a key quantity called the intersection multiplicity introduced below and use this discussion to contrast our results with those of \cite{PTW} and \cite{ThornerZaman}.

\section*{Acknowledgements}

The authors would like to thank William Duke, John Friedlander, Boris Hasselblatt, Peter Humphries, George McNinch, Philippe Michel, Ken Ono, David Rohrlich, Kannan Soundararajan, Frank Thorne, Akshay Venkatesh, and Jiuya Wang for helpful conversations. RJLO was partially supported by NSF grant DMS-1601398.

\section{Main results: holomorphy and non-vanishing of Artin $L$-functions} \label{sec:approximate-artin}

We now describe the ideas leading to \cref{thm:artin-pnt-intro}.  As above, let $k$ be a number field and let $G$ be a finite group.  We let $\mathfrak{F}_{k}^{G}$ denote the family of number fields $K$ inside a fixed choice of the algebraic closure $\bar k$ that are normal extensions of $k$ with Galois group $\Gal(K/k)$ isomorphic to $G$.  For any $Q \geq 1$, let
\[
	\mathfrak{F}_{k}^{G}(Q) := \{ K\in\mathfrak{F}_k^G\colon D_K \leq Q\},
\]
where $D_K$ denotes the absolute discriminant of $K$ over $\mathbb{Q}$.

\subsection{Summary of preceding work}

It is instructive to briefly review the ideas in \cite{PTW,ThornerZaman}.  Let $k=\Q$.  The approach of Pierce, Turnage-Butterbaugh, and Wood in \cite{PTW} relied on the zero density estimate for $L$-functions of families of cuspidal automorphic representations proved by Kowalski and Michel \cite{KM}.  Assuming the strong Artin conjecture, the $L$-function $\zeta_K(s)/\zeta_{\Q}(s)$ associated to each $K\in\mathfrak{F}_{\Q}^G$ is  the $L$-function associated to an isobaric {\it non-cuspidal} automorphic representation $\Pi_K$ defined over $\Q$.  Note that if $K_1,K_2\in\mathfrak{F}_{\Q}^G(Q)$ are distinct, then $\Pi_{K_1}$ and $\Pi_{K_2}$ might have some cuspidal constituents in common (as would happen if $K_1$ and $K_2$ share a common subfield), in which case there exists a cuspidal automorphic representation $\pi_0$ such that $L(s,\Pi_{K_1})/L(s,\pi_0)$ and $L(s,\Pi_{K_2})/L(s,\pi_0)$ are entire. 
If one studies the zeros of the $L$-functions $L(s,\Pi_K)$ with $K\in\mathfrak{F}_{\Q}^G(Q)$ and many of the $L(s,\Pi_K)$ share a particular common factor $L(s,\pi_0)$, then the zeros of $L(s,\pi_0)$ are counted with high multiplicity.  If this multiplicity is too high, then the zero density estimate is rendered trivial.

Let $\mathscr{R}_G$ denote a particular condition on the primes that tamely ramify in a normal extension $K/\Q$ with $\Gal(K/\Q)\cong G$, and let $\mathfrak{F}_{\Q}^G(Q,\mathscr{R}_G)$ be the subset of $\mathfrak{F}_{\Q}^G(Q)$ whose fields satisfy $\mathscr{R}_G$.  For certain groups $G$, Pierce, Turnage-Butterbaugh, and Wood find conditions $\mathscr{R}_G$ that enable them to relate the distribution of fields $K_1,K_2\in\mathfrak{F}_{\Q}^G(Q,\mathscr{R}_G)$ such that $\Pi_{K_1}$ and $\Pi_{K_2}$ share some cuspidal constituents to the arithmetic-statistical problem of counting number fields $K\in\mathfrak{F}_{\Q}^G(Q,\mathscr{R}_G)$ that share a given discriminant.  In the situations where their approach works (see Section 6.3 and Theorems 3.1 and 3.3 in \cite{PTW}), they prove that if one assumes the strong Artin conjecture for $G$, then for all $\epsilon>0$ and all $K\in\mathfrak{F}_{\Q}^G(Q,\mathscr{R}_G)$ with at most
\[
O_{|G|,\epsilon}(Q^{\epsilon}\max_{D\leq Q}|\{K\in\mathfrak{F}_{\Q}^G(Q,\mathscr{R}_G)\colon D_K=D\}|)
\]
exceptions, the ratio $\zeta_K(s)/\zeta_{\Q}(s)$ (and all of the Artin $L$-functions in its factorization, which are assumed to be automorphic, hence entire) is nonvanishing in the region 
\[
\re(s)\geq 1-\epsilon,\qquad |\im(s)|\leq(\log D_K)^{2/\epsilon}.
\]
A ``discriminant multiplicity conjecture'' of Duke \cite{Duke98} implies that
\begin{equation}
\label{eqn:disc_mult}
\max_{D\leq Q}|\{K\in\mathfrak{F}_{\Q}^G(Q,\mathscr{R}_G)\colon D_K=D\}|\ll_{G,\epsilon} Q^{\epsilon}.
\end{equation}
Since the strong Artin conjecture is assumed in \cite{PTW}, this would imply that for all except $O_{|G|,\epsilon}(Q^{\epsilon})$ of the $K\in\mathfrak{F}_{\Q}^G(Q,\mathscr{R}_G)$, the $L$-functions of the Artin representations attached to $\Gal(K/\Q)$ have a strong zero-free region for the low-lying zeros, provided that there exists a constant $\tau=\tau(G)>0$ such that $|\mathfrak{F}_{\Q}^G(Q,\mathscr{R}_G)|\gg_{|G|}Q^{\tau}$.  Much of the work in \cite{PTW} centers around making progress towards \eqref{eqn:disc_mult} for groups $G$ for which the strong Artin conjecture holds and for which $|\mathfrak{F}_{\Q}^G(Q,\mathscr{R}_G)|\gg_{|G|}Q^{\tau}$ with $\tau$ suitably large, including cyclic groups, dihedral groups of order $2p$ for odd primes $p$, $S_3$, and $S_4$.

Note, however, that if the ramification restriction $\mathscr{R}_G$ is non-empty then the results of Pierce, Turnage-Butterbaugh, and Wood do not quantify the number of exceptional fields in the full family $\mathfrak{F}_{\Q}^G(Q)$, even assuming the full force of Duke's discriminant multiplicity conjecture \eqref{eqn:disc_mult} and the strong Artin conjecture.
See \cref{ssec:multiplicity} for a discussion of the limitations.

Let $\rho_K$ be the Artin representation of $\Gal(K/\Q)$ such that $L(s,\rho_K)=\zeta_K(s)/\zeta_{\Q}(s)$.  The approach of Thorner and Zaman \cite{ThornerZaman} removes the need to assume the strong Artin conjecture by proving the first unconditional large sieve for the Artin representations $\rho_K$ as $K\in\mathfrak{F}_{\Q}^G(Q)$ varies.  They used character theory for the tensor products $\rho_{K_1}\otimes\rho_{K_2}$ and Galois theory in lieu of automorphy.  In the process, they simplified the arithmetic-statistical problem that one must solve to address the subfield problem.  Defining the ``intersection multiplicity''
\begin{equation}
\label{eqn:intersection_multiplicity}
\mathfrak{m}_{k}^G(Q):=\max_{K_1\in\mathfrak{F}_{k}^G(Q)}|\{K_2\in\mathfrak{F}_{k}^G(Q)\colon K_1\cap K_2\neq k\}|,
\end{equation}
they unconditionally proved that for all $K\in\mathfrak{F}_{\Q}^G(Q)$ with at most $O_{G,\epsilon}(\mathfrak{m}_{\Q}^G(Q)Q^{\epsilon})$ exceptions, the ratio $\zeta_K(s)/\zeta_{\Q}(s)$ is nonvanishing in a region containing the box
\[
1-\frac{\epsilon}{10^8|G|^3}\leq \re(s)\leq 1,\qquad |\im(s)|\leq D_K^{1000}.
\]
This result is nontrivial if there exists a constant $\delta>0$ such that
\[
\mathfrak{m}_{\Q}^G(Q)\ll_{G,\delta}Q^{-\delta}|\mathfrak{F}_{\Q}^G(Q)|.
\]
Since two normal extensions meet in a normal extension, if $G$ is simple, then $\mathfrak{m}_{\Q}^G(Q)=1$.  Otherwise, the best bounds on $\mathfrak{m}_{\Q}^G(Q)$ follow from progress toward \eqref{eqn:disc_mult} and typically also require restrictions $\mathscr{R}_G$ on ramification.  Therefore, the need for the strong Artin conjecture is removed, but the subfield problem still remains unaddressed apart from a handful of special cases.  If a suitable bound for $\mathfrak{m}_{\Q}^G(Q)$ is known, then for all $K\in\mathfrak{F}_{\Q}^G(Q)$ with few exceptions, the large zero-free region of $\zeta_K(s)/\zeta_{\Q}(s)$ will translate to a large region of holomorphy and non-vanishing for all of $L$-functions associated to the nontrivial Artin representations of $\Gal(K/\Q)$.  This last step crucially uses Artin induction to express $L(s,\rho)$ in terms of $L$-functions of one-dimensional representations of cyclic subgroups of $G$, each of which inherits the zero-free region of $\zeta_K(s)/\zeta_{\Q}(s)$.

For an arbitrary group $G$, and in particular those for which Malle's conjecture is not known, it seems quite difficult to prove that there exists a constant $\delta=\delta(|G|,k)>0$ such that $\mathfrak{m}_{k}^G(Q)\ll_{|G|,k}Q^{-\delta}|\mathfrak{F}_{k}^G(Q)|$, regardless of whether $k=\Q$.  (Recall that Malle's conjecture predicts an asymptotic formula for the growth of $|\mathfrak{F}_k^G(Q)|$ as $Q \to \infty$, and this is known only in few cases.)  This is the technical limitation that we discussed in \cref{sec:intro}.  

If $G$ is simple, then this issue disappears as $\mathfrak{m}_k^G(Q) = 1$, but if $G$ is not simple, then one should expect $\mathfrak{m}_k^G(Q) \gg_{G,k} Q^{c}$ for some constant $c$ depending on $G$ (see Conjecture \ref{conj:mgn-pi} below).  This is the strucutral limitation discussed in \cref{sec:intro}.  For example, for the family $\mathfrak{F}_k^{S_n}(Q)$, it follows from \cite[Theorem 1.3]{LLOT} that
\begin{equation}
	\label{eqn:structural}
\mathfrak{m}_{k}^{S_n}(Q)\gg_{n,k} Q^{(\frac{1}{8}-\frac{27}{32n})\frac{1}{n!}}.	
\end{equation}
In fact, it follows from Conjecture \ref{conj:mgn-pi} below that we should expect that for all $\epsilon>0$, we have $\mathfrak{m}_k^{S_n}(Q) \gg_{k,n,\epsilon} Q^{-\epsilon}|\mathfrak{F}_k^{S_n}(Q)|$.  Thus, the desired bound $\mathfrak{m}_k^{G}(Q) \ll_{G,k,\epsilon}Q^{-\delta}|\mathfrak{F}_{k}^G(Q)|$ is not expected to hold when $G = S_n$, nor is it expected to hold in many other natural situations.

\subsection{Changing the average}
\label{subsec:component_1}

As we mentioned earlier, our new approach that circumvents the technical and structural limitations of the approach in \cite{ThornerZaman} has two independent components.  To describe them, let $N \trianglelefteq G$ be a nontrivial normal subgroup and, for any $K \in \mathfrak{F}_{k}^{G}$, let $K^N$ denote the subfield of $K$ fixed by $N$ under the given isomorphism $\mathrm{Gal}(K/k)\simeq G$.  Note that $N=G$ is permissible, in which case $K^N=k$.  In our current setting, we no longer require our base field $k$ to equal $\Q$.

The first component of our new method is to study the $L$-function $\zeta_K(s)/\zeta_{K^N}(s)$ as $K\in\mathfrak{F}_k^G(Q)$ varies instead of the $L$-functions $\zeta_K(s)/\zeta_k(s)$.
The quotient $\zeta_{K}(s)/\zeta_{K^N}(s)$ is entire by the Aramata--Brauer theorem.  With suitable modifications to the ideas in \cite{ThornerZaman}, the work of Brauer in \cite{Brauer} and Galois theory will once again alleviate the need for unproven analytic hypotheses such as the Artin conjecture.  The crux of our new average is that we trade the intersection multiplicity $\mathfrak{m}_k^G(Q)$ in \eqref{eqn:intersection_multiplicity} that arises in \cite{ThornerZaman} for a new multiplicity, namely
\begin{equation}
\label{eqn:normalized_intersection_multiplicity}
\mathfrak{m}_k^{G,N}(Q):=\max_{K_1\in\mathfrak{F}_k^G(Q)}|\{K_2\in\mathfrak{F}_k^G(Q)\colon K_1\cap K_2\neq K_1^N\cap K_2^N\}|.
\end{equation}
The first component of our new approach is summarized in the following theorem.
\begin{theorem}
\label{thm:zero-density-intro}
Let $Q\geq 1$.  Let $G$ be a finite group, $N\trianglelefteq G$ be a nontrivial normal subgroup, and $k$ be a number field.  Let $\mathfrak{m}_k^{G,N}(Q)$ be as in \eqref{eqn:normalized_intersection_multiplicity}.  There exists an absolute and effectively computable constant $\Cl[abcon]{ZFR}>0$ such that for all $\epsilon>0$ and all number fields $K\in\mathfrak{F}_k^G(Q)$ with $O_{|G|,[k:\Q],\epsilon}(\mathfrak{m}_{k}^{G,N}(Q) Q^\epsilon)$ exceptions, the quotient $\zeta_K(s)/\zeta_{K^N}(s)$ is non-vanishing in the region $\Omega_K(\epsilon)$ defined by
\begin{equation}
\label{eqn:nice_ZFR}
1-\re(s)\leq \begin{cases}
\displaystyle\frac{\epsilon}{10|G|}\frac{ \log D_K}{\log D_K+[k:\Q]\log(3+|\im(s)|)}&\mbox{if $|\im(s)|\leq \exp(D_K^{\epsilon/(6C_{G}[k:\Q])})$,}\\ \\
\displaystyle\frac{\Cr{ZFR}}{\log D_K+|G|[k:\Q] \log(3+|\im(s)|)}&\mbox{if $|\im(s)|>\exp(D_K^{\epsilon/(6C_{G}[k:\Q])})$.}
\end{cases}
\end{equation}
The constant $C_G$, which depends at most on $|G|$, is the same as in \cref{thm:ZDE} below.
\end{theorem}
\begin{remark}
It follows from work of Lagarias and Odlyzko \cite[Section 8]{LO} that $\zeta_K(s)/\zeta_{K^N}(s)$ does not vanish in the region
\begin{equation}
	\label{eqn:Std_ZFR-intro}
	1-\re(s)\leq \frac{\Cr{ZFR}}{\log D_K+|G|[k:\Q]\log(3+|\im(s)|)}
\end{equation}
apart from at most one exceptional zero of $\zeta_K(s)$, which (if it exists) is necessarily real and simple.  By comparison, the zero-free region $\Omega_K(\epsilon)$ defined by \eqref{eqn:nice_ZFR} contains the box
\[
1-\frac{\epsilon}{20|G|}\leq \re(s)\leq 1,\qquad |\im(s)|\leq D_K^{1000/[k:\Q]}
\]
when $D_K$ is sufficiently large with respect to $\epsilon$ and $|G|$.  Thus, the zero-free region $\Omega_K(\epsilon)$ constitutes a substantial improvement over \eqref{eqn:Std_ZFR-intro}, and it applies for the vast majority of $K\in\mathfrak{F}_k^G(Q)$ when we can prove that $\mathfrak{m}_k^{G,N}(Q)$ is small.
\end{remark}

Note that $\mathfrak{m}_k^{G,G}(Q)=\mathfrak{m}_k^{G}(Q)$.  Therefore, our new results subsume all of the work in \cite{ThornerZaman}.  The novelty of \cref{thm:zero-density-intro} is that when $N\neq G$, the quantity $\mathfrak{m}_{k}^{G,N}(Q)$ can often be controlled independently of one's ability to estimate the size of the family $\mathfrak{F}_{k}^{G}(Q)$.  Notably, if $G$ has a {\it unique minimal nontrivial normal subgroup} $N$, then $\mathfrak{m}_{k}^{G,N}(Q) = 1$ for all $Q\geq 1$.  This is the case for many groups of interest, including:

\begin{itemize}
	\item all transitive permutation groups of prime degree, including the dihedral groups $D_p$;
	\item the affine linear group $\mathrm{AGL}_d(\mathbb{F}_p)$ for any integer $d \geq 1$ and any prime $p$;
	\item the full symmetric group $S_n$ for any $n \geq 2$; and
	\item the non-simple alternating group $A_4$.
\end{itemize}

In fact, generalizing these examples, any primitive permutation group has at most $2$ minimal normal subgroups \cite[Theorem 4.3B]{DixonMortimer}.  The determination of which of these has a unique minimal normal is a consequence of the O'Nan--Scott theorem \cite[Theorem 4.1A]{DixonMortimer}.  As a consequence of this, for example, if $n$ is not equal to $|T|^k$ for some nonabelian simple group $T$ and integer $k \geq 1$, then every primitive group of degree $n$ has a unique minimal normal subgroup. Many imprimitive groups have unique minimal normal subgroups as well, corresponding to subgroups stabilizing the nontrivial blocks, but these are not our main focus.

\subsection{Character theory with restricted components}
\label{subsec:component_2}

Let $K \in \mathfrak{F}_k^G$.  Suppose that $\zeta_K(s) / \zeta_{K^N}(s)$ is non-vanishing in some region $\Omega \subseteq \mathbb{C}$, e.g., the region $\Omega_K(\epsilon)$ in \cref{thm:zero-density-intro}.  In light of the factorization
	\begin{equation} \label{eqn:K-over-KN-factorization}
		\frac{\zeta_K(s)}{\zeta_{K^N}(s)} 
			= \prod_{\substack{ \rho \in \mathrm{Irr}(G) \\ \ker \rho \not\supseteq N}} L(s,\rho)^{\dim \rho},
	\end{equation}
it is reasonable to hope that the zero-free region $\Omega$ extends to each of the irreducible Artin $L$-functions $L(s,\rho)$ whose associated representation has kernel not containing $N$.  However, since it is not known that the $L$-functions $L(s,\rho)$ appearing in this factorization are holomorphic, this does not immediately follow, nor does it follow in general from any existing result on Artin $L$-functions.  The second component of our approach, therefore, is a new conjecture in the character theory of finite groups that would accommodate such a transfer of zero-free regions.  We prove this conjecture in many cases of interest (including for the symmetric group and groups of prime degree), and we provide a reduction that suggests a general attack based on the classification of finite simple groups.

To describe this conjecture, we begin by recalling classical work of Artin.  In particular, for each representation $\rho$ of a finite group $G$, none of whose irreducible consituents is trivial, Artin showed that there are 
 rational constants $c_{\rho,\chi}$ such that
	\begin{equation} \label{eqn:artin-induction-1}
		\mathop{\mathrm{tr}} \rho 
			= \sum_{H} \sum_{ \substack{ \chi \in \mathrm{Irr}(H): \\ \dim \chi = 1}} c_{\rho,\chi} \mathrm{Ind}_H^G \chi,
	\end{equation}
where the summation over $H$ runs over the cyclic subgroups of $G$.  This leads to a corresponding factorization of the Artin $L$-function in terms of Hecke $L$-functions,
	\begin{equation} \label{eqn:artin-factorization-1}
		L(s,\rho)
			= \prod_H \prod_{ \substack{ \chi \in \mathrm{Irr}(H): \\ \dim \chi = 1}} L(s,\chi)^{c_{\rho,\chi}},
	\end{equation}
from which Artin deduced that some integral power of $L(s,\rho)$ possesses meromorphic continuation to all of $\mathbb{C}$.  This work was later extended by Brauer to show that each $L(s,\rho)$ itself is meromorphic by using a different class of subgroups (namely, so called ``elementary'' subgroups) but the decomposition of $L(s,\rho)$ into Hecke $L$-functions, or equivalently of $\mathop{\mathrm{tr}}\rho$ into the induction of $1$-dimensional characters, remains essentially the only general way of inferring analytic properties of Artin $L$-functions outside the region of absolute convergence.

To infer consequences toward non-vanishing, we note that in the factorization of the quotient $\zeta_K(s)/\zeta_k(s)$ in \eqref{eqn:artin-factorization-1}, every $L(s,\chi)$ associated to a cyclic subgroup can be taken to appear with a positive exponent.  It follows that if $\zeta_K(s)/\zeta_k(s)$ is non-vanishing in a region $\Omega$, then every $L(s,\chi)$ is non-vanishing in $\Omega$ as well.  Consequently, every $L(s,\rho)$ is both holomorphic and non-vanishing in $\Omega$ by means of the factorization \eqref{eqn:artin-factorization-1}.  For our purposes, however, we wish to assume only that $\zeta_K(s)/\zeta_{K^N}(s)$ is non-vanishing in some region, and infer the holomorphy and non-vanishing of every irreducible Artin $L$-function $L(s,\rho)$ whose kernel does not contain $N$.  From \eqref{eqn:artin-factorization-1}, this only follows for those irreducible representations that are induced from $N$, which is typically a small subset of the irreducible representations at hand.  

We propose the following hypothesis on the finite group $G$ and normal subgroup $N \unlhd G$.

\begin{namedtheorem}[\hypertarget{hyp:T}{Hypothesis $\mathrm{T}(G,N)$}] 
	If $\rho$ is an irreducible representation of $G$ such that $N\not\subseteq\ker\rho$, then for each subgroup $H \subseteq G$ and each one-dimensional character $\chi$ of $H$ for which $H \cap N \not\subseteq \ker \chi$, there exists $c_{\rho,\chi}\in\Q$ such that
	\[
	\mathop{\mathrm{tr}}\rho = \sum_{H\subseteq G}\sum_{\substack{\chi\in\mathrm{Irr}(H) \\ \dim\chi=1 \\ H\cap N\not\subseteq\ker\chi}}c_{\rho,\chi}\mathrm{Ind}_H^G\chi.
	\]	
\end{namedtheorem}

\begin{conjecture}
\label{conj:restricted-induction}
For all finite groups $G$ and normal subgroups $N \unlhd G$, Hypothesis  \hyperlink{hyp:T}{$\mathrm{T}(G,N)$} holds.
\end{conjecture}

Two remarks are in order.  First, Conjecture \ref{conj:restricted-induction} lies much deeper than \eqref{eqn:artin-induction-1}, the proof of which is almost immediate from a modern perspective.  In particular, \cref{conj:restricted-induction} is not amenable to standard techniques exploiting the adjointness of the induction and restriction maps.  For this reason, it appears to be much more group theoretic in nature than \eqref{eqn:artin-induction-1}.  Additionally, the analogue of Conjecture \ref{conj:restricted-induction} does not hold if the subgroups $H$ are required to be cyclic or even abelian.  For this reason, a wider set of subgroups is required, analogous to how Brauer enlarged the set of subgroups to obtain a version of \eqref{eqn:artin-induction-1} with integral coefficients.  We note that most, or perhaps all, proofs of Brauer induction rely on the module structure of the character ring of $G$ and proceed by finding a representation of the trivial character in terms of inductions of characters from elementary subgroups.  However, such an approach cannot work for Conjecture \ref{conj:restricted-induction}---there is less inherent module structure at play and the trivial character is not in the subspace under consideration.

Second, the $1$-dimensional characters appearing in Conjecture \ref{conj:restricted-induction}, namely those $\chi$ such that $\ker \chi \not\supseteq H \cap N$, are precisely those whose induction to $G$ may be decomposed solely in terms of the characters of irreducible representations $\rho$ of $G$ whose kernel does not contain $N$.  Thus, this is the largest set of $1$-dimensional characters to which a zero-free region of $\zeta_K(s)/\zeta_{K^N}(s)$ might plausibly transfer.  Indeed, we show in Lemma \ref{lem:abelian-transfer} that such a transfer always occurs; thus, the importance of Conjecture \ref{conj:restricted-induction} is made clear by the following result.

\begin{theorem}\label{thm:zfr-transfer}
	Let $G$ be a finite group, $N\trianglelefteq G$ be a nontrivial normal subgroup for which Hypothesis \hyperlink{hyp:T}{$\mathrm{T}(G,N)$} holds, and $k$ be a number field.  For any $K \in \mathfrak{F}_{k}^{G}$, if $\zeta_K(s)/\zeta_{K^N}(s)$ is non-zero in a region $\Omega \subseteq \mathbb{C}$, then for each irreducible Artin representation $\rho$ of $\mathrm{Gal}(K/k)$ whose kernel does not contain $N$, the Artin $L$-function $L(s,\rho)$ is holomorphic and non-vanishing on $\Omega$.  
\end{theorem}

\begin{remark}
It is worthwhile to compare Theorem \ref{thm:zfr-transfer} with a well known theorem of Stark \cite[Theorem 3]{Stark}.  Let $K/k$ be a normal extension of number fields.  Stark shows that if $\zeta_K(s)$ has a simple zero, then this zero must be inherited from the Dedekind zeta function of a cyclic extension of $k$, and that it is not a zero or pole of any Artin $L$-function that does not factor through this cyclic extension.  (If this simple zero is also real, then the cyclic extension must in fact be at most a quadratic extension.  This is how Stark's theorem is most commonly invoked.)  Since cyclic extensions of $k$ contained in $K$ correspond to normal subgroups $N$ for which $G/N$ is cyclic, Stark's theorem may be interpreted as showing that simple zeros of $\zeta_K(s)$ are constrained to arise from zeros of Dedekind zeta functions $\zeta_{K^N}(s)$ for cyclic extensions $K^N/k$, and that these zeros do not propagate to Artin $L$-functions attached to representations whose kernel does not contain $N$.  By contrast, a simple consequence of Theorem \ref{thm:zfr-transfer} is that if $s_0$ is a zero of $\zeta_{K}(s)$ of any order that is ``explained'' by $\zeta_{K^N}(s)$ in the sense that $\zeta_K(s)/\zeta_{K^N}(s)$ is analytic and non-vanishing at $s_0$, then this zero does not propagate, i.e. every Artin $L$-function that does not factor through $K^N$ must be analytic and non-vanishing at $s_0$.
\end{remark}

We now record our progress toward Conjecture \ref{conj:restricted-induction}.

\begin{theorem}
	\label{thm:conj}
	Let $G$ be a finite group and let $N \unlhd G$ be a nontrivial normal subgroup.  Then Hypothesis \hyperlink{hyp:T}{$\mathrm{T}(G,N)$} holds if the index $[G:N]$ is a prime power or if $N$ is solvable (and thus also if $G$ is solvable).
	Additionally, it holds if the order of $G$ is at most $2000$ or if $G$ is a transitive permutation group of degree at most $31$.
\end{theorem}
\begin{proof}
	When $N$ is solvable or the index $[G:N]$ is a prime power, this follows from \cref{thm:orthogonal-span} below.  The remaining claims follow from a computation in \verb^Magma^.
\end{proof}

The following corollary to \cref{thm:conj} is crucial for the applications in \cref{sec:Applications}.

\begin{corollary}
\label{cor:approx_Sn_Sp}
Let $n\geq 2$ be an integer and $p$ be prime.  Then Hypothesis \hyperlink{hyp:T}{$\mathrm{T}(G,N)$} holds when $G$ is the symmetric group $S_n$ or a transitive group of degree $p$, and $N$ is any normal subgroup of $G$.
\end{corollary}
\begin{proof}
	If $n \neq 4$, the only nontrivial normal subgroups of $S_n$ are $S_n$ itself and the alternating group $A_n$.  Both of these have prime power index.  If $n=4$, there is also the Klein four subgroup, which is abelian (hence solvable).  If $p$ is prime, then the transitive groups of degree $p$ have been classified; see Lemma \ref{lem:transitive-p} below.  It follows from this classification that either $G$ is solvable or the unique minimal normal subgroup of $G$ has prime power index, in which case every nontrivial normal subgroup will also have prime power index.
\end{proof}

Finally, we have succeeded in reducing the general conjecture, \cref{conj:restricted-induction}, to the case that $N$ is a minimal normal subgroup and $G \subseteq \mathrm{Aut}(N)$.  Exploiting the characterization of minimal normal subgroups, along with the techniques used in Theorem \ref{thm:conj}, gives the following.

\begin{theorem} \label{thm:minimal-normal}
	Let $T$ be a nonabelian simple group and let $p$ be a prime.  If Hypothesis \hyperlink{hyp:T}{$\mathrm{T}(G,N)$} holds for all groups $G \subseteq \mathrm{Aut}(T) \wr C_p$ containing $N = T^p$ for which the quotient $G/N$ is cyclic, then it holds for all finite groups.
\end{theorem}

\subsection{Holomorphy and non-vanishing}
\label{subsec:holomorphy_non-vanishing}

Our main result follows from the combination of \cref{thm:zero-density-intro,thm:zfr-transfer}.  In what follows, we define $\chi_\rho(\mathfrak{p}) := \mathop{\mathrm{tr}}\rho(\mathrm{Frob}_\mathfrak{p})$.

\begin{theorem}
\label{thm:main_result}
Let $G$ be a finite group, and let $N\trianglelefteq G$ be a nontrivial normal subgroup such that Hypothesis \hyperlink{hyp:T}{$\mathrm{T}(G,N)$} holds.  Let $Q\geq 1$, let $k$ be a number field, and recall  $\mathfrak{m}_k^{G,N}(Q)$ is defined by \eqref{eqn:normalized_intersection_multiplicity}.  For all $\epsilon>0$, there exists an effectively computable constant $\Cr{main}=\Cr{main}(|G|,[k:\Q],\epsilon)>0$ such that for all $K \in \mathfrak{F}_{k}^{G}(Q)$ apart from at most $O_{|G|,[k:\Q],\epsilon}(\mathfrak{m}_{k}^{G,N}(Q) Q^\epsilon)$ exceptions, the following properties hold for the 
Artin representations $\rho$ of $\Gal(K/k)$ whose kernel does not contain $N$:
\begin{enumerate}
	\item $L(s,\rho)$ is holomorphic and non-vanishing in the region $\Omega_K(\epsilon)$ defined by \eqref{eqn:nice_ZFR}, and
	\item if $x\geq (\log D_K)^{81|G|/\epsilon}$, then
	\begin{equation}
	\label{eqn:pnt_a.e.}
	\Big|\sum_{\N_{k/\Q}\kp\leq x}\chi_{\rho}(\kp)\Big|\ll_{|G|,[k:\Q],\epsilon} x\exp(-\Cr{main}\sqrt{\log x}).
	\end{equation}
\end{enumerate}
If $N$ is the unique minimal nontrivial normal subgroup of $G$, then $\mathfrak{m}_{k}^{G,N}(Q)=1$.
\end{theorem}
\begin{remark}
One of the primary benefits of a strong zero-free region for an $L$-function is that one can typically prove a correspondingly strong analogue of the prime number theorem, provided that one can suitably bound the logarithmic derivative.  This is usually done by exploiting the full analytic continuation of the $L$-function, but this is not something afforded by Theorem \ref{thm:zfr-transfer}.  However, using the central ideas of its proof, we are still able to prove the strong effective prime number theorem \eqref{eqn:pnt_a.e.} for the $L(s,\rho)$ considered in \cref{thm:zfr-transfer}.  This is why we list the region of holomorphy and non-vanishing separately from the effective prime number theorem in \cref{thm:main_result}.
\end{remark}

\begin{proof}[Proof of \cref{thm:artin-pnt-intro}]
Let $G$ be $S_n$ for some integer $n\geq 2$ or a transitive subgroup of $S_p$ for some prime $p$.  By \cref{cor:approx_Sn_Sp}, $G$ satisfies \cref{conj:restricted-induction}.  By \cref{lem:degree-p}, $G$ has a unique nontrivial minimal normal subgroup, say $N(G)$.  Since $N(G)$ is nontrivial, it cannot be contained in the kernel of a faithful representation $\rho$ of $G$ since $\ker\rho$ is trivial.  Therefore, since $\mathfrak{m}_{k}^{G,N(G)}(Q)=1$, \cref{thm:artin-pnt-intro} follows from \cref{thm:main_result}.
\end{proof}

Let $G$ be a finite group, and let $N\unlhd G$ be a nontrivial normal subgroup. We observe that if $F$ is a subfield of $K\in\mathfrak{F}_k^G(Q)$ for which $F\cap K^N=k$, then the kernel of the Artin representation $\rho_F$ of $\Gal(K/k)$ associated to the Artin $L$-function $\zeta_F(s)/\zeta_k(s)$ does not contain $N$.  Therefore, \cref{thm:main_result} is applicable.  We emphasize two widely applicable cases of this in our next result; these will enable us to prove the applications in \cref{sec:Applications}. 

\begin{corollary} \label{cor:approximate-dedekind} 
Let $k$ be a number field and $Q\geq 1$.  For a field $F$, let $\tilde{F}$ denote its normal closure over $k$ and let $\rho$ be the Artin representation satisfying $L(s,\rho_F) = \zeta_F(s)/\zeta_k(s)$. 
\begin{enumerate}
	\item Let $n\geq 3$.  For all $\epsilon>0$, there exists an effectively computable constant $\Cl[abcon]{approxDedekind_n}=\Cr{approxDedekind_n}(n,[k:\Q],\epsilon)>0$ such that for all except $O_{n,[k:\Q],\epsilon}(Q^\epsilon)$ fields $F \in \mathscr{F}_{k}^{n,S_n}(Q)$,
	\begin{enumerate}
		\item  $L(s,\rho_F)$ is holomorphic and non-vanishing in the region $\Omega_{\tilde{F}}(\epsilon/n!)$, and 
		\item  if $x\geq (\log D_F)^{81(n!)^2/\epsilon}$, then
		\begin{equation}
		\label{eqn:partialsums_n}
		\Big|\sum_{\N_{k/\Q}\kp\leq x}\chi_{\rho_F}(\kp)\Big|\ll_{n,[k:\Q],\epsilon}x\exp(-\Cr{approxDedekind_n}\sqrt{\log x}).
		\end{equation}
	\end{enumerate}
	\item Let $p$ be prime.  For all $\epsilon>0$, there exists an effectively computable constant $\Cl[abcon]{approxDedekind_p}=\Cr{approxDedekind_p}(p,[k:\Q],\epsilon)>0$ such that for all except $O_{p,[k:\Q],\epsilon}(Q^\epsilon)$ fields $F \in \mathscr{F}_{k}^{p}(Q)$,
	\begin{enumerate}
		\item  $L(s,\rho_F)$ is holomorphic and non-vanishing in the region $\Omega_{\tilde{F}}(\epsilon/p!)$, and
		\item  if $x\geq (\log D_F)^{81(p!)^2/\epsilon}$, then
		\begin{equation}
		\label{eqn:partialsums_p}
		\Big|\sum_{\N_{k/\Q}\kp\leq x}\chi_{\rho_F}(\kp)\Big|\ll_{p,[k:\Q],\epsilon}x\exp(-\Cr{approxDedekind_p}\sqrt{\log x}).
		\end{equation}
	\end{enumerate}

\end{enumerate}
\end{corollary}

\begin{proof}~
(1) When $G=S_n$, let $N$ be the unique minimal nontrivial normal subgroup of $G$, which is either $A_n$ or $V_4$.  In either case, $N$ is transitive, and thus its interesection with a stabilizer subgroup of $S_n$ has index $n$ in $N$.  Additionally, since $N$ is the unique minimal, $\mathfrak{m}_k^{G,N}(R)=1$ for all  $R\geq 1$, and Hypothesis \hyperlink{hyp:T}{$\mathrm{T}(G,N)$} is satisfied via \cref{cor:approx_Sn_Sp}.  Finally, there exists an effectively computable constant $\Cl[abcon]{disc_c}=\Cr{disc_c}(n,[k:\Q])>0$ such that $D_{\tilde{F}}\leq \Cr{disc_c}D_F^{[\tilde{F}:k]}\leq \Cr{disc_c}D_F^{[F:k]!}$.  Thus, for each $F\in\mathscr{F}_k^{n,S_n}(Q)$, the normal closure $\tilde{F}$ over $k$ lies in $\mathfrak{F}_k^{S_n}(\Cr{disc_c}Q^{[F:k]!})$.  If $F\in \mathscr{F}_k^{n,S_n}(Q)$, then $F\cap \tilde{F}^N =k$ because $F$ is the fixed field of a stabilizer subgroup.  The result now follows from \cref{thm:main_result}.

(2) The proof is the same as the previous part, except that we combine the contributions from all of the transitive subgroups of $S_p$ (of which there are $O_p(1)$).  We invoke \cref{lem:degree-p} and \cref{cor:approx_Sn_Sp} to each of these transitive subgroups in order to apply \cref{thm:main_result}.
\end{proof}

\begin{remark}
The lower bounds in \eqref{eqn:lower_bounds} ensure that \cref{cor:approximate-dedekind} is not vacuous.
\end{remark}

Our proofs for \cref{thm:extremal-class-number,thm:ell-torsion,thm:fibered-chebotarev} rely on the bounds \eqref{eqn:partialsums_n} and \eqref{eqn:partialsums_p}.  Our proof of \cref{thm:subconvexity-nice-families}, and hence our proof of \cref{thm:torus-orbits}, uses the strong zero-free region in \cref{cor:approximate-dedekind}.

\subsection{Further examples}

As is made clear in the previous sections, our results are strongest for groups $G$ possessing a unique minimal normal subgroup $N$ for which $\mathrm{T}(G,N)$ holds.  In particular, we obtain the following analogue of Theorem \ref{thm:artin-pnt-intro} for such groups $G$.

\begin{theorem}\label{thm:artin-pnt-later}
Let $k$ be a number field.  Let $G$ be a finite group with a unique minimal normal subgroup $N$ such that Hypothesis \hyperlink{hyp:T}{$\mathrm{T}(G,N)$} holds. Let $Q\geq 1$.  For all $\epsilon>0$, there exists an  effectively computable constant $\Cr{main}=\Cr{main}(|G|,[k:\Q],\epsilon)>0$ such that for all except $O_{|G|,[k:\Q],\epsilon}(Q^\epsilon)$ normal extensions $K/k$ with $\mathrm{Gal}(K/k) \simeq G$ and absolute discriminant $D_K$ at most $Q$, each irreducible faithful Artin representation $\rho$ of $\mathrm{Gal}(K/k)$ satisfies 
	\[
	\Big|\sum_{ \mathrm{N}_{k/\mathbb{Q}} \mathfrak{p} \leq x} \mathop{\mathrm{tr}} \rho(\mathrm{Frob_\mathfrak{p}})\Big|\ll_{|G|,[K:\Q],\epsilon} x \exp(-\Cr{main} \sqrt{\log x})
	\]
	for all $ x \geq (\log D_K)^{81|G|/\epsilon}$.
\end{theorem}

There are many groups $G$ satisfying both hypotheses of Theorem \ref{thm:artin-pnt-later}.  For example, there are $50$ transitive groups of degree $8$, of which $42$ are subject to Theorem \ref{thm:artin-pnt-later}, and of the $1954$ transitive groups of degree $16$, there are $1706$ subject to Theorem \ref{thm:artin-pnt-later}.  We therefore do not aim to provide an exhaustive list of such groups.  Instead, we highlight a few systematic examples beyond the symmetric groups $S_n$ and transitive groups of prime degree that have played a role earlier in this paper.

\begin{itemize}
	\item All simple groups $G$, with $N=G$.
	\item The non-simple alternating group $A_4$, with $N=V_4$, the Klein four subgroup.
	\item The affine general linear group $\mathrm{AGL}_d(\mathbb{F}_p)$, with $N=(\mathbb{Z}/p\mathbb{Z})^d$, and more generally any group of the form $G_0 \rtimes (\mathbb{Z}/p\mathbb{Z})^d$ with $G_0$ an irreducible subgroup of $\mathrm{GL}_d(\mathbb{F}_p)$.
	\item Almost simple groups $G$ whose socle has prime power index in $G$.
	\item Wreath products $S_3 \wr H$ and $S_4 \wr H$ for transitive permutation groups $H$ of degree $d$, with $N = A_3^d$ and $N=V_4^d$, respectively.
\end{itemize}

All but the last of these are primitive permutation groups, and consequently it is straightforward to obtain a version of Theorem \ref{thm:ell-torsion} for these groups equal in quality to that for $S_n$ and groups of prime degree.  It is also possible to obtain a version for imprimitive groups, but with a somewhat worse bound on the $\ell$-torsion subgroup.  Additionally, for any group $G$ with a unique minimal normal subgroup $N$ for which Hypothesis \hyperlink{hyp:T}{$\mathrm{T}(G,N)$} holds, primitive or otherwise, an analogue of Theorem \ref{thm:fibered-chebotarev-Sn} holds with the quadratic resolvent replaced by the subfield $K^N$ fixed by $N$.

\begin{theorem}\label{thm:fibered-chebotarev}
	Let $G$ be a finite group with a unique minimal normal subgroup $N$ for which Hypothesis \hyperlink{hyp:T}{$\mathrm{T}(G,N)$} holds.  Let $k$ be a number field.  Let $\mathcal{C} \subseteq G$ be a conjugacy class and let $[\mathcal{C}]_{G/N}$ denote the associated conjugacy class in $G/N$.  Let $Q\geq 1$.  For all $\epsilon>0$, there exists a constant $\Cr{main}=\Cr{main}(|G|,[k:\Q],\epsilon)>0$ such that for all except $O_{n,[k:\Q],\epsilon}(Q^{\epsilon})$ fields $K \in \mathfrak{F}^{G}_k(Q)$, one has that for any $x \geq (\log D_K)^{81|G|/\epsilon}$, there holds
		\[
			\pi_\mathcal{C}(x; K/k)
				= \frac{|\mathcal{C}|}{|N|\cdot|[\mathcal{C}]_{G/N}|} \pi_{[\mathcal{C}]_{G/N}}(x; K^N/k) + O_{n,[k:\Q],\epsilon}(x \exp(-\Cr{main}\sqrt{\log x})),
		\]
	where $K^N$ denotes the subfield of $K$ fixed by $N$.
\end{theorem}

\section{Preliminaries on Artin $L$-functions} \label{sec:artin-prelim}

We recall the definition of an Artin $L$-function following \cite[Chapter 2]{MurtyMurty}.  Let $K/k$ be a Galois extension of number fields with Galois group $G=\Gal(K/k)$. Let $\mathcal{O}_k$ be the ring of integers of $k$.  For each prime $\kp$ of $k$ and each prime $\mathfrak{P}$ of $K$ lying over $\mathfrak{p}$, let $D_{\mathfrak{P}}=\Gal(K_{\mathfrak{P}}/k_{\kp})$, where $K_{\mathfrak{P}}$ and $k_{\kp}$ are the completions of $K$ and $k$ at $\mathfrak{P}$ and $\kp$, respectively.  Let $F_\mathfrak{P}$ and $F_\mathfrak{p}$ denote the residue fields of $\mathfrak{P}$ and $\mathfrak{p}$.  There is a map from $D_{\mathfrak{P}}$ to $\Gal(F_{\mathfrak{P}}/F_{\kp})$ that is surjective by Hensel's lemma.  Define $I_{\mathfrak{P}}$ to be the kernel of this map; we then have an exact sequence
\[
1\to I_{\mathfrak{P}}\to D_{\mathfrak{P}}\to \Gal(F_{\mathfrak{P}}/F_{\kp})\to 1.
\]
The group $\Gal(F_{\mathfrak{P}}/F_{\kp})$ is cyclic with generator $x\mapsto x^{\N\kp}$.  Choose $\sigma_{\mathfrak{P}}\in D_{\mathfrak{P}}$ whose image in $\Gal(F_{\mathfrak{P}}/F_{\kp})$ is this generator; it is only defined modulo $I_{\mathfrak{P}}$.  We have $I_{\mathfrak{P}}=1$ for all unramified $\kp$, so for these $\kp$, $\sigma_{\mathfrak{P}}$ is well-defined.   If we choose another prime $\mathfrak{P}'$ above $\kp$, then $I_{\mathfrak{P}'}$ and $D_{\mathfrak{P}'}$ are conjugates of $I_{\mathfrak{P}}$ and $D_{\mathfrak{P}}$.  For $\kp$ unramified, we denote by $\sigma_{\kp}$ the conjugacy class of Frobenius automorphisms at primes $\mathfrak{P}$ above $\kp$.

Let $\rho\colon G \to\mathrm{GL}_n(\mathbb{C})$ be a complex representation of $G$, and let $V$ be the underlying complex vector space on which $\rho$ acts.  We may restrict this action to the decomposition group $D_{\mathfrak{P}}$ and see that the quotient $D_{\mathfrak{P}}/I_{\mathfrak{P}}$ acts on the subspace $V^{I_{\mathfrak{P}}}$ of $V$ on which $I_{\mathfrak{P}}$ acts trivially.  Any $\sigma_{\mathfrak{P}}$ will have the same characteristic polynomial on this subspace.  For $\re(s)>1$, we define
\[
L_\kp(s,\rho)=\det(1 - \rho(\sigma_{\mathfrak{P}})|V^{I_{\mathfrak{P}}}\N \kp^{-s})^{-1} = \prod_{j=1}^{n}(1-\alpha_{j,\rho}(\kp)\N \kp^{-s})^{-1}.
\]
Note that the matrix $\rho(\sigma_{\mathfrak{P}})|V^{I_{\mathfrak{P}}}$ remains the same if one changes the prime $\mathfrak{P}$ lying above $\kp$; indeed, if $\kp$ is unramified, then $\rho(\sigma_{\mathfrak{P}})|V^{I_{\mathfrak{P}}}=\rho(\sigma_{\kp})$.  We then define
 \begin{equation}
 \label{eqn:Euler_prod}
L(s,\rho) = \prod_{\kp}L_{\kp}(s,{\rho})=\sum_{\kn}\frac{\lambda_{\rho}(\kn)}{\N\kn^s}.
\end{equation}
We have that $|\alpha_{j,\rho}(\kp)|\leq 1$ for all $j$ and $\kp$, so $L(s,\rho)$ has an absolutely convergent Dirichlet series and Euler product for $\re(s)>1$.

Let $\Gamma_{\R}(s) := \pi^{-s/2}\Gamma(s/2)$.  For each archimedean place $v$ of $k$, we define
\[
L_v(s,\rho) = \begin{cases}
	\Gamma_{\R}(s)^{n}\Gamma_{\R}(s+1)^n&\mbox{if $k_v=\mathbb{C}$,}\\
\Gamma_{\R}(s)^a \Gamma_{\R}(s+1)^{n-a}&\mbox{if $k_v=\R$,}
\end{cases}
\]
where $a = a(\rho)$ is the dimension of the $+1$ eigenspace of complex conjugation.  We define the numbers $\mu_{\rho}(j)$ by the identity
\[
L_{\infty}(s,\rho)=\prod_{\textup{$v$ archim.}}L_v(s,\rho)=\prod_{j=1}^{n[k:\Q]}\Gamma_{\R}(s+\mu_{\rho}(j)).
\]
Let the integral ideal $\kq_{{\rho}}\subseteq\mathcal{O}_k$  denote the conductor of ${\rho}$ over $k$.  The completed $L$-function is defined by
\begin{equation}
\label{eqn:ext_Lambda_def}
\Lambda(s,{\rho}):= (D_k^n \N_{k/\Q}\kq_{{\rho}})^{s/2}L(s,{\rho})L_{\infty}(s,{\rho}).
\end{equation}
There exists $W(\rho) \in \mathbb{C}$ of modulus one such that
\[
	\Lambda(s,\rho) = W(\rho) \Lambda(1-s,\overline{\rho})
\]
for all $s\in\mathbb{C}$ at which $\Lambda(s,\rho)$ is holomorphic, where $\overline{\rho}$ is the complex conjugate of $\rho$.  We define the analytic conductor of $\rho$ by
\begin{equation}
\label{eqn:analytic_conductor}
C(\rho,t):=D_k^n \N_{k/\Q}\kq_{\rho}\prod_{j=1}^{n[k:\Q]}(3+|\mu_{\rho}(j)+it|),\qquad C(\rho):=C(\rho,0).
\end{equation}
We observe that
\begin{equation}
\label{eqn:t-bound}
C(\rho,t)\ll_{n,[k:\Q]}D_k^n \mathrm{N}_{k/\Q}\kq_{\rho}(3+|t|)^{n[k:\Q]}.
\end{equation}

\begin{lemma}
	\label{lem:convexity}
	If $\rho$ is an $n$-dimensional Artin representation over $k$ whose $L$-function $L(s,\rho)$ is entire, then for $\re(s)\geq 1/2$, then
	\[
	|L(s,\rho)|\ll_{n,[k:\Q]}C(\rho,t)^{\max\{\frac{1-\sigma}{2},0\}}(\log C(\rho,t))^{d[k:\Q](2\sigma-1)}.
	\]
\end{lemma}
\begin{proof}
	The bound $|L(1+it,\rho)|\ll_{n,[k:\Q]}(\log C(\rho,t))^{n[k:\Q]}$ follows by proceeding as in the proof of \cite[Theorem 2]{MR1276986}.  The bound $|L(\frac{1}{2}+it,\rho)|\ll_{n,[k:\Q]}C(\rho,t)^{n[k:\Q]/4}$ follows from the convexity bound due to Heath-Brown \cite{HB}.  
\end{proof}

If $\rho$ is $1$-dimensional, then Artin reciprocity shows that $L(s,\rho)$ is a Hecke $L$-function and is thus entire if $\rho$ is nontrivial.  The Artin conjecture asserts that $L(s,\rho)$ is entire for every nontrivial irreducible representation of $G=\mathrm{Gal}(K/k)$, but this is unknown in general.  The best general result is due to Brauer, and is a consequence of his induction theorem.

\begin{lemma}[Brauer induction]\label{lem:brauer-induction}
For any complex representation $\rho$ of $\mathrm{Gal}(K/k)$, the Artin $L$-function $L(s,\rho)$ has a meromorphic continuation to $\mathbb{C}$.
\end{lemma}

We shall speak interchangably about the Artin $L$-function associated to $\rho$ and to its character $\chi_\rho$ (i.e., its trace).  For example, if $\mathbf{1}_G$ denotes the character of the trivial representation of $G$, then $L(s,\mathbf{1}_G) = \zeta_k(s)$.  Also, for any $\rho$, if there are rational coefficients $c_i$ such that
\[
	\chi_\rho = \sum_{i=1}^k c_i \chi_{\rho_i},
\]
where for each $i \leq k$, $\rho_i$ is complex representation of $G$, then Artin showed
\[
	L(s,\rho) = L(s,\chi_\rho) = \prod_{i=1}^k L(s,\chi_{\rho_i})^{c_i} = \prod_{i=1}^k L(s,\rho_i)^{c_i}. 
\]
Finally, if $\chi$ is a character of a subgroup $H \subset G$, let $\chi^* = \mathrm{Ind}_H^G \chi$ denote the character of $G$ induced by $\chi$.  Then there is an equality of the associated $L$-functions,
\[
L(s,\chi) = L(s,\chi^*),
\]
where $L(s,\chi)$ is an Artin $L$-function associated to $K/K^H$ and $K^H$ denotes the subfield fixed by $H$.

\section{Transfer of zero-free regions}
\label{sec:transfer}

The goal of this section is to prove Theorems \ref{thm:zfr-transfer} and \ref{thm:conj}, whose setup we briefly recall.  Let $K/k$ be a normal extension of number fields with Galois group $G$.  Let $N \trianglelefteq G$ be a normal subgroup and let $K^N$ denote the fixed field of $N$.  We wish to show that if $\zeta_K(s)/\zeta_{K^N}(s)$ is non-zero in a region $\Omega \subseteq \mathbb{C}$, then $L(s,\rho)$ is holomorphic and non-vanishing on $\Omega$ for every irreducible Artin representation $\rho$ of $K/k$ whose kernel does not contain $N$.  
We shall make this more precise shortly, but in loose terms, the idea is to first show that a zero-free region for $\zeta_K(s)/\zeta_{K^N}(s)$ transfers to $L$-functions attached to certain $1$-dimensional characters of subgroups $H \subseteq G$.  We do so in Lemma \ref{lem:abelian-transfer} for the largest possible set of characters for which this conclusion could reasonably hold, namely those whose kernel does not contain $H \cap N$.  The next step is to show that this transfer of zero-free region to $L$-functions associated to $1$-dimensional characters suffices.  
This transfer is equivalent to showing that the inductions of these $1$-dimensional characters generate all characters of $G$ whose kernel does not contain $N$.   This is Hypothesis \hyperlink{hyp:T}{$\mathrm{T}(G,N)$}, which we restate in an equivalent manner below.  Conjecture \ref{conj:restricted-induction} asserts that Hypothesis \hyperlink{hyp:T}{$\mathrm{T}(G,N)$} holds for all finite groups $G$ and normal subgroups $N \unlhd G$.
If Hypothesis \hyperlink{hyp:T}{$\mathrm{T}(G,N)$} holds, then Theorem \ref{thm:zfr-transfer} readily follows.  The verification of Hypothesis \hyperlink{hyp:T}{$\mathrm{T}(G,N)$} is the most subtle piece of the argument.  We are unable to prove Hypothesis \hyperlink{hyp:T}{$\mathrm{T}(G,N)$} holds in general; however, we do show in Theorem \ref{thm:orthogonal-span} that it holds in many natural cases, most notably when either the index $[G:N]$ is a prime power or $N$ is solvable.  We then provide the proof of Theorem \ref{thm:minimal-normal} that reduces the general Conjecture \ref{conj:restricted-induction} to proving Hypothesis \hyperlink{hyp:T}{$\mathrm{T}(G,N)$} for a concrete set of groups.
We close this section by illuminating the ideas of this paper with the example $G=S_5$ and $N=A_5$.

Our approach is inspired by Brauer's approach to the Aramata--Brauer theorem that the quotient $\zeta_K(s)/\zeta_k(s)$ is entire.  
The ideas of Brauer's proof quickly establish the theorem in the special case $N = G$, and thus with $K^N = k$.  In fact, as Brauer's work is also an ingredient in our proof of the general case, we find it useful to briefly summarize his ideas.

\subsection{The Aramata--Brauer theorem and Theorem \ref{thm:zfr-transfer} when $N=G$}

We first introduce some notation.  Given a subgroup $H \subseteq G$, we let $\mathbf{1}_H$ denote the trivial character of $H$.  Given any character $\chi$ of $H$ and any subgroup $H^\prime \supseteq H$, we let $\mathrm{Ind}_H^{H^\prime} \chi$ denote the character of $H^\prime$ induced by $\chi$.  We begin with the following lemma of Brauer \cite{Brauer}.

\begin{lemma}[Brauer]\label{lem:brauer}
Let $N$ be a finite group.  For each nontrivial character $\chi$ of a cyclic subgroup of $N$, there is a positive rational constant $c_\chi$ such that
\[
\mathrm{Ind}_1^N \mathbf{1} - \mathbf{1}_N = \sum_{\substack{g \in N \\ g \neq \mathrm{id}}} \sideset{}{^\prime}\sum_{\chi \in \widehat{\langle g\rangle}} c_\chi \mathrm{Ind}_{\langle g\rangle}^N \chi,
\]
where $\chi$ runs over the nontrivial characters of the cyclic group $\langle g\rangle$.
\end{lemma}
\begin{proof}
Brauer shows that if $\chi$ is a nontrivial character of the cyclic subroup $\langle g \rangle$, one may take
\[
c_\chi
    = \frac{1}{|N|} \sum_{\substack{h \in \langle g\rangle \\ \langle h \rangle = \langle g \rangle}} (1 - \bar\chi(h)).
\]
This is rational (it is invariant under the action of Galois) and positive (because $\chi$ is nontrivial).
\end{proof}

From Lemma \ref{lem:brauer}, it is straightforward to deduce the Aramata--Brauer theorem and a first consequence toward non-vanishing.

\begin{lemma}\label{lem:cyclic-transfer}
The quotient $\zeta_K(s)/\zeta_{K^N}(s)$ is entire.  Moreover, let $H \subseteq N$ be a cyclic subgroup and let $\chi$ be a nontrivial irreducible character of $H$.  If $\zeta_K(s)/\zeta_{K^N}(s)$ is non-zero in a region $\Omega \subseteq \mathbb{C}$, then $L(s,\chi)$ is holomorphic and non-vanishing in $\Omega$.
\end{lemma}
\begin{proof}
By Lemma \ref{lem:brauer}, it follows for $\re(s) > 1$ that
\[
\frac{\zeta_K(s)}{\zeta_{K^N}(s)} = \prod_{\substack{g \in N \\ g \neq \mathrm{id}}} \sideset{}{^\prime}\prod_{\chi \in \widehat{\langle g \rangle}} L(s,\chi)^{c_\chi}.
\]
Each character $\chi$ is abelian, hence each $L(s,\chi)$ is entire. Since $\zeta_K(s)/\zeta_{K^N}(s)$ is meromorphic and since each $c_\chi$ is positive, it follows that $\zeta_K(s)/\zeta_{K^N}(s)$ is entire.  The statement on non-vanishing similarly follows from the positivity of each $c_\chi$ and the fact that each $L(s,\chi)$ is entire.
\end{proof}

For the next lemma, we introduce some additional notation.  For any finite group $G$, let $\mathcal{R}_\mathbb{Z}(G)$ denote the ring of virtual characters of $G$, i.e. integral linear combinations of the irreducible characters of $G$, and let $\mathcal{R}_{\mathbb{Q}}(G) = \mathcal{R}_\mathbb{Z}(G) \otimes \mathbb{Q}$ and $\mathcal{R}_\mathbb{C}(G) = \mathcal{R}_\mathbb{Z}(G) \otimes \mathbb{C}$.  The space $\mathcal{R}_\mathbb{C}(G)$ is the space of class functions on $G$, 
on which there is an inner product
\[
\langle f_1, f_2 \rangle_G = \frac{1}{|G|} \sum_{g \in G} f_1(g) \overline{f_2(g)}.
\]
The irreducible characters of $G$ form an orthonormal basis for $\mathcal{R}_\mathbb{C}(G)$ with respect to this inner product.  The orthogonal complement of the trivial character is the space of class functions with mean $0$ on $G$.  By Frobenius reciprocity, this condition is invariant under induction.  We also note that this inner product is defined on $\mathcal{R}_\mathbb{Q}(G)$ as well.

Our next lemma is essentially the Artin induction theorem. 

\begin{lemma}[Artin induction]\label{lem:artin-induction}
The orthogonal complement of the trivial character in $\mathcal{R}_\mathbb{Q}(G)$ is spanned by the induction of nontrivial characters of cyclic subgroups of $G$.
\end{lemma}
\begin{proof}
It suffices to prove the analogous statement over $\mathbb{C}$.  Given $g \neq \mathrm{id}$, let $f_g$ be the class function on the cyclic group $H=\langle g\rangle$ defined by
\[
f_g(h) = \begin{cases} 1, & \mbox{if $h=g$,}\\
-1 & \mbox{if $h = \mathrm{id}$,}\\
0 & \mbox{if $h\notin\{g, \mathrm{id}\}$.}
\end{cases}
\]
As $f_g$ has mean $0$ on $H$, it may be expressed as a linear combination of the non-trivial characters of $H$.  Its induction $\mathrm{Ind}_H^G f_g$ thus also has mean $0$, and is supported on the identity and the conjugacy class of $g$.  Varying over all $g \neq \mathrm{id}$, such functions naturally span the orthogonal complement of the trivial character.
\end{proof}

Combining Lemmas \ref{lem:cyclic-transfer} and \ref{lem:artin-induction}, we obtain:

\begin{lemma}\label{lem:normal-transfer}
Suppose that $\zeta_{K}(s)/\zeta_{K^N}(s)$ is non-zero in a region $\Omega \subseteq \mathbb{C}$.  Let $\psi \in \mathcal{R}_\mathbb{Q}(N)$ lie in the orthogonal complement of the trivial character of $N$.  Then the Artin $L$-function $L(s,\psi)$ is holomorphic and non-vanishing in $\Omega$.
\end{lemma}
\begin{proof}
Applying Lemma \ref{lem:artin-induction} to $N$, there are rational numbers $c_{\chi,\psi}$ such that 
\[
	\psi = \sum_{\substack{g \in N\\ g \neq \mathrm{id}}} \sum_{\substack{ \chi \in \widehat{\langle g\rangle} \\ \chi \neq 1}} c_{\chi,\psi} \mathrm{Ind}_{\langle g \rangle}^N \chi,
\]
where the sum runs over the nontrivial characters of $\langle g \rangle$.  Consequently, we find the factorization
\[
	L(s,\psi)
		= \prod_{\substack{g \in N\\ g \neq \mathrm{id}}} \prod_{\substack{\chi \in \widehat{\langle g\rangle} \\ \chi \neq 1}} L(s,\chi)^{c_{\chi,\psi}}.
\]
By Lemma \ref{lem:cyclic-transfer}, each $L(s,\chi)$ is holomorphic and non-vanishing in $\Omega$, so the same must hold for $L(s,\psi)$ as well.
\end{proof}

Lemma \ref{lem:normal-transfer} yields Theorem \ref{thm:zfr-transfer} in the case $N=G$.  The general case is apparently more subtle, however.

\subsection{Inductions of characters with restricted components}

Recall that there is a natural injection $\mathcal{R}_\mathbb{Q}(G/N) \hookrightarrow \mathcal{R}_\mathbb{Q}(G)$ given by pullback, so we may regard $\mathcal{R}_\mathbb{Q}(G/N)$ as a subgroup of $\mathcal{R}_\mathbb{Q}(G)$.  Exploiting the inner product on $\mathcal{R}_\mathbb{Q}(G)$, we may thus consider an orthogonal decomposition
\begin{equation}\label{eqn:orthogonal-decomp}
    \mathcal{R}_\mathbb{Q}(G) 
        = \mathcal{R}_\mathbb{Q}(G/N) \oplus \mathcal{R}_\mathbb{Q}(G/N)^\perp,
\end{equation}
where $\mathcal{R}_\mathbb{Q}(G/N)^\perp$ denotes the orthogonal complement of $\mathcal{R}_\mathbb{Q}(G/N)$.  
So doing, $\mathcal{R}_\mathbb{Q}(G/N)$ is spanned by the irreducible constituents of the character $\mathrm{Ind}_N^G \mathbf{1}_N$.  Consequently, $\mathcal{R}_\mathbb{Q}(G/N)^\perp$ is most naturally spanned by the irreducible characters whose kernel does not contain $N$.  These are exactly the characters of the type considered in Theorem \ref{thm:zfr-transfer}.  However, as in the previous subsection, the analytic properties of the $L$-functions attached to such characters are not easily accessed directly and are only understood by comparison to $L$-functions attached to $1$-dimensional characters.

Hypothesis \hyperlink{hyp:T}{$\mathrm{T}(G,N)$} asserts a refinement of the Artin induction theorem that respects the orthogonal decomposition \eqref{eqn:orthogonal-decomp}.  It follows from Lemma \ref{lem:artin-induction} that $\mathcal{R}_\mathbb{Q}(G/N)$ is spanned by the induction of one-dimensional characters whose induction itself lies in $\mathcal{R}_\mathbb{Q}(G/N)$; in fact, characters of cyclic subgroups suffice.  In other words, it follows that
\[
	\mathcal{R}_\mathbb{Q}(G/N) = \mathrm{span}_\mathbb{Q} \bigcup_{H \subseteq G} \{ \mathrm{Ind}_H^G \chi \colon \chi \in \mathrm{Irr}(H),~\mathrm{dim} \chi = 1, \text{ and } \mathrm{Ind}_H^G \chi \in \mathcal{R}_\mathbb{Q}(G/N)\},
\]
where $\mathrm{Irr}(H)$ denotes the set of irreducible complex representations of a subgroup $H \subset G$.

Hypothesis \hyperlink{hyp:T}{$\mathrm{T}(G,N)$} is equivalent to the corresponding statement for $\mathcal{R}_\mathbb{Q}(G/N)^\perp$.  

\begin{hypothesis}[Equivalent formulation of Hypothesis \hyperlink{hyp:T}{$\mathrm{T}(G,N)$}]
	The space $\mathcal{R}_{\mathbb{Q}}(G/N)^\perp$ is spanned by the induction of one-dimensional characters whose induction lies in $\mathcal{R}_{\mathbb{Q}}(G/N)^\perp$.  In other words, we have that
	\[
		\mathcal{R}_\mathbb{Q}(G/N)^\perp = \mathrm{span}_\mathbb{Q} \bigcup_{H \subseteq G} \{ \mathrm{Ind}_H^G \chi \colon \chi \in \mathrm{Irr}(H),~ \dim\chi = 1, \text{ and } \mathrm{Ind}_H^G \chi \in \mathcal{R}_\mathbb{Q}(G/N)^\perp\}.
	\]
\end{hypothesis}

The main result of this section, Theorem \ref{thm:orthogonal-span} below, establishes Hypothesis \hyperlink{hyp:T}{$\mathrm{T}(G,N)$} when $N$ is solvable or the index $[G:N]$ is a prime power.  Its proof concludes that of \cref{thm:conj}.  

The subtlety of Hypothesis \hyperlink{hyp:T}{$\mathrm{T}(G,N)$} is revealed upon noting that the orthogonal decomposition \eqref{eqn:orthogonal-decomp} is preserved neither upon restriction to subgroups $H \subseteq G$ nor upon induction from subgroups.  For example, the trivial character of $H$ always lies in $\mathcal{R}_\mathbb{Q}(H/H\cap N)$, but its induction to $G$ need not lie in $\mathcal{R}_{\mathbb{Q}}(G/N)$.  Thus, typical elementary approaches exploiting the adjointness of induction and restriction do not obviously apply.  Moreover, it is not the case that the induction of characters of cyclic, or even abelian, subgroups suffice to span $\mathcal{R}_\mathbb{Q}(G/N)^\perp$ in general.

We first prove a non-vanishing result for $L$-functions attached to the characters we shall use.  For this, we note that if $\chi$ is a $1$-dimensional character of a subgroup $H \subseteq G$, then $\mathrm{Ind}_H^G \chi \in \mathcal{R}_\mathbb{Q}(G/N)^\perp$ if and only if the kernel of $\chi$ does not contain $H \cap N$.

\begin{lemma}
\label{lem:abelian-transfer}
Let $H \subseteq G$ be a subgroup and suppose that $\chi$ is a $1$-dimensional character of $H$ whose kernel does not contain $H \cap N$.  If $\zeta_K(s)/\zeta_{K^N}(s)$ is non-vanishing in a region $\Omega \subseteq \mathbb{C}$, then $L(s,\chi)$ is holomorphic and non-vanishing in $\Omega$.
\end{lemma}
\begin{proof}
Since $\chi$ is $1$-dimensional, it follows from Artin reciprocity that $L(s,\chi)$ is a Hecke $L$-function, and since $\chi$ is nontrivial, it is thus entire.  To show that $L(s,\chi)$ is non-vanishing in $\Omega$, let $\chi^* = \mathrm{Ind}_H^G \chi$, and consider the character $\theta_N = \mathrm{Ind}_N^G \mathbf{1}_N - \mathbf{1}_G$ of $G$.  Then $\chi^* + \theta_N \chi^* = \chi^* \mathrm{Ind}_N^G \mathbf{1}_N = \mathrm{Ind}_N^G \chi^*_N$, where $\chi^*_N := \chi^*|_N$ denotes the restriction of $\chi^*$ to $N$.  
It follows that
\[
L(s,\chi) L(s,\chi^* \otimes \theta_N) = L(s,\chi^*_N).
\]
By our assumption on $\chi$, the character $\chi^*_N$ of $N$ is orthogonal to the trivial character.  Thus, by Lemma \ref{lem:normal-transfer}, $L(s,\chi^*_N)$ is holomorphic and non-vanishing on $\Omega$.  Consequently, the lemma will follow provided we show that $L(s, \chi^* \otimes \theta_N)$ is entire.

By Lemma \ref{lem:brauer} applied to $G/N$, we may express $\theta_N$ as a positive rational linear combination of the induction $\psi^*$ of characters $\psi$ of cyclic subgroups of $G/N$.  Since $L(s,\chi^* \otimes \theta_N)$ is meromorphic by Lemma \ref{lem:brauer-induction}, it suffices to show that each $L(s, \chi^* \otimes \psi^*)$ is entire.  Let $K^\psi \subseteq K^N$ be the cyclic subextension from which $\psi$ is induced and let $F^\prime/F$ be the cyclic extension corresponding to $\chi$.  Let $\psi_F$ denote the restriction of $\psi$ to the subgroup $\mathrm{Gal}(FK^N/FK^\psi) \subseteq \mathrm{Gal}(K^N/K^\psi)$.  

If $F^\prime$ is linearly disjoint from $FK^\psi$ over $F$, then we may regard both $\chi$ and $\psi_F$ as characters of the abelian Galois group $\mathrm{Gal}(F^\prime K^N/FK^\psi)$.  Thus, the product $\chi \psi_F$ is well-defined, and we have $L(s,\chi^*\otimes \psi^*) = L(s,\chi \psi_F)$.  Our assumption on $\chi$ ensures that $F^\prime$ is not contained in $FK^N$, so the character $\chi\psi_F$ is non-trivial, and hence $L(s,\chi^* \otimes \psi^*)$ is entire in this case.

If $F^\prime$ is not linearly disjoint from $FK^\psi$ over $F$, the character $\chi \psi_F$ need not be defined.  However, we may regard $\chi$ as a character of $G_F = \mathrm{Gal}(K^NF^\prime / F)$ and $\psi_F$ as a character of its subgroup $G_{\psi} = \mathrm{Gal}(K^NF^\prime / FK^\psi)$.  So doing, we find $\chi \mathrm{Ind}_{G_\psi}^{G_F} \psi_F = \mathrm{Ind}_{G_\psi}^{G_F} (\psi_F \chi_\psi)$, where $\chi_\psi := \chi|_{G_\psi}$ denotes the restriction of $\chi$ to $G_\psi$.  Then $L(s, \chi^* \otimes \psi^*) = L(s, \psi_F \chi_\psi)$.  As above, by our assumption on $\chi$, the (abelian) character $\psi_F \chi_\psi$ is non-trivial, and we again conclude that $L(s,\chi^* \otimes \psi^*)$ is entire.  This completes the proof of the lemma.
\end{proof}

\begin{theorem}\label{thm:orthogonal-span}
Let $G$ be a finite group and $N \trianglelefteq G$ a normal subgroup.  If $N$ is solvable or the index $[G:N]$ is a prime power, then Hypothesis \hyperlink{hyp:T}{$\mathrm{T}(G,N)$} holds.
\end{theorem}
\begin{proof}
Notice that the induction of a character $\chi$ from a subgroup $H$ is in the orthogonal complement of $\mathcal{R}_{\Q}(G/N)$ if and only if its kernel does not contain $H \cap N$.  Thus, the conclusion of the theorem holds if every irreducible character of $G$ is the induction of a $1$-dimensional character of a subgroup of $G$, i.e. if $G$ is monomial.  Consequently, the theorem holds whenever $G$ is abelian or, more generally, whenever $G$ is nilpotent \cite[Chapter 8]{Serre_linear}.  
In fact, all monomial groups are solvable, so this fact is subsumed by the statement of the theorem.  We find it convenient to note also that all $p$-groups are nilpotent, so the statement is known in this case.  It also holds when $N=G$ by Lemma \ref{lem:artin-induction}.  Finally, we note that it is also known when $N$ is abelian by work of Deligne and Henniart \cite[Proposition 2.2]{DeligneHenniart}.

We first consider the case that $N$ is solvable.  We proceed by induction on the order of $G$, taking as a base case the statement for abelian $G$.  If $N$ is trivial, then the conclusion is vacuous, so we may assume that $N$ is nontrivial.  Let $N^\prime$ be a nontrivial minimal normal subgroup of $G$ contained in $N$.  Since $N$ is solvable, $N^\prime$ must be solvable as well.  Additionally, since $N^\prime$ must be characteristically simple, it is isomorphic to a group of the form $T^d$ for a simple group $T$.  Since $N^\prime$ is solvable, $T$ must be cyclic of prime order, so $N^\prime$ is abelian.  Consequently, by \cite[Proposition 2.2]{DeligneHenniart}, it follows that {$\mathrm{T}(G,N^\prime)$} holds.

If $N^\prime = N$, then this establishes the result.  If $N^\prime \subsetneq N$, then $\mathcal{R}_\mathbb{Q}(G/N)$ is a proper subspace of $\mathcal{R}_\mathbb{Q}(G/N^\prime)$ and $\mathcal{R}_\mathbb{Q}(G/N^\prime)^\perp$ is a subspace of $\mathcal{R}_\mathbb{Q}(G/N)^\perp$, where all spaces are viewed as subspaces of $\mathcal{R}_\mathbb{Q}(G)$.  Moreover,
	\[
		\mathcal{R}_\mathbb{Q}(G/N)^\perp 
			= \mathcal{R}_\mathbb{Q}(G/N^\prime)^\perp \oplus \big( \mathcal{R}_\mathbb{Q}(G/N^\prime) \cap \mathcal{R}_\mathbb{Q}(G/N)^\perp \big).
	\]
Since {$\mathrm{T}(G,N^\prime)$} holds, $\mathcal{R}_\mathbb{Q}(G/N^\prime)^\perp$ is spanned by the induction of $1$-dimensional characters whose inductions lie in $\mathcal{R}_\mathbb{Q}(G/N^\prime)^\perp$, which must also lie in $\mathcal{R}_\mathbb{Q}(G/N)^\perp$.  It remains to show that $\mathcal{R}_\mathbb{Q}(G/N^\prime) \cap \mathcal{R}_\mathbb{Q}(G/N)^\perp$ is spanned by such characters.  However, since $N^\prime \subset N$, $G/N \simeq (G/N^\prime) / (N/N^\prime)$, and the elements of $\mathcal{R}_\mathbb{Q}(G/N^\prime) \cap \mathcal{R}_\mathbb{Q}(G/N)^\perp$ are obtained by pullback from $\mathcal{R}_\mathbb{Q}((G/N^\prime)/(N/N^\prime))^\perp \subseteq \mathcal{R}_\mathbb{Q}(G/N^\prime)$, where we regard $\mathcal{R}_\mathbb{Q}(G/N^\prime)$ in this final statement as its own space and not as a subspace of $\mathcal{R}_\mathbb{Q}(G)$.  Since $N^\prime$ is nontrivial, $\mathrm{T}(G/N^\prime, N/N^\prime)$ must hold by the inductive hypothesis, and we conclude that $\mathcal{R}_\mathbb{Q}((G/N^\prime)/(N/N^\prime))^\perp$ is spanned by the induction of $1$-dimensional characters whose kernel does not contain $N/N^\prime$.  Pulling back a basis of such characters, we deduce that  $\mathcal{R}_\mathbb{Q}(G/N^\prime) \cap \mathcal{R}_\mathbb{Q}(G/N)^\perp$ is spanned by the induction of $1$-dimensional characters whose kernel does not contain $N$.  This yields the result in the case $N$ is solvable.

We now consider the case that $[G:N]$ is a prime power.  We once again induct on the order of $G$, exploiting the subgroups of $G$.  It suffices to work with complex coefficients, that is, to show the orthogonal complement of $\mathcal{R}_\mathbb{C}(G/N)$ inside $\mathcal{R}_\mathbb{C}(G)$ is spanned by the induction of such characters.  We therefore consider a class function $f$ of $G$ that is orthogonal to the induction of all such characters, with the goal of showing that it lies in $\mathcal{R}_\mathbb{C}(G/N)$.  Equivalently, we wish to show that $f$ is constant on cosets of $N$, that is, that $f(a) = f(b)$ whenever $ab^{-1} \in N$.  Consider two such elements $a,b \in G$.

By Frobenius reciprocity and the inductive hypothesis, it follows that if $a$ and $b$ lie in the same proper subgroup $H$, then $f(a)=f(b)$.  In particular, if $a$ and $b$ do not generate $G$, then we may take $H$ to be the subgroup generated by $a$ and $b$.  Thus, we may assume that $a$ and $b$ generate $G$.  This implies that $G/N$ is cyclic, since we have assumed that $aN = bN$.

Next, since the index of $N$ in $G$ is a power of a prime, we may write $[G:N]=p^k$.  By our inductive hypothesis applied to the cyclic subgroup generated by $a$, or noting that such a group is abelian, we see that $f(a) = f(a^m)$ whenever $m \equiv 1 \pmod{p^k}$.  Thus, replacing $a$ by a suitable power if necessary, we may assume that the order of $a$ is a power of $p$.  Similarly, we assume that the order of $b$ is a power of $p$.  Let $H$ be a Sylow $p$-subgroup of $G$ containing $a$, possibly equal to $G$ itself.  Since $f$ is a class function, we may conjugate $b$ if necessary to assume that $b \in H$ as well.  But $H$ is a $p$-group, so the theorem holds for $H$.  As we have $ab^{-1} \in H \cap N$, it follows that $f(a)=f(b)$.  This establishes the theorem in the case $[G:N]$ is a prime power.
\end{proof}

\begin{proof}[Proof of Theorem \ref{thm:zfr-transfer}]

With these facts in hand, the proof of Theorem \ref{thm:zfr-transfer} is now straightforward.  
Any irreducible representation $\rho$ of $G$ whose kernel does not contain $N$ lies in the orthogonal complement of $\mathcal{R}_\mathbb{Q}(G/N)$.  
If Hypothesis \hyperlink{hyp:T}{$\mathrm{T}(G,N)$} holds, and in particular under the assumptions of Theorem \ref{thm:orthogonal-span},
there are rational constants $c_\chi(\rho)$ such that
\[
L(s,\rho) = \prod_\chi L(s,\chi)^{c_\chi(\rho)},
\]
where the product runs over $1$-dimensional characters $\chi$ of subgroups $H$ of $G$ whose kernel does not contain $N$.  By Lemma \ref{lem:abelian-transfer}, if $\zeta_K(s)/\zeta_{K^N}(s)$ is non-vanishing in the region $\Omega \subseteq \mathbb{C}$, then each $L(s,\chi)$ is holomorphic and non-vanishing in $\Omega$.  Thus, the same must be true for $L(s,\rho)$.  This completes the proof.
\end{proof}

\begin{proof}[Proof of Theorem \ref{thm:minimal-normal}]
	The proof of Theorem \ref{thm:orthogonal-span} given above shows it is sufficient to understand Hypothesis \hyperlink{hyp:T}{$\mathrm{T}(G,N)$} for finite groups $G$ with nonabelian minimal normal subgroups $N$ such that $G/N$ is cyclic of non-prime power order.  We begin by showing it is further possible to reduce to the case that $G$ is a subgroup of $\mathrm{Aut}(N)$.
		
	Let $K$ be the kernel of the map $G \to \mathrm{Aut}(N)$ given by conjugation.  Since $N$ is nonabelian and minimal, $K \cap N$ is trivial, and this realizes $G$ as the fiber product $G/N \times_{G/NK} G/K$.  The quotient $G/K$ is naturally a subgroup of $\mathrm{Aut}(N)$ containing $N$, and we suppose that Hypothesis $\mathrm{T}(G/K, N)$ is true.  Since $G/N$ is abelian, any irreducible character of $G$ may be decomposed as a product $(\chi_1, \chi_2)$, where $\chi_1$ is a character of $G/N$ and $\chi_2$ is a character of $G/K$, with two such products being the same if they differ by a factor $(\psi,\bar\psi)$ for a character $\psi$ of the common quotient, $G/NK$.  Let $(\chi_1,\chi_2)$ be an irreducible character of $G$ not contained in $\mathcal{R}_\mathbb{Q}(G/N)$, so in particular $\chi_2$ is an irreducible character of $G/K$ lying in $\mathcal{R}_\mathbb{Q}(G/KN)^\perp$.  Since we have assumed $\mathrm{T}(G/K,N)$ holds, there are constants $c_\chi \in \mathbb{Q}$ such that
		\[
			\chi_2
				= \sum_{H \subseteq G/K} \sum_{\substack{ \chi \in \mathrm{Irr}(H) \\ \chi(1) = 1 \\ \ker \chi \not \supseteq H \cap N}} c_\chi \mathrm{Ind}_H^{G/K} \chi.
		\]
	For a subgroup $H \subseteq G/K$, let $H_G$ denote the corresponding subgroup of $G$, i.e. $H_G = \{ (a,h) \in G/N \times_{G/NK} G/K  : h \in H\}$.  Then any $1$-dimensional character $\chi$ of $H$ is also a character of $H_G$ and the character $\chi_1$ restricts to a character on $H_G$.  Moreover, if $N \cap H \not\subseteq \ker \chi$, then $N \cap H_G \not\subseteq \ker (\chi \cdot \chi_1\vert_{H_G})$ as well.	 Additionally, we have
		\[
			(\chi_1,\chi_2)
				= \sum_{H \subseteq G/K} \sum_{\substack{ \chi \in \mathrm{Irr}(H) \\ \chi(1) = 1 \\ \ker \chi \not \supseteq H \cap N}} c_\chi \mathrm{Ind}_{H_G}^{G} (\chi \cdot \chi_1\vert_{H_G}),
		\]
	which shows that$\mathrm{T}(G,N)$ holds provided $\mathrm{T}(G/K, N)$ does.  This reduces Conjecture \ref{conj:restricted-induction} to the case that $G$ is a subgroup of $\mathrm{Aut}(N)$, where $N$ is a nonabelian minimal normal subgroup of $G$.
	
	Since $N$ is minimal, it is characteristically simple, and hence of the form $T^d$ for some nonabelian simple group $T$, in which case $G$ is a subgroup of the wreath product $\mathrm{Aut}(T) \wr S_d$.  Applying the reductions in the proof of Theorem \ref{thm:orthogonal-span}, we may assume $G/N$ is cyclic.  This implies that $d$ must be prime, since otherwise the degree $d$ permutation action of $G$ would have nontrivial blocks, violating the assumption that $N$ is minimal.  This yields the theorem.
\end{proof}

\subsection{An example: The symmetric group $S_5$}
\label{subsec:S5}

We illustrate the proof of Theorem \ref{thm:zfr-transfer}, in particular Lemma \ref{lem:abelian-transfer} and Theorem \ref{thm:orthogonal-span}, with one of the simplest interesting examples.  Let $G=S_5$ and $N=A_5$, let $K/k$ be a normal extension with Galois group $S_5$, and assume that $\zeta_K(s)/\zeta_{K^{A_5}}(s)$ is non-vanishing in some region $\Omega \subseteq \mathbb{C}$.  Apart from the trivial character, $S_5$ also admits the sign character, which we denote $\chi_{\mathrm{sgn}}$, whose kernel is equal to $A_5$.  Every other irreducible complex representation of $S_5$ is faithful.  We denote these $\rho_4$, $\rho_4 \otimes \chi_{\mathrm{sgn}}$, $\rho_5$, $\rho_5 \otimes \chi_{\mathrm{sgn}}$, and $\rho_6$, where the subscript indicates the dimension of the representation; for the purposes of this discussion, the choice of which $4$-dimensional representation we denote $\rho_4$ and which we denote $\rho_4 \otimes \chi_\mathrm{sgn}$ is neither important nor illuminating, nor is the choice of $\rho_5$ and $\rho_5 \otimes \chi_\mathrm{sgn}$.
Theorem \ref{thm:zfr-transfer} asserts that the $L$-functions attached to these five faithful representations should be holomorphic and non-vanishing in $\Omega$.

For $L(s,\rho_6)$, this follows from Lemma \ref{lem:normal-transfer} since $\rho_6$ is its own twist by $\chi_{\mathrm{sgn}}$ and is therefore induced from an irreducible representation of $A_5$.  In general, a representation is induced from $A_5$ if and only if it is its own twist by $\chi_\mathrm{sgn}$, and the $A_5$-representation is orthogonal to the trivial character if and only if the associated $L$-function is holomorphic at $s=1$.  Thus, for the other faithful representations, it only follows from Lemma \ref{lem:normal-transfer} that the products $L(s,\rho_4)L(s,\rho_4 \otimes \chi_{\mathrm{sgn}})$ and $L(s,\rho_5)L(s,\rho_5 \otimes \chi_{\mathrm{sgn}})$ are holomorphic and non-vanishing in $\Omega$.  

Let now $\chi_{C_6}$ be either primitive character of the cyclic subgroup $C_6 := \langle(123)(45)\rangle \subseteq S_5$.  The kernel of $\chi_{C_6}$ is trivial, but $C_6 \cap A_5$ is not, so $L(s,\chi_{C_6})$ is subject to Lemma \ref{lem:abelian-transfer}.  Indeed, the proof of Lemma \ref{lem:abelian-transfer} in this case amounts to observing that the character $\chi_{C_6} \otimes \chi_{\mathrm{sgn}} = \chi_{C_6}^4$ is another nontrivial character of $C_6$ whose kernel does not contain $C_6 \cap A_5$.  The product $L(s,\chi_{C_6})L(s,\chi_{C_6} \otimes \chi_{\mathrm{sgn}})$ is then induced from a representation of $A_5$ that must be orthogonal to the trivial representation of $A_5$ since the $L$-function does not have a pole at $s=1$.  It follows from Lemma \ref{lem:normal-transfer} that this product is holomorphic and non-vanishing in $\Omega$, and since each factor is entire, the same holds for both $L(s,\chi_{C_6})$ and $L(s,\chi_{C_6} \otimes \chi_{\mathrm{sgn}})$.

By expressing the induction of $\chi_{C_6}$ in terms of the irreducible representations of $S_5$, we now compute that
\[
    L(s,\chi_{C_6})
        = L(s,\rho_4) L(s,\rho_5)L(s,\rho_5 \otimes \chi_{\mathrm{sgn}}) L(s,\rho_6).
\]
Observe that $L(s,\rho_4)$ may therefore be expressed as a quotient of functions that are holomorphic and non-vanishing in $\Omega$, so it must be holomorphic and non-vanishing as well.  By instead considering $L(s,\chi_{C_6} \otimes \chi_\mathrm{sgn})$, we similarly find that $L(s,\rho_4 \otimes \chi_{\mathrm{sgn}})$ is holomorphic and non-vanishing in $\Omega$.

To show that $L(s,\rho_5)$ and $L(s,\rho_5 \otimes \chi_\mathrm{sgn})$ are holomorphic and non-vanishing in $\Omega$, it is necessary to work with non-abelian subgroups of $S_5$, since a computation reveals that the induction of characters from other abelian subgroups gives at most the same information as $\chi_{C_6}$ above.  Thus, we consider the $2$-Sylow subgroup $D_4 \subseteq S_5$.  The sign character restricts to a non-trivial character of $D_4$, but there are two other quadratic characters of $D_4$ that are twists of each other by $\chi_\mathrm{sgn}$, thus witnessing again the proof of Lemma \ref{lem:abelian-transfer}.  Let $\chi_{D_4}$ denote one of these characters.  Then $L(s,\chi_{D_4})$ is holomorphic and non-vanishing in $\Omega$, and we compute
\[
    L(s,\chi_{D_4})
        = L(s,\rho_4) L(s,\rho_5) L(s,\rho_6).
\]
From the above discussion, both $L(s,\rho_4)$ and $L(s,\rho_6)$ are holomorphic and non-vanishing in $\Omega$, so it follows that $L(s,\rho_5)$ must be too.  Finally, working instead with $L(s,\chi_{D_4} \otimes \chi_\mathrm{sgn})$, we conclude that $L(s,\rho_5 \otimes \chi_\mathrm{sgn})$ is holomorphic and non-vanishing in $\Omega$ as well.

In general, an analysis of the proof of Theorem \ref{thm:orthogonal-span} shows that it always suffices to consider characters of abelian subgroups and of subgroups of $p$-Sylow subgroups for $p \mid [G:N]$.  
Unlike the above example, however, the proof of Theorem \ref{thm:orthogonal-span} proceeds indirectly via class functions to avoid needing a careful understanding of the character theory of $G$ and its subgroups. 
We leave open the questions of whether there is a more direct proof of Theorem \ref{thm:orthogonal-span}, whether Hypothesis \hyperlink{hyp:T}{$\mathrm{T}(G,N)$} holds in general, and whether there is an explicit description of a sufficient set of characters.

\section{The zero density estimate} \label{sec:zero-density}

We now turn to establishing the main zero density estimate, from which Theorem \ref{thm:zero-density-intro} will ultimately follow.  
Let $k$ be a number field, let $G$ be a finite group, and let $N \trianglelefteq G$ be a normal subgroup.  Recall the definitions of $\mathfrak{F}_k^G$ and $\mathfrak{F}_k^G(Q)$.

Given $K\in\mathfrak{F}_{k}^{G}$, define
	\[
	N_{K/K^N}(\sigma,T):=\#\Big\{\beta+i\gamma\colon \beta\geq\sigma,~|\gamma|\leq T,~\frac{\zeta_K(\beta+i\gamma)}{\zeta_{K^N}(\beta+i\gamma)}=0\Big\}.
	\]
	Recall that
	\[
	\mathfrak{m}_{k}^{G,N}(Q)=\max_{K_1\in\mathfrak{F}_{k}^{G}(Q)}\#\{K_2\in\mathfrak{F}_{k}^{G}(Q)\colon K_1\cap K_2\neq K_1^N\cap K_2^N\}.
	\]
We will prove the following zero density estimate.

\begin{theorem}
	\label{thm:ZDE}
	Let $0\leq\sigma\leq 1$, $Q\geq 1$, $T\geq 2$, and $\epsilon>0$. There exists a constant $C_G>0$, depending at most on $G$, such that
	\[
	\sum_{K\in\mathfrak{F}_{k}^{G}(Q)}N_{K/K^N}(\sigma,T)\ll_{|G|,[k:\Q],\epsilon} \mathfrak{m}_{k}^{G,N}(Q)(QT^{|G|[k:\Q]})^{(8+\epsilon)|G|(1-\sigma)}(\log QT)^{C_G [k:\Q]}.
	\]
\end{theorem}

\subsection{Tensor products of certain Artin representations}

For future convenience, given $K \in \mathfrak{F}_{k}^{G}$, let $\psi_K$ denote the character of $G$ given by $\psi_K := \mathrm{Reg}_G - \mathrm{Ind}_N^G \mathbf{1}_N$, where $\mathrm{Reg}_G$ denotes the character of the regular representation of $G$.  Thus,
\[
L(s,\psi_K) = \zeta_K(s)/\zeta_{K^N}(s)
\]
for any $K \in \mathfrak{F}_k^G$.

\begin{lemma}\label{lem:brauer-2}
Let $G_1$ and $G_2$ be finite groups.  Suppose that $\psi_1$ and $\psi_2$ are characters of $G_1$ and $G_2$, respectively, that are positive rational linear combinations of characters induced from non-trivial one-dimensional representations of cyclic subgroups of $G_1$ and $G_2$.

Then the character $\psi_1 \psi_2$ of the direct product $G_1 \times G_2$ is a positive rational linear combination of characters induced from non-trivial one-dimensional representations of subgroups of $G_1 \times G_2$.
\end{lemma}
\begin{proof}
If $\chi_1$ and $\chi_2$ are non-trivial characters of subgroups $H_1 \subseteq G_1$ and $H_2 \subseteq G_2$, then $\chi_1\chi_2$ is a character of $H_1 \times H_2$, and $(\mathrm{Ind}_{H_1}^{G_1} \chi_1)(\mathrm{Ind}_{H_2}^{G_2} \chi_2) = \mathrm{Ind}_{H_1 \times H_2}^{G_1 \times G_2} \chi_1\chi_2$, where the equality is taken as characters of $G_1 \times G_2$.  The result follows.
\end{proof}

\begin{lemma}\label{lem:entire-two}
Let $K_1$ and $K_2$ be distinct normal extensions of $k$ with Galois group $G$.  If $K_1 \cap K_2 = K_1^N \cap K_2^N$ then the Artin $L$-function
\[
L(s,\psi_{K_1} \otimes \psi_{K_2})
\]
is entire, $\kq_{\psi_{K_1}\otimes\psi_{K_2}}$ divides $\mathfrak{D}_{K_1/k}^{[K_2:k]} \mathfrak{D}_{K_2/k}^{[K_1:k]}$, and $\N_{k/\Q}\kq_{\psi_{K_1}\otimes\psi_{K_2}}$ divides $D_{K_1}^{[K_2:k]} D_{K_2}^{[K_1:k]}D_k^{-[K_1:k][K_2:k]}$.
\end{lemma}
\begin{proof}
Let $F = K_1 \cap K_2$.  By assumption, both $K_1^N$ and $K_2^N$ contain $F$, and we may regard the quotients $\zeta_{K_1}(s)/\zeta_{K_1^N}(s)$ and $\zeta_{K_2}(s)/\zeta_{K_2^N}(s)$ as $L$-functions over $F$, say
\[
\frac{\zeta_{K_1}(s)}{\zeta_{K_1^N}(s)} =: L(s,\psi_{K_1/F}) \quad \text{and} \quad \frac{\zeta_{K_2}(s)}{\zeta_{K_2^N}(s)} =: L(s,\psi_{K_2/F}).
\]
By Lemma \ref{lem:brauer}, the characters $\psi_{K_1/F}$ and $\psi_{K_2/F}$ are non-negative linear combinations of characters induced from nontrivial one-dimensional representations.  Moreover, we have $\mathrm{Gal}(K_1K_2/F) \simeq \mathrm{Gal}(K_1/F) \times \mathrm{Gal}(K_2/F)$, so it follows from Lemma \ref{lem:brauer-2} that the same holds for the character $\psi_{K_1/F}\psi_{K_2/F}$.  In particular, the $L$-function
\begin{equation}\label{eqn:product-over-F}
L(s,\psi_{K_1/F} \otimes \psi_{K_2/F}) = \frac{\zeta_{K_1K_2}(s)\zeta_{K_1^NK_2^N}(s)}{\zeta_{K_1K_2^N}(s) \zeta_{K_1^NK_2}(s)}
\end{equation}
is entire.  Next, if we set $G_F = \mathrm{Gal}(K_1K_2/F)$ and $G_k = \mathrm{Gal}(K_1K_2/k)$, then the characters $\psi_{K_1}$ and $\psi_{K_2}$ of $G_k$ are induced by the characters $\psi_{K_1/F}$ and $\psi_{K_2/F}$ of $G_F$.  It follows that
\begin{align*}
\psi_{K_1} \psi_{K_2} 
	&= \Big(\mathrm{Ind}_{G_F}^{G_k} \psi_{K_1/F}\Big) \Big(\mathrm{Ind}_{G_F}^{G_k} \psi_{K_2/F}\Big) \\
	&= \mathrm{Ind}_{G_F}^{G_k} \Big(\psi_{K_1/F}\cdot \mathrm{Res}_{G_F}^{G_k} \mathrm{Ind}_{G_F}^{G_k} \psi_{K_2/F}\Big) = \sum_{\sigma \in G_k/G_F} \mathrm{Ind}_{G_F}^{G_k} \psi_{K_1/F} \psi_{K_2/F}^\sigma,
\end{align*}
where $\psi_{K_2/F}^\sigma$ is the character conjugate to $\psi_{K_2/F}$ via $\sigma$.
However, since $K_2/k$ is Galois and the character $\psi_{K_2/F}$ is valued in $\mathbb{Z}$, $\psi_{K_2/F}^{\sigma} = \psi_{K_2/F}$ for every $\sigma \in G_k/G_F$.  Thus, we conclude that $\psi_{K_1} \psi_{K_2} = [F:k] \mathrm{Ind}_{G_F}^{G_k} \psi_{K_1/F} \psi_{K_2/F}$ and $L(s,\psi_{K_1} \otimes \psi_{K_2}) = L(s,\psi_{K_1/F} \otimes \psi_{K_2/F})^{[F:k]}$.  This shows that it is entire.

To bound the conductor, we rewrite
\[
	L(s,\psi_{K_1/F} \otimes \psi_{K_2/F})
		= \frac{ \zeta_{K_1K_2}(s) / \zeta_{K_1K_2^N}(s) }{\zeta_{K_1^NK_2}(s)/\zeta_{K_1^NK_2^N}(s)}.
\]
Both the numerator and denominator are entire, by the Aramata--Brauer theorem.  Since $\psi_{K_1/F}\psi_{K_2/F}$ is a character of $G_F$, it follows that the conductor of $L(s,\psi_{K_1/F} \otimes \psi_{K_2/F})$ divides that of the numerator, which in turn divides $\mathfrak{D}_{K_1K_2/F}$, the relative discriminant of the compositum $K_1K_2/F$.  We have $\mathfrak{D}_{K_1K_2/F} \mid \mathfrak{D}_{K_1/F}^{[K_2:F]} \mathfrak{D}_{K_2/F}^{[K_1:F]}$, so the result follows by the conductor-discriminant formula and taking norms to $k$. Additionally taking norms to $\mathbb{Q}$ gives the result on the absolute discriminants.
\end{proof}

\subsection{Dirichlet series for completely multiplicative functions}
\label{sec:prelim_large_sieve}

Given an $n$-dimensional Artin $L$-function $L(s,\rho)$ over $k$ and a parameter $z$ depending at most on $n$, we introduce the {\it completely multiplicative Artin $L$-function}
\begin{equation}
\label{eqn:Lc_def}
L_z(s,\rho):=\prod_{\N\kp>z}(1-\lambda_{\rho}(\kp)\N\kp^{-s})^{-1}=\sum_{\kn}\frac{a_{\rho}(\kn)}{\N\kn^s},
\end{equation}
where $z>0$ is a parameter that we will choose to depend at most on $n$.  By construction, $a_{\rho}(\kn)$ is completely multiplicative and satisfies $a_{\rho}(\kp)=\lambda_{\rho}(\kp)$ for $\N\kp>z$ while $a_{\rho}(\kp)=0$ for $\N\kp\leq z$.
\begin{lemma}
	\label{lem:Daileda}
	Let $\rho$ be an $n$-dimensional Artin representation over $k$ whose $L$-function $L(s,\rho)$ is entire, and let $L_z(s,\rho)$ be as in \eqref{eqn:Lc_def}.  If $z$ is sufficiently large with respect to $n$, then there exists an Euler product
	\[
	H_z(s,\rho)=\prod_{\kp}H_{\kp}(s,\rho)
	\]
	such that if $\re(s)>\frac{1}{2}$, then:
	\begin{enumerate}
		\item $H_z(s,\rho)$ converges uniformly and absolutely,
		\item $L_z(s,\rho) = H_z(s,\rho)L(s,\rho)$,
		\item $H_z(s,\rho)$ is non-vanishing, and
		\item there exists a constant $B=B(n)>0$ such that $H_z(s,\rho)\ll_{n}(\re(s)-\tfrac{1}{2})^{-B_n}$.
	\end{enumerate}
\end{lemma}
\begin{proof}
The claimed identity for $L_z(s,\rho)$ holds for $\re(s)>1$ once we define
	\[
	H_{\kp}(s,\rho) = \begin{cases}
		\frac{L_{\kp}(s,\rho)^{-1}}{1-\lambda_{\rho}(\kp)\N\kp^{-s}}&\mbox{if $\N\kp>z$,}\\
		L_{\kp}(s,\rho)^{-1}&\mbox{otherwise.}
	\end{cases}
	\]
	Because $L(s,\rho)$ is assumed to be an entire Artin $L$-function and $z$ depends at most on $n$, it follows that
	\[
	\prod_{\N\kp\leq z}|H_{\kp}(s,\rho)|\ll_n 1
	\]
	for $\re(s)>\frac{1}{2}$.  Since each $\alpha_{j,\rho}(\kp)$ has modulus at most 1, it follows that if $\N\kp\leq z$, then $H_{\kp}(s,\rho)$ has no zero in the region $\re(s)>\frac{1}{2}$.
	
	When $\N\kp>z$, a tedious calculation similar to that in the proof of \cite[Proposition 2]{DK} shows that there exists a degree $n-2$ polynomial $f$ whose coefficients depend only on the $\alpha_{j,\rho}(\kp)$ such that
	\[
	H_{\kp}(s,\rho) = 1+\frac{\N\kp^{-2s} f(\N\kp^{-s})}{1-\lambda_{\rho}(\kp)\N\kp^{-s}}.
	\]
	If $z\geq 4n^2$, then for $\N\kp>z$ and $\re(s)>\frac{1}{2}$, we have $|1-\lambda_{\rho}(\kp)\N\kp^{-s}|\geq\frac{1}{2}$.  Since each $\alpha_{j,\rho}(\kp)$ has modulus at most one, it follows that $|f(\N\kp^{-s})|\ll_n 1$ when $\re(s)>\frac{1}{2}$.  Therefore,
	\[
	H_{\kp}(s,\rho) = 1+O_n(\N\kp^{-2\re(s)}).
	\]
	If $z$ is sufficiently large with respect to $n$, then we ensure that $H_{\kp}(s,\rho)\neq 0$ for $\re(s)>\frac{1}{2}$.  Furthermore, in this region, there exists a constant $A(n)>0$ depending at most on $n$ such that
	\[
	|H_{\kp}(s,\rho)|\leq (1+\N\kp^{-2\re(s)})^{A(n)}.
	\]
	The desired result now follows from the bound
	\[
	\prod_{\N\kp>z}|H_{\kp}(s,\rho)|\ll_n \zeta(2\re(s))^{A(n)},\qquad \re(s)>\frac{1}{2}.
	\]
\end{proof}

Let $K_1,K_2\in\mathfrak{F}_{k}^{G}$, and let $N$ be a normal subgroup of $G$.  We write
\begin{equation}
\label{eqn:The_Euler_Product_EP}
L_{\kp}(s,\psi_K) = \prod_{j=1}^{d}\Big(1-\frac{\alpha_{j,\psi_K}(\kp)}{\N\kp^s}\Big)^{-1}=:1+\sum_{j=1}^{\infty}\frac{\lambda_{\psi_K}(\kp^j)}{\N\kp^{js}},
\end{equation}
where
\[
d:=|G|-|G|/|N|.
\]
We define $\lambda_{\psi_K}(\kn)$ by
\[
L(s,\psi_K)=\prod_{\kp}L_{\kp}(s,\psi_K)=:\sum_{\kn}\frac{\lambda_{\psi_K}(\kn)}{\N\kn^s}.
\]
In particular, $\lambda_{\psi_K}(\kn)$ is a multiplicative function.

Suppose that $K_1\cap K_2=K_{1}^N\cap K_{2}^N$.  Since $\psi_{K_1}\otimes \psi_{K_2}$ is an Artin representation, there exist complex numbers $\alpha_{j_1,j_2,K_1\times K_2}(\kp)$ with modulus at most 1 such that
\[
L_{\kp}(s,\psi_{K_1}\otimes \psi_{K_2}) = \prod_{j_1=1}^{d}\prod_{j_2=1}^{d}\Big(1 - \frac{\alpha_{j_1,j_2,K_1\times K_2}(\kp)}{\N\kp^s} \Big)^{-1}=:1+\sum_{j=1}^{\infty}\frac{\lambda_{\psi_{K_1}\otimes\psi_{K_2}}(\kp^j)}{\N\kp^{js}}.
\]
It follows from our proof of \cref{lem:entire-two} that if $\kp\nmid \kD_{K_1/k}\kD_{K_2/k}$, then
\begin{equation}
\label{eqn:RS_Dirichlet_series}
\{\alpha_{j_1,j_2,K_1\times K_2}(\kp)\colon 1\leq j_1,j_2\leq d\}=\{\alpha_{j_1,K_1}(\kp)\alpha_{j_2,K_2}(\kp)\colon 1\leq j_1,j_2\leq d\}.
\end{equation}
As with $L(s,\psi_K)$, we write
\begin{equation}
	\label{eqn:dirichlet_RS}
	L(s,\psi_{K_1}\otimes\psi_{K_2})=\prod_{\kp}L_{\kp}(s,\psi_{K_1}\otimes\psi_{K_2})=:\sum_{\kn}\frac{\lambda_{\psi_{K_1}\otimes\psi_{K_2}}(\kn)}{\N\kn^s}.
\end{equation}

Parallel with \cref{lem:Daileda}, we provide a convenient factorization of $L(s,\psi_{K_1}\otimes\psi_{K_2})$ when $K_1\cap K_2=K_1^N\cap K_2^N$.  We define
\begin{equation}
\label{eqn:LRS_def}
L_z^{RS}(s,\psi_{K_1}\otimes\psi_{K_2}):=\sum_{\kn}\frac{a_{\psi_{K_1}}(\kn)a_{\psi_{K_2}}(\kn)}{\N\kn^s}.
\end{equation}
Recall that the definitions of the completely multiplicative functions $a_{\psi_{K_j}}(\kn)$ depend on the choice of $z$ in \cref{lem:Daileda}, which we assume is sufficiently large and depending at most on $d$.
\begin{lemma}
	\label{lem:Daileda_2}
	Suppose that $K_1\cap K_2=K_1^N\cap K_2^N$.  Let $z$ in \cref{lem:Daileda} be sufficiently large (depending at most on $d$).  There exists an Euler product
	\[
	H_z^{RS}(s,\psi_{K_1}\otimes\psi_{K_2})=\prod_{\kp}H_{\kp}^{RS}(s,\psi_{K_1}\otimes\psi_{K_2})
	\]
	such that if $\re(s)>\frac{1}{2}$, then:
	\begin{enumerate}
		\item $H_z^{RS}(s,\psi_{K_1}\otimes\psi_{K_2})$ converges uniformly and absolutely,
		\item $L_z^{RS}(s,\psi_{K_1}\otimes\psi_{K_2}) = H_z^{RS}(s,\psi_{K_1}\otimes\psi_{K_2})L(s,\psi_{K_1}\otimes\psi_{K_2})$, and
		\item there exists a constant $B_d'>0$, depending at most on $d$, such that
		\[
		|H_z^{RS}(s,\psi_{K_1}\otimes\psi_{K_2})|\ll_{d} (D_{K_1}D_{K_2})^{\epsilon/2}(\re(s)-\tfrac{1}{2})^{-B_d'}.
		\]
	\end{enumerate}
\end{lemma}
\begin{proof}
This is similar to \cref{lem:Daileda}.  Observe by complete multiplicativity that
\[
L_z^{RS}(s,\psi_{K_1}\otimes\psi_{K_2}) = \prod_{\N\kp>z}(1-\lambda_{\psi_{K_1}}(\kp)\lambda_{\psi_{K_2}}(\kp)\N\kp^{-s})^{-1}.
\]
We thus compute
\[
H_{\kp}^{RS}(s,\psi_{K_1}\otimes\psi_{K_2})=\begin{cases}
\displaystyle\frac{L_{\kp}(s,\psi_{K_1}\otimes\psi_{K_2})^{-1}}{1-\lambda_{\psi_{K_1}}(\kp)\lambda_{\psi_{K_2}}(\kp)\N\kp^{-s}}&\mbox{if $\N\kp>z$,}\vspace{.05in}\\
L_{\kp}(s,\psi_{K_1}\otimes\psi_{K_2})^{-1}&\mbox{otherwise.}
\end{cases}
\]
Because $K_1\cap K_2=K_1^N\cap K_2^N$ by hypothesis, $L(s,\psi_{K_1}\otimes\psi_{K_2})$ is an entire Artin $L$-function.  Since $z$ depends at most on $d$, it follows that
\[
\prod_{\N\kp\leq z}|H_{\kp}^{RS}(s,\psi_{K_1}\otimes\psi_{K_2})|\ll_d 1,\qquad\re(s)>\frac{1}{2}.
\]
We proceed as in \cite[Proposition 2]{DK} using \cref{lem:entire-two}.  For all $\epsilon>0$, we have
\[
\Big|\prod_{\substack{\N\kp>z \\ \textup{$\kp$ ramified}}}H_{\kp}^{RS}(s,\psi_{K_1}\otimes\psi_{K_2})\Big|\ll_{d,\epsilon}(D_{K_1}D_{K_2})^{\epsilon/2},\qquad\re(s)>\frac{1}{2}.
\]
For the other Euler factors, there exists a degree $d^2-2$ polynomial $f$ whose coefficients depend only on the $\alpha_{j_1,j_2,\psi_{K_1}\otimes\psi_{K_2}}(\kp)$ such that
	\[
	H_{\kp}^{RS}(s,\psi_{K_1}\otimes\psi_{K_2}) = 1+\frac{\N\kp^{-2s} f(\N\kp^{-s})}{1-\lambda_{\rho}(\kp)\N\kp^{-s}}.
	\]
	If $z\geq 4n^4$, then $|1-\lambda_{\rho}(\kp)\N\kp^{-s}|\geq\frac{1}{2}$ for $\N\kp>z$ and $\re(s)>\frac{1}{2}$.  Since $\alpha_{j_1,j_2,\psi_{K_1}\otimes\psi_{K_2}}(\kp)$ has modulus at most one, it follows that $|f(\N\kp^{-s})|\ll_d 1$ when $\re(s)>\frac{1}{2}$.  Therefore,
	\[
	H_{\kp}^{RS}(s,\psi_{K_1}\otimes\psi_{K_2}) = 1+O_d(\N\kp^{-2\re(s)}),\qquad \re(s)>\frac{1}{2}.
	\]
	Thus, there exists a constant $B_d'>0$ depending at most on $d$ such that
	\[
	|H_{\kp}^{RS}(s,\psi_{K_1}\otimes\psi_{K_2})|\leq (1+\N\kp^{-2\re(s)})^{B_d'}.
	\]
	The desired result now follows from the bound
	\[
	\prod_{\N\kp>z}|H_{\kp}^{RS}(s,\psi_{K_1}\otimes\psi_{K_2})|\leq  \zeta(2\re(s))^{B_d'}.
	\]
	for $\re(s)>\frac{1}{2}$.
\end{proof}

\subsection{Partial sums}

Let $\phi(t)$ be a smooth test function supported in a compact subset of $[-2,2]$, and suppose that $\phi(t)=1$ for $t\in[0,1]$ and $\phi(t)\in [0,1)$ otherwise.  The Laplace transform of $\phi$ is
\begin{equation}
\label{eqn:mellin}
\widehat{\phi}(s) = \int_{\R}\phi(y)e^{sy}dy.
\end{equation}
By construction, $\widehat{\phi}(s)$ is an entire function of $s$, and upon integrating by parts, we find for any integer $m\geq 0$ that
\begin{equation}
\label{eqn:test_fcn_bounds}
\widehat{\phi}(s)\ll_{\phi,m}e^{2|\mathrm{Re}(s)|}|s|^{-m}.
\end{equation}
Let $T\geq 1$. By Fourier inversion, for any $x>0$ and any $c\in\R$, we have the identity
\[
\phi(T\log x)=\frac{1}{2\pi i T}\int_{c-i\infty}^{c+i\infty}\widehat{\phi}(s/T)x^{-s}ds.
\]

\begin{lemma}
\label{lem:local_density}
Let $T,x\geq 1$ and $\epsilon>0$.  If $D_{K_1},D_{K_2}\leq Q$, then
	\begin{multline*}
	\Big|\sum_{\kn} a_{\psi_{K_1}}(\kn)a_{\psi_{K_2}}(\kn)\phi\Big(T\log\frac{\N\kn}{x}\Big)\Big|\\
	\ll_{|G|,[k:\Q],\phi,\epsilon}\begin{cases}
	\sqrt{x}(\log x)^{B_d'}Q^{|G|/2+\epsilon}  T^{d^2[k:\Q]}&\mbox{if $K_1\cap K_2=K_{1}^N\cap K_{2}^N$,}\\
	\frac{1}{T}x(\log x)^{d^2[k:\Q]-1}+\sqrt{x}T^{d^2[k:\Q]}&\mbox{if $K_1\cap K_2\neq K_1^N\cap K_2^N$.}
\end{cases}
	\end{multline*}
\end{lemma}
\begin{proof}
Let $\tau_m(n)$ be the $n$-th Dirichlet coefficient of $\zeta(s)^m$.  By the definition of $a_{\psi_K}(\kn)$ and the fact that $|\alpha_{j,\psi_K}(\kp)|\leq 1$, the sum we want to estimate is bounded in modulus by 
\begin{align*}
\sum_{n}\tau_{d^2[k:\Q]}(n)\phi\Big(T\log\frac{n}{x}\Big)&=\frac{1}{2\pi i T}\int_{3-i\infty}^{3+i\infty}\zeta(s)^{d^2[k:\Q]}\widehat{\phi}(s/T)x^s ds\\
&=\frac{1}{T}\mathop{\mathrm{Res}}_{s=1}\zeta(s)^{d^2[k:\Q]}\widehat{\phi}(s/T)x^s+\frac{1}{2\pi iT}\int_{\frac{1}{2}-i\infty}^{\frac{1}{2}+i\infty}\zeta(s)^{d^2[k:\Q]}\widehat{\phi}(s/T)x^sds\\
&=\frac{x}{T} P_{d^2[k:\Q]}(\log x)+O_{d,[k:\Q]} (\sqrt{x}T^{d^2[k:\Q]}).
\end{align*}
Here, $P_{d^2[k:\Q]}$ is a polynomial of degree $d^2[k:\Q]-1$ whose coefficients depend the Laurent series expansion of $\zeta(s)$ centered at $s=1$ and $\phi^{(j)}(1/T)$ for $0\leq j\leq d^2[k:\Q]-1$.  Since $T\geq 1$, it follows from \eqref{eqn:test_fcn_bounds} that $P_{d^2[k:\Q]}(\log x)\ll_{d,[k:\Q]} (\log x)^{d^2[k:\Q]-1}$.  This result holds for all $K_1,K_2\in\mathfrak{F}_k^G(Q)$; in particular, the result when $K_1\cap K_2\neq K_1^N\cap K_2^N$ follows.

Now, assume that $K_1\cap K_2=K_{1}^N\cap K_{2}^N$.  Note that in view of the preceding analysis, our proposed bound is trivial if $x\leq Q^{|G|/2+\epsilon}T^{d^2[k:\Q]}$.  Thus, we may assume that $x>Q^{|G|/2+\epsilon}T^{d^2[k:\Q]}$.  By \eqref{eqn:LRS_def} and \cref{lem:Daileda_2}, the sum we want to estimate equals
	\[
		\Big|\frac{1}{2\pi i T}\int_{1/2+\frac{1}{\log x}-i\infty}^{1/2+\frac{1}{\log x}+i\infty}H_z^{RS}(s,\psi_{K_1}\otimes\psi_{K_2})L(s,\psi_{K_1}\otimes\psi_{K_2})\widehat{\phi}(s/T)x^s ds\Big|.
	\]
	Thus, by \cref{lem:convexity} and \eqref{eqn:test_fcn_bounds}, the integral is
	\begin{align*}
	&\ll_{d,[k:\Q],\phi}  \frac{\sqrt{x}}{T}(\log x)^{B_d'}Q^{\frac{|G|}{2}+\epsilon} \int_{-\infty}^{\infty}(2+|t|)^{\frac{d^2[k:\Q]}{4}}\Big|\widehat{\phi}\Big(\frac{1}{T}\Big(\frac{1}{2}+\frac{1}{\log x}+it\Big)\Big)\Big|dt\\
	&\ll_{d,[k:\Q],\phi}  \frac{\sqrt{x}}{T}(\log x)^{B_d'}Q^{\frac{|G|}{2}+\epsilon} \int_{-\infty}^{\infty}(2+|t|)^{\frac{d^2[k:\Q]}{4}}\min\Big\{1,\frac{T^{d^2[k:\Q]+2}}{(2+|t|)^{d^2[k:\Q]+2}}\Big\}dt,
	\end{align*}
which is bounded as claimed.
\end{proof}

\subsection{A large sieve inequality for Artin representations}
\label{sec:large_sieve}

We use the results from the preceding subsections to prove a large sieve inequality for the coefficients $a_{\psi_K}(\kn)$ of $L_z(s,\psi_K)$.  We then apply our large sieve to bound the mean value of a certain Dirichlet polynomial which naturally arises in Montgomery's method of detecting zeros of $L$-functions \cite{MR0249375}.

Recall that $\mathfrak{F}_{k}^{G}$ is a set of distinct number fields $K$ which are Galois extensions of $k$, each with Galois group isomorphic to a fixed group $G$, $\mathfrak{F}_{k}^{G}(Q)=\{K\in\mathfrak{F}\colon D_K\leq Q\}$, and
\[
\mathfrak{m}_{k}^{G,N}(Q) := \max_{K_1 \in \mathfrak{F}_{k}^{G}(Q)}\#\{K_2\in\mathfrak{F}_{k}^{G}(Q)\colon K_1\cap K_2\neq K_{1}^N\cap K_{2}^N\}.
\]
Let $b:\mathcal{O}_k\to\mathbb{C}$ and $\beta: \mathfrak{F}_k^G(Q)\to\mathbb{C}$ be functions with $\ell^2$ norms $\|b\|_2$ and $\|\beta\|_2$ defined by
\[
\|b\|_2^2=\sum_{x<\N\kn\leq xe^{1/T}}|b(\kn)|^2,\qquad \|\beta\|_2^2=\sum_{K\in\mathfrak{F}_k^G(Q)}|\beta(K)|^2.
\]

\begin{theorem}
\label{thm:pre_large_sieve}
	Let $Q,T,x\geq 1$ and $\epsilon>0$.  Define
	\[
	C(Q,T,x):=\sup_{\substack{b \\ \|b\|_2\neq 0}}\sum_{K\in\mathfrak{F}_{k}^{G}(Q)}\Big|\sum_{\N\kn\in(x,xe^{1/T}]}a_{\psi_K}(\kn)b(\kn)\Big|^2\Bigg\slash \sum_{\N\kn\in(x,xe^{1/T}]}|b(\kn)|^2
	\]
There exists a constant $B_{G}''>0$, depending at most on $|G|$, such that
	\begin{align*}
	C(Q,T,x)\ll_{|G|,[k:\Q]}  (\log x)^{B_{G}'' [k:\Q]}\Big(\mathfrak{m}_{k}^{G,N}(Q)\frac{x}{T}+\sqrt{x}Q^{|G|+\epsilon}T^{|G|^2[k:\Q]}\Big).
	\end{align*}
\end{theorem}

\begin{proof}
By the duality principle of finite-dimensional Hilbert spaces, $C(Q,T,x)$ equals the supremum over all functions $\beta\colon \mathfrak{F}_k^G(Q)\to\mathbb{C}$ such that $\|\beta\|_2=1$ of
\begin{equation}
\label{eqn:ratio_2}
\sum_{\N\kn\in(x,xe^{1/T}]}\Big|\sum_{K\in\mathfrak{F}_{k}^{G}(Q)}a_{\psi_K}(\kn)\beta(K)\Big|^2.
\end{equation}
Fix a smooth function $\phi$ supported on a compact subset of $[-2,2]$, such that $\phi(T\log\frac{t}{x})$ is a pointwise upper bound for the indicator function of the interval $(x,xe^{1/T}]$.  Then \eqref{eqn:ratio_2} is
\begin{equation}
\label{eqn:ratio4}
\leq \sum_{\kn}\Big|\sum_{K\in\mathfrak{F}_{k}^{G}(Q)}a_{\psi_K}(\kn)\beta(K)\Big|^2 \phi\Big(T\log\frac{\N\kn}{x}\Big).
\end{equation}
Expanding the square and swapping the order of summation, we find that \eqref{eqn:ratio4} equals
\[
    \sum_{K_1,K_2\in\mathfrak{F}_{k}^{G}(Q)}\beta(K_1)\overline{\beta(K_2)}\sum_{\kn} a_{\psi_{K_1}}(\kn)\overline{a_{\psi_{K_2}}(\kn)}\phi\Big(T\log\frac{\N\kn}{x}\Big)
\]
Since $\psi_K$ is real-valued, we have $a_{\psi_K}(\kn)\in\R$.  Since $|\beta(K_1)\overline{\beta(K_2)}|\leq\frac{1}{2}(|\beta(K_1)|^2+|\beta(K_2)|^2)$ by the inequality of arithmetic and geometric means and $\|\beta\|_2=1$, the above display is
\[
\leq \max_{K_1\in\mathfrak{F}_{k}^{G}(Q)}\sum_{K_2\in\mathfrak{F}_{k}^{G}(Q)}\Big|\sum_{\kn}a_{\psi_{K_1}}(\kn)a_{\psi_{K_2}}(\kn)\phi\Big(T\log\frac{\N\kn}{x}\Big)\Big|.
\]

By \cref{lem:local_density}, the definition of $\mathfrak{m}_{k}^{G,N}(Q)$, and the fact that $d\leq|G|-1$, there exists a constant $B_{G}''>0$ (depending at most on $|G|$) such that
\[
\ll_{|G|,[k:\Q],\phi,\epsilon} (\log x)^{B_{G}''[k:\Q]}\Big(\mathfrak{m}_{k}^{G,N}(Q)\frac{x}{T}\hat{\phi}(1/T)+\sqrt{x}Q^{|G|/2+\epsilon}T^{|G|^2[k:\Q]}\#\mathfrak{F}_{k}^{G}(Q)\Big).
\]
Once we fix $\phi$, we have that $\widehat{\phi}(1/T)\ll 1$ by  \eqref{eqn:mellin}.  Since Schmidt \cite{Schmidt} proved that 
\[
\#\mathfrak{F}_{k}^{G}(Q)\ll_{|G|,[k:\Q]}(Q/D_k)^{(|G|+2)/4}
\]
and $|G|\geq 2$, the theorem follows.
\end{proof}

Let $\mu_{k}(\kn)$ be the $\kn$-th Dirichlet coefficient of $\zeta_k(s)^{-1}$.  Then $\mu_k(\kn)=0$ unless $\kn$ is squarefree, in which case $\mu_{k}(\kn)=(-1)^{\#\{\kp|\kn\}}$.

\begin{corollary}
\label{cor:MVT_primes}
	Let $Q,T\geq 1$.  There exists a constant $B_{G}'''>0$, depending at most on $|G|$, such that if $X:= Q^{2(|G|+\epsilon)}T^{2|G|^2[k:\Q]}$ and $\log Y\asymp_{|G|,[k:\Q]} \log X$, then
	\begin{align*}
		\sum_{K\in\mathfrak{F}_k^G(Q)}\int_{-T}^{T}\Big|\sum_{\substack{X<\N\kn\leq X^{\log Y}}}\frac{a_{\psi_K}(\kn)\mu_k(\kn)}{\N\kn^{1+\frac{1}{\log Y}+iv}}\Big|^2 dv&\ll_{|G|,[k:\Q],\epsilon}\mathfrak{m}_k^{G,N}(Q)(\log X)^{B_{G}'''[k:\Q]},\\
		\sum_{K\in\mathfrak{F}_k^G(Q)}\int_{-T}^{T}\Big|\sum_{\N\kn\leq X}\frac{a_{\psi_K}(\kn)\mu_k(\kn)}{\N\kn^{\frac{1}{2}+\frac{1}{\log Y}+iv}}\Big|^2 dv&\ll_{|G|,[k:\Q],\epsilon}\mathfrak{m}_k^{G,N}(Q) X(\log X)^{B_{G}'''[k:\Q]}.
	\end{align*}
\end{corollary}
\begin{proof}
	We prove the first bound; the second is proved identically.  A formal generalization of a result of Gallagher \cite[Theorem 1]{Gallagher} to number fields tells us that if $c(\kn)$ is a complex-valued function supported on the integral ideals of $F$ such that $\sum_{\kn}|c(\kn)|<\infty$, then
	\[
	\int_{-T}^{T}\Big|\sum_{\kn}c(\kn)\N\kn^{-it}\Big|^2 dt\ll T^2\int_0^{\infty}\Big|\sum_{\N\kn\in(x,xe^{1/T}]}c(\kn)\Big|^2\frac{dx}{x}.
	\]
		Let $X=Q^{2|G|}T^{2|G|^2[k:\Q]}$,  $\log Y\asymp_{|G|,[k:\Q]}\log X$, and
	\[
	b(\kn)=\begin{cases}
		\mu_k(\kn)\N\kn^{-1-\frac{1}{\log Y}}&\mbox{if $\N\kn\in[X,X^{\log Y}]$,}\\
		0&\mbox{otherwise.} 
	\end{cases}
	\]
	If $c(\kn)=a_{\psi_{K}}(\kn)b(\kn)$, then
	{\small\begin{equation}
	\label{eqn:comint}
	\sum_{K\in\mathfrak{F}_{k}^{G}(Q)}\int_{-T}^{T}\Big|\sum_{\N\kn\in(X,X^{\log Y}]}\frac{a_{\psi_K}(\kn)\mu_{k}(\kn)}{\N\kn^{1+\frac{1}{\log Y}+iv}}\Big|^2 dv\ll T^2\int_0^{\infty}\sum_{K\in\mathfrak{F}_{k}^{G}(Q)}\Big|\sum_{\N\kn\in(x,xe^{1/T}]}a_{\psi_{K}}(\kn)b(\kn)\Big|^2\frac{dx}{x}.
	\end{equation}}%
	We apply \cref{thm:pre_large_sieve} and bound the above display by
	\[
	\ll_{|G|,[k:\Q],\epsilon} \mathfrak{m}_{k}^{G,N}(Q)\sum_{\kn}|b(\kn)|^2 \N\kn\Big((\log \N\kn)^{B_{G}''[k:\Q]}+\frac{X}{\N\kn^{1/2}}\Big).
	\]
	Since $\log Y\asymp\log X$, Lemma 2.4 of \cite{Weiss} implies that \eqref{eqn:comint} is bounded by $\mathfrak{m}_k^{G,N}(Q)$ times
	\begin{align*}
	(\log X)^{B_{G}''[k:\Q]}\sum_{\N\kn\in[X,X^{\log Y}]}\frac{1}{\N\kn}\ll_{|G|,[k:\Q]}(\log X)^{B_{G}''[k:\Q]+1}\mathop{\mathrm{Res}}_{s=1}\zeta_k(s).
	\end{align*}
	The residue is $\ll_{[k:\Q]}(\log D_k)^{[k:\Q]-1}$ \cite[Theorem 1]{Louboutin}, and the result follows.
\end{proof}

\subsection{Proof of \cref{thm:ZDE}}

Let $K\in\mathfrak{F}_k^G(Q)$.  The bound
\[
\#\{\beta+i\gamma\colon \beta\geq 0,~|\gamma-t|\leq 1,~L(\beta+i\gamma,\psi_K)=0\}\ll_{|G|,[k:\Q]} \log Q+\log(|t|+2)
\]
holds for all $t\in\R$ by proceeding as in \cite[Proposition 5.7]{IK}.  Given $\sigma\in(\frac{1}{2},1)$, we decompose the rectangle $[\sigma,1]\times[-T,T]$ into disjoint boxes of the shape $[\sigma,1]\times[u,u+2(\log X)^2]$, where $X$ is as in \cref{cor:MVT_primes}.  Each of these boxes contains $\ll_{|G|,[k:\Q]}(\log Q T)^3$ zeros.  Writing $\eta_{\psi_K}$ for the number of smaller boxes which contain at least one zero of $L(s,\psi_K)$, then
\[
N_{K/K^N}(\sigma,T)\ll_{|G|,[k:\Q]} (\log Q T)^3 \eta_{\psi_K}.
\]

Since $L(s,\psi_K)$ is entire, it follows from \cref{lem:Daileda} if $\gamma\in\R$ and $\beta>\frac{1}{2}$, then $\beta+i\gamma$ is a zero of $L(s,\psi_K)$ if and only if it is a zero of $L_z(s,\psi_K)$.  Thus, we will detect the zeros of $L_z(s,\psi_K)$.   Since the coefficients $a_{\psi_K}(\kn)$ are completely multiplicative, the $\kn$-th Dirichlet coefficient of $L_z(s,\psi_K)^{-1}$ is $a_{\psi_K}(\kn)\mu_k(\kn)$.  With $X$ as in \cref{cor:MVT_primes}, we define
\[
Y := X^{\frac{2}{3-2\sigma}},\qquad M_X(s,\psi_K) := \sum_{\N\kn\leq X}\frac{a_{\psi_K}(\kn)\mu_{k}(\kn)}{\N\kn^s}.
\]
A straightforward computation shows that if $\beta+i\gamma$ is a nontrivial zero of $L(s,\psi_K)$ with $\beta\geq\sigma\geq \frac{1}{2}+(\log Y)^{-1}$, then
\begin{align*}
	e^{-1/Y}&=\frac{1}{2\pi i}\int_{1-\beta+\frac{1}{\log Y}-i\infty}^{1-\beta+\frac{1}{\log Y}+i\infty}(1-L_z(\beta+i\gamma+w,\psi_K)M_X(\beta+i\gamma+w,\psi_K))\Gamma(w)Y^w dw\\
	&+\frac{1}{2\pi i}\int_{\frac{1}{2}-\beta+\frac{1}{\log Y}-i\infty}^{\frac{1}{2}-\beta+\frac{1}{\log Y}+i\infty}L_z(\beta+i\gamma+w,\psi_K)M_X(\beta+i\gamma+w,\psi_K)\Gamma(w)Y^w dw.
\end{align*}
If such a zero exists, then at least one of the two integrals above must be large since $e^{-1/Y}=1+O(Y^{-1})$.  As in  Montgomery's method \cite{MR0249375} for Dirichlet characters, one uses Stirling's formula to show that
\begin{align*}
	\eta_{\psi_K}&\ll_{|G|,[k:\Q]} Y^{2(1-\sigma)}(\log Y)^2\int_{-T}^{T}|1-L_z(1+\tfrac{1}{\log Y}+iv,\psi_K)M_X(1+\tfrac{1}{\log Y}+iv,\psi_K)|^2 dv\\
	&+ Y^{\frac{1}{2}-\sigma}\int_{-T}^{T}|L_z(\tfrac{1}{2}+\tfrac{1}{\log Y}+iv,\psi_K)M_X(\tfrac{1}{2}+\tfrac{1}{\log Y}+iv,\psi_K)| dv\\
	&\ll_{|G|,[k:\Q]} Y^{2(1-\sigma)}(\log Y)^2\int_{-T}^{T}|1-L_z(1+\tfrac{1}{\log Y}+iv,\psi_K)M_X(1+\tfrac{1}{\log Y}+iv,\psi_K)|^2 dv\\
	&+ Y^{\frac{1}{2}-\sigma}\int_{-T}^{T}|L_z(\tfrac{1}{2}+\tfrac{1}{\log Y}+iv,\psi_K)|^2 dv + Y^{\frac{1}{2}-\sigma}\int_{-T}^{T}|M_X(\tfrac{1}{2}+\tfrac{1}{\log Y}+iv,\psi_K)|^2 dv.
\end{align*}
The second bound follows from the inequality of arithmetic and geometric means.

\cref{lem:convexity,lem:Daileda} imply that
\[
\int_{-T}^{T}|L_z(\tfrac{1}{2}+\tfrac{1}{\log Y}+iv,\psi_K)|^2 dv \ll_{|G|,[k:\Q]} Q^{1/2}T^{1+d[k:\Q]/2}(\log Y)^{2A_d}.
\]
Furthermore, by \cref{lem:convexity,lem:Daileda} again, if we temporarily write $L=L_z(1+\frac{1}{\log Y}+iv,\psi_K)$ and $M_X=M_X(1+\frac{1}{\log Y}+iv,\psi_K)$, then we have
\begin{align*}
	|1-LM_X|^2&=|L|^2|L^{-1}-M_X|^2\\
	&\leq |L|^2(|L^{-1}-M_{X^{\log Y}}|+|M_{X^{\log Y}}-M_X|)^2\\
	&\ll_{|G|,[k:\Q]} (\log QT)^{d[k:\Q]}(|L^{-1}-M_{X^{\log Y}}|^2+|M_{X^{\log Y}}-M_X|^2).
\end{align*}
A straightforward partial summation shows that $|L^{-1}-M_{X^{\log Y}}|^2\ll_{|G|,[k:\Q]}1$, and thus
\begin{multline*}
	N_{K/K^N}(\sigma,T)\ll_{|G|,[k:\Q]} \Big(Y^{\frac{1}{2}-\sigma}\Big[Q^{1/2}T^{1+|G|[k:\Q]/2}(\log QT)^{2A_d}+\int_{-T}^{T}\Big|\sum_{\N\kn\leq X}\frac{a_{\psi_K}(\kn)\mu_{k}(\kn)}{\N\kn^{\frac{1}{2}+\frac{1}{\log Y}+iv}}\Big|^2 dv\Big]\\
	+Y^{2(1-\sigma)}(\log QT)^{|G|[k:\Q]+2}\Big[1+\int_{-T}^{T}\Big|\sum_{X<\N\kn\leq X^{\log Y}}\frac{a_{\psi_K}(\kn)\mu_{k}(\kn)}{\N\kn^{1+\frac{1}{\log Y}+iv}}\Big|^2\Big]\Big)(\log QT)^{3}.
\end{multline*}
Finally, we sum over $K\in\mathfrak{F}_k^G(Q)$ and apply \cref{cor:MVT_primes}.  By the Schmidt bound $\#\mathfrak{F}_k^G(Q)\ll_{|G|,[k:\Q]}Q^{(|G|+2)/4}$ and our choices of $X$ and $Y$, we find that there exists a constant $C_{G}>0$, depending at most on $G$, such that
\[
\sum_{K\in\mathfrak{F}_k^G(Q)}N_{K/K^N}(\sigma,T)\ll_{[k:\Q],|G|,\epsilon} \mathfrak{m}_k^{G,N}(Q)(\log QT)^{C_{G}[k:\Q]}(Y^{\frac{1}{2}-\sigma}X+Y^{2(1-\sigma)}).
\]
which is bounded as desired when $\sigma>\frac{1}{2}$.  For $\sigma<\frac{1}{2}$, our results are trivial in view of the generalized Riemann--von Mangoldt asymptotic for the count of all zeros up to height $T$ (see \cite[Theorem 5.8]{IK}).

\section{Proofs of \cref{thm:zero-density-intro,thm:main_result}}
\label{sec:zero-free} 

We begin by recalling the zero-free region of Lagarias and Odlyzko \cite[Section 8]{LO}: $\zeta_K(s)$ does not vanish in the region
\begin{equation}
	\label{eqn:Std_ZFR}
	\re(s)\geq 1-\frac{\Cr{ZFR}}{\log(D_K(|\im(s)|+3)^{[K:\Q]})},
\end{equation}
apart from the possibility of a single real simple zero.  Since $\zeta_K(s)/\zeta_{K^N}(s)$ is entire, it follows that $\zeta_K(s)/\zeta_{K^N}(s)\neq 0$ in the region \eqref{eqn:Std_ZFR}, apart from the possibility of a single real simple zero.  

We will use \cref{thm:ZDE} to prove \cref{thm:zero-density-intro}, showing that for all $K\in\mathfrak{F}_{k}^{G}(Q)$ apart from a small exceptional subset, the ratio of Dedekind zeta functions $\zeta_K(s)/\zeta_{K^N}(s)$ has a much stronger zero-free region than \eqref{eqn:Std_ZFR}.  In particular, since we may assume that $C_G\geq 6$ in \cref{thm:ZDE}, we may choose $\epsilon=\frac{1}{3}-\frac{5}{3C_G[k:\Q]}$ in \cref{thm:ZDE}, which leads to
\begin{equation}
\label{eqn:zde_simplified}
\sum_{K\in\mathfrak{F}_{k}^{G}(Q)}N_{K/K^N}(\sigma,T)\ll_{|G|,[k:\Q]} \mathfrak{m}_{k}^{G,N}(Q)(QT^{|G|[k:\Q]})^{(\frac{25}{3}-\frac{5}{3C_G[k:\Q]})|G|(1-\sigma)}(\log QT)^{C_G [k:\Q]}.
\end{equation}

\begin{proof}[Proof of \cref{thm:zero-density-intro}]
Let $0<\epsilon<1$, and define $\delta := \epsilon/(20|G|)$.  For each integer $2\leq j\leq Q^{\epsilon/(6C_{G}[k:\Q])}+1$, we iteratively apply \eqref{eqn:zde_simplified} with
\[
T = T_j=e^j-3,\qquad \sigma=\sigma_j:=1-\frac{2\delta\log Q}{\log Q+|G|[k:\Q]\log(T_j+3)},
\]
discarding $O_{|G|,[k:\Q],\epsilon}(\mathfrak{m}_{k}^{G,N}(Q)Q^{2(\frac{25}{3}-\frac{5}{3C_G[k:\Q]})\delta|G|+\frac{\epsilon}{6}})$ exceptions at most $Q^{\frac{\epsilon}{6C_G[k:\Q]}}$ times.  This dyadically builds a zero-free region for all except $O_{|G|,[k:\Q],\epsilon}(\mathfrak{m}_{k}^{G,N}(Q)Q^{\epsilon})$ of the fields $K\in\mathfrak{F}_{k}^{G}(Q)$. Thus, for all except $O_{|G|,[k:\Q],\epsilon}(\mathfrak{m}_{k}^{G,N}(Q)Q^{\epsilon})$ fields $K\in\mathfrak{F}_{k}^{G}(Q)$, the ratio $\zeta_K(s)/\zeta_{K^N}(s)$ is holomorphic and non-vanishing in the region
\[
\re(s)\geq1-\frac{2\delta\log Q}{\log Q+|G|[k:\Q]\log(|\im(s)|+3)},\qquad |\im(s)|\leq \exp(Q^{\epsilon/(6C_{G}[k:\Q])}).
\]
Since $D_K\leq Q$, we may replace the above region with the more restrictive region
\begin{equation}
\label{eqn:super_nice_ZFR}
\re(s)\geq1-\frac{2\delta\log D_K}{\log D_K+|G|[k:\Q]\log(|\im(s)|+3)},\qquad |\im(s)|\leq \exp(D_K^{\epsilon/(6C_{G}[k:\Q])}).
\end{equation}
For $|\im(s)|>\exp(D_K^{\epsilon/(6C_{G}[k:\Q])})$, we have the zero-free region \eqref{eqn:Std_ZFR}.  The theorem follows once we combine \eqref{eqn:super_nice_ZFR} with \eqref{eqn:Std_ZFR}.
\end{proof}

\begin{proposition}
\label{prop:good_PNT} Let $K/k$ be a Galois extension of number fields. Assume $\chi$ is a character of $G = \Gal(K/k)$ which is the induction of a non-trivial 1-dimensional cyclic character of a subgroup of $G$. Let $0 < \epsilon < 1$, and let the region  $\Omega_{K}(\epsilon)$ be given by \eqref{eqn:nice_ZFR}.  If $D_K$ is sufficiently large with respect to $[k:\Q]$, $|G|$, and $\epsilon$, and $L(s,\chi)$ does not vanish in the region $\Omega_K(\epsilon)$, then there exists an effectively computable constant $\Cr{main}=\Cr{main}(|G|,[k:\Q],\epsilon)>0$ such that
	\[
	\Big|\sum_{\N\kp\leq x}\chi(\kp)\Big|\ll_{|G|,[k:\Q],\epsilon} x\exp(-\Cr{main}\sqrt{\log x}),\qquad x\geq (\log D_K)^{81|G|/\epsilon}.
	\]
\end{proposition}

\begin{proof}[Proof of \cref{thm:main_result} assuming \cref{prop:good_PNT}] 
Since we have proved \cref{thm:zero-density-intro,thm:zfr-transfer}, it suffices for us to prove \eqref{eqn:pnt_a.e.} using \cref{prop:good_PNT}.  Suppose that $\zeta_K(s)/\zeta_{K^N}(s)$ is non-vanishing in the region $\Omega_K(\epsilon)$ and that $\rho$ is an irreducible Artin representation of $K/k$ whose kernel does not contain $N$.  Assuming \cref{conj:restricted-induction}, there exist rational constants $c_\chi(\rho)$ such that
  \[
    L(s,\rho) = \prod_{H \subseteq G} \prod_{\substack{ \chi \in \mathrm{Irr}(H) \\ \mathrm{dim}\chi = 1 \\ \ker\chi \not\supseteq H \cap N}} L(s, \mathrm{Ind}_H^G \chi)^{c_\chi(\rho)},
  \]
  the inner summation running over $1$-dimensional characters of subgroups $H$ whose kernel does not contain $H \cap N$.  Note that $c_\chi(\rho) \ll_G 1$ for each $\chi$.  Taking logarithmic derivatives, it follows for all $x \geq 3$ that
  \begin{equation}
  \label{eqn:identity_of_partial_sums}
    \sum_{\N\kn \leq x} \chi_{\rho}(\kn) \Lambda_k(\kn)=\sum_{H\subseteq G}~ \sum_{\substack{\chi \in \mathrm{Irr}(H) \\ \mathrm{dim}\chi = 1 \\ \ker\chi \not\supseteq H \cap N}} c_{\chi}(\rho)\sum_{\mathrm{N}\kn\leq x}(\mathrm{Ind}_H^G\chi)(\kn)\Lambda_k(\kn),
  \end{equation}
  where $\Lambda_k(\kp^j)=\log\N\kp$ and $\Lambda_k(\kn)=0$ otherwise. For each  $\chi$ in the righthand sum, \cref{lem:abelian-transfer} implies that   $L(s,\mathrm{Ind}_H^G\chi)$ is  non-vanishing in $\Omega_K(\epsilon)$. The desired result now follows by \cref{prop:good_PNT}, partial summation, and \eqref{eqn:identity_of_partial_sums}. 
\end{proof}

To prove \cref{prop:good_PNT}, we use the following smooth function to count prime ideals $\kp$ with $\N\kp\leq x$.

\begin{lemma}\label{lem:WeightChoice}

For all $x \geq 3$ and $\Delta \in (0,1/4)$, there exists a continuous real-variable function $f(t) = f_{x,\Delta}(t)$ such that:
\begin{enumerate}
	\item $0 \leq f(t) \leq 1$ for all $t \in \R$, and $f(t) \equiv 1$ for $\tfrac{1}{2} \leq t \leq 1$.
	\item The support of $f$ is contained in the interval $[\tfrac{1}{2} - \frac{\Delta}{\log x}, 1 +  \frac{\Delta}{\log x}]$. 
	\item Its Laplace transform $F(z) = \int_{\R} f(t) e^{-zt}dt$ is entire and is given by
			\begin{equation*}	
				F(z) = e^{-(1+ \frac{\Delta}{\log x})z} \cdot \Big( \frac{1-e^{(\frac{1}{2}+\frac{\Delta}{\log x})z}}{-z} \Big) \Big( \frac{1-e^{\frac{\Delta z}{2  \log x}}}{- \frac{\Delta z}{2\log x}} \Big)^{2}.
			\end{equation*}
	\item Let $s = \sigma + i t, \sigma > 0,$ and $t \in \R$. Then
	\[
	|F(-s\log x)|\leq 
		e^{\sigma \Delta} x^{\sigma} \min\Big\{ 1,  \frac{1 + x^{-\sigma/2} }{|s|\log x}  \Big( \frac{4}{\Delta|s|} \Big)^{2} \Big\}. 
	\]
	Moreover, $1/2 < F(0) < 3/4$ and
		  \begin{equation}
		F(-\log x)  = \frac{x}{\log x} + O\Big(\frac{\Delta x  + x^{1/2}}{\log x} \Big).
		\label{eqn:WeightChoice_MainTerm} 
		  \end{equation}
	\item Let $s = -\tfrac{1}{2}+it$ with $t \in \R$. Then
	\[
	|F(-s\log x)| \leq  \frac{5 x^{-1/4}}{\log x} \Big( \frac{4}{\Delta}\Big)^{2} (1/4+t^2)^{-1}.
	\]
\end{enumerate}
\end{lemma}	
\begin{proof}
	This lemma and its proof can be found by taking $\ell=2$  in \cite[Lemma 2.2]{TZ3}. 	
\end{proof}

Our next lemma uses the weight function constructed in \cref{lem:WeightChoice} to establish a preliminary form of the prime number theorem for $L(s,\chi)$.  Define $\delta=\epsilon/(20|G|)$.  For convenience, we rewrite the zero-free region \eqref{eqn:nice_ZFR} as
\[
\Omega_K(\epsilon):=\{s\in\mathbb{C}\colon \re(s)\geq 1-\omega_K(|\im(s)|+3)\},
\]
where for $t\geq 3$ we define
\begin{equation}
\label{eqn:omega-K}
\omega_K(t)=\omega_K(t;\epsilon)=\begin{cases}
	\frac{2\delta\log D_K}{\log D_K+[k:\Q]\log t}&\mbox{if $\log 3\leq t\leq D_K^{\epsilon/(6C_G[k:\Q])}$,}\\
	\frac{\Cr{ZFR}}{\log D_K+|G|[k:\Q]\log t}&\mbox{if $\log t>D_K^{\epsilon/(6C_G[k:\Q])}$}
\end{cases}
\end{equation}

\begin{lemma}
\label{lem:error_term_data}
Recall the notation and hypotheses of  \cref{prop:good_PNT}. Let $\omega_K(t) = \omega_K(t;\epsilon)$ be as in \eqref{eqn:omega-K}, and define
	\[
	\eta_K(x) = \inf_{t\geq 3}(\omega_K(t)\log x+\log t).
	\]
	If $L(s,\chi)$ does not vanish in the region $\Omega_K(\epsilon)$, then
	\[
	\Big|\sum_{\N\kp\leq x}\chi(\kp)\Big|\ll_{[K:\Q]} \frac{x}{\log x}e^{-\eta_K(x)/8}\log(eD_K)+\frac{x^{3/4}}{\log x},\qquad x\geq \max\{3,(\log D_K)^{4}\}.
	\]
\end{lemma}

\begin{proof}
		To start, we record a basic observation  that will be often used:
		\begin{equation}
		[K:\Q] \ll \log D_K \leq  x^{1/4}.
		\label{eqn:BasicObservation}
		\end{equation}	
		The first bound is Minkowski's inequality; the second bound holds by assumption. 
		
		Select the weight function $f(\, \cdot \,) = f_{x,\Delta}( \, \cdot \, )$ from \cref{lem:WeightChoice} for any $x \geq 3$ and with
		\[
		\Delta = x^{-1/4} + \min\{ \tfrac{1}{8}, 8 e^{-\eta_{K}(x)/4}  \}.
		\]
		A calculation identical to that in \cite[Lemma 2.3]{TZ3} shows that
		\[
		\sum_{\N\kn\leq x}\Lambda_k(\kn)\chi(\kn)=\sum_{\kn}\Lambda_k(\kn)\chi(\kn)f\Big(\frac{\log\N\kn}{\log x}\Big)+O_{[K:\Q]}(\sqrt{x}+\Delta x).
		\]
		Since $\chi$ is the induction of a non-trivial 1-dimensional character, $L(s,\chi)$ is a Hecke $L$-function and hence entire.  Thus, by Mellin inversion, we have
		\begin{align*}
		&\sum_{\N\kn\leq x}\Lambda_k(\kn)\chi(\kn)\\
		&=\frac{\log x}{2\pi i} \int_{2-i\infty}^{2+i\infty} -\frac{L'}{L}(s,\chi) F(-s\log x) ds+O_{[K:\Q]}(\sqrt{x}+\Delta x)\\
		&=\log x\sum_{\rho}F(-\rho\log x)+\frac{\log x}{2\pi i}\int_{-\frac{1}{2}-i\infty}^{-\frac{1}{2}+i\infty} -\frac{L'}{L}(s,\chi) F(-s\log x) ds+O_{[K:\Q]}(\sqrt{x}+\Delta x),
		\end{align*}
		where $\rho$ ranges over the nontrivial zeros of $L(s,\chi)$. The standard bound
		\[
		-\frac{L'}{L}(s,\chi)\ll_{[K:\Q]} \log(D_K(|\im(s)|+3)),\qquad \re(s)=-\frac{1}{2}
		\]
		and the lower bound $\Delta\geq x^{-1/4}$ imply via \cref{lem:WeightChoice}(5) that
		\[
		\sum_{\kn}\Lambda_k(\kn)\chi(\kn)=\log x\sum_{\rho}F(-\rho\log x)+O_{[K:\Q]}(\sqrt{x}+\Delta x).
		\]
		By \cref{lem:WeightChoice}(4), we have
		\[
		\log x\sum_{|\rho| \leq 1/4}F(-\rho\log x)\ll \log x\sum_{|\rho| \leq 1/4} x^{1/4} \ll_{[K:\Q]} x^{1/4}(\log x) \log D_K \ll_{[K:\Q]} \sqrt{x}\log x
		\]
		by \eqref{eqn:BasicObservation}. For the zeros $\rho = \beta+i\gamma$ of $L(s,\chi)$ with $|\rho| \geq 1/4$, observe that our assumed zero-free region for $L(s,\chi)$ implies that
		\[
		\frac{x^{-(1-\beta)}}{(|\gamma|+3)} =e^{-((1-\beta)\log x+\log(|\gamma|+3))}\leq e^{-\eta_{K}(x)}.
		\]
		by definition of $\eta_{K}$.  Hence \cref{lem:WeightChoice}(4) and our choice of $\Delta$ yields the estimate
		\[
			(\log x) |F(-\rho \log x)| \ll \frac{x^{\beta}}{(|\gamma|+3)} \cdot \frac{\Delta^{-2}}{(|\gamma|+3)^2} \ll x e^{-\eta_{K}(x)} \cdot \frac{ e^{\eta_{K}(x)/2}}{(|\gamma|+3)^2}
		\]
		for $|\rho| \geq 1/4$.	Thus, summing over all zeros $\rho$ of $L(s,\chi)$, it follows that
		\[
		\log x \sum_{\rho} |F(-\rho \log x)| \ll_{[K:\Q]} x e^{-\eta_{K}(x)/2} \sum_{\rho} \frac{1}{(|\gamma|+3)^2} +  \sqrt{x}\log x. 
		\]
		Since $L(s,\chi)$ is a Hecke $L$-function with $C(\chi)\ll_{[K:\Q]} D_K$, it follows by standard estimates for it zeros \cite[Lemma 2.5]{TZ3} and \eqref{eqn:BasicObservation} that the above expression is
		\begin{equation*}
		\begin{aligned}		
		&\ll_{[K:\Q]} x e^{-\eta_{K}(x)/2} \sum_{T=1}^{\infty}~ \sum_{\substack{T-1 \leq |\im(\rho)| \leq T}} \frac{\log D_{K}+\log(T+3)}{T^2}  + \sqrt{x}\log x\\
		&\ll_{[K:\Q]}  x e^{-\eta_{K}(x)/2} \log(e D_{K}) + \sqrt{x}\log x.
		\end{aligned}
		\end{equation*}
		By our choice of $\Delta$ and  \eqref{eqn:BasicObservation}, this implies that
		\begin{equation}
			\label{eqn:Artin-computedpsi}
	\Big|\sum_{\N\kn\leq x}\Lambda_k(\kn)\chi(\kn)\Big| \ll_{[K:\Q]} xe^{-\eta_{K}(x)/4} \log(eD_{K})   +  x^{3/4}. 
		\end{equation}
		The contribution from the prime powers and ramified primes is $O(\sqrt{x}+\log(eD_K))$, so by partial summation \cite[Lemma 2.1 and Equation 5.3]{TZ3}, it follows that
		\[
		\Big|\sum_{\N\kp\leq x}\chi(\kp)\Big|\ll_{[K:\Q]} \frac{x}{\log x} \sup_{\sqrt{x} \leq y \leq x}( e^{-\eta_{K}(y)/4} ) \log(eD_{K})  + \frac{x^{3/4}}{\log x} + \log(eD_K). 
		\]
		From the definition of $\eta_K$, one can see that $\eta_K(y)$ is an increasing function of $y$ and also $\eta_K(x^{1/2}) \geq \frac{1}{2} \eta_K(x)$. Hence, as $\log D_K \leq x^{1/4}$, we conclude the desired result. 
\end{proof}

\begin{lemma}
\label{lem:error_term_data_bound}
Recall the notation and hypotheses of \cref{prop:good_PNT}.  There exist effectively computable constants $\Cr{main}=\Cr{main}(|G|,[k:\Q],\epsilon)>0$ and $\Cl[abcon]{main2}=\Cr{main2}(|G|,[k:\Q],\epsilon)>0$ such that if $D_K\geq \Cr{main2}$ is sufficiently large with respect to $|G|$, $[k:\Q]$, and $\epsilon$, and $x \geq (\log D_K)^{81|G|/\epsilon}$, then
\[
\frac{\log(eD_K)}{\log x}e^{-\eta_K(x)/8}\ll_{|G|,[k:\Q],\epsilon}\exp(-\Cr{main}\sqrt{\log x}).
\]
\end{lemma}
\begin{proof}
	For notational compactness, we introduce $\epsilon_0=\epsilon/(6C_{G}[k:\Q])$.  By the definition of $\eta_K(x)$ and \eqref{eqn:omega-K}, we have that
	\[
	\eta_K(x)\geq \min\Big\{\inf_{0\leq u\leq D_K^{\epsilon_0}}\Big(\frac{\epsilon(\log D_K)\log x}{10|G|(\log D_K+[k:\Q]u)}+u\Big),\inf_{u\geq D_K^{\epsilon_0}}\Big(\frac{\Cr{ZFR}\log x}{\log D_K+|G|[k:\Q]u}+u\Big)\Big\}.
	\]
	Define
	\[
	\phi_1(u,x):=\frac{\epsilon(\log D_K)\log x}{10|G|(\log D_K+[k:\Q]u)}+u,\qquad \phi_2(u,x):=\frac{\Cr{ZFR}\log x}{\log D_K+|G|[k:\Q]u}+u.
	\]

	Notice that the global infimum of $\phi_1(u,x)$ over $u\in (-(\log D_K)/[k:\Q],\infty)$ is at
	\[
	u=u_1 := \sqrt{\frac{\epsilon\log D_K}{10|G|[k:\Q]}\log x}-\frac{\log D_K}{[k:\Q]}.
	\]
	Thus, the value of $u\in [0,D_K^{\epsilon_0}]$ at which $\phi_1(u,x)$ attains its infimum lies in $\{0,u_1,D_K^{\epsilon_0}\}\cap[0,D_K^{\epsilon_0}]$.  Observe that $u_1\geq 0$ if and only if $x\geq D_K^{(10|G|)/(\epsilon[k:\Q])}$, in which case $\phi_1(u_1,x)\leq \phi_1(D_K^{\epsilon_0},x)$ because $u_1$ is the global minimum.  For this range of $x$, we compute
\[
\phi_1(u_1,x)=\sqrt{\frac{2\epsilon\log D_K}{5|G|[k:\Q]}\log x}-\frac{\log D_K}{[k:\Q]}\geq \sqrt{\frac{\epsilon\log D_K}{10|G|[k:\Q]}\log x}.
\]
We also compute that $\phi_1(0,x)\geq \phi_1(D^{\epsilon_0},x)$ if and only if
\[
x\geq D_K^{\frac{10|G|}{\epsilon[k:\Q]}}e^{\frac{10|G|D_K^{\epsilon_0}}{\epsilon}},
\]
a range in which we already established that $u_1\geq 0$.  Therefore, since $\phi_1(0,x)=\frac{\epsilon}{10|G|}\log x$, we conclude that
\[
\inf_{0\leq u\leq D_K^{\epsilon_0}}\phi_1(u,x)\geq \min\Big\{\frac{\epsilon}{10|G|}\log x,\sqrt{\frac{\epsilon\log D_K}{10|G|[k:\Q]}\log x}\Big\}.
\]

Next, notice that the global minimum of $\phi_2(u,x)$ over $u\in(-(\log D_K)/[k:\Q],\infty)$ is at
\[
u_2 = \sqrt{\frac{\Cr{ZFR}\log x}{|G|[k:\Q]}}-\frac{\log D_K}{|G|[k:\Q]}.
\]
Thus, $\phi_2(u,x)$ attains its infimum over $u\geq D_K^{\epsilon_0}$ at $u=\max\{D_K^{\epsilon_0},u_2\}$.  It follows from a straightforward calculation that
\[
\inf_{u\geq D_K^{\epsilon_0}}\phi_2(u,x)\geq \sqrt{\frac{\Cr{ZFR}\log x}{|G|[k:\Q]}}+D_K^{\epsilon_0}.
\]
We conclude from the analysis for $u\leq D_K^{\epsilon_0}$ and $u\geq D_K^{\epsilon_0}$ that if $D_K\geq\Cr{main2}$, then
\begin{multline*}
\frac{\log(eD_K)}{\log x}e^{-\eta_K(x)/8} \\
\leq \frac{\log(eD_K)}{\log x}\exp\Big(- \frac{1}{8}\min\Big\{\frac{\epsilon\log x}{10|G|},\sqrt{\frac{\epsilon\log D_K}{10|G|[k:\Q]}\log x},\sqrt{\frac{\Cr{ZFR}\log x}{|G|[k:\Q]}}+D_K^{\epsilon_0}\Big\}\Big).
\end{multline*}
The desired result follows once we ensure that $x\geq(\log D_K)^{81|G|/\epsilon}$.
\end{proof}

\begin{proof}[Proof of \cref{prop:good_PNT}] 
This follows from \cref{lem:error_term_data,lem:error_term_data_bound} with the same constant $\Cr{main}$.
\end{proof}

\section{Application to prime degree extensions}
\label{sec:degree-p}

In this section, we show how our results apply to the family of prime degree $p$ extensions of a number field $k$.  The normal closures $K/k$ of such fields have Galois groups $G$ that are transitive subgroups of the symmetric group $S_p$.  The properties of such groups are well understood.  For example, they must be primitive permutation groups, and are thus subject to many of the results in Dixon and Mortimer \cite{DixonMortimer}.  More than this, such groups are \emph{classified}; see Lemma \ref{lem:transitive-p} below.  However, as the properties of such groups are of vital importance to many of our applications, we begin by providing a succinct but complete proof of the classical fact that such groups always have a unique minimal normal subgroup, and that the fixed field $K^N$ of this normal subgroup is linearly disjoint from the degree $p$ extension that we started with.

\begin{lemma}\label{lem:degree-p}
If $G\subseteq S_p$ is a transitive subgroup, then $G$ has a unique minimal nontrivial normal subgroup $N$.  Also, if $H \subseteq G$ is the stabilizer of a point, then $[H : H \cap N] = [G : N]$.
\end{lemma}
\begin{proof}
	Since $G$ is a transitive group of prime degree, it is primitive.  Consequently, any nontrivial normal subgroup acts transitively, and thus has an element of order $p$.  If $G$ had two minimal normal subgroups, say $N_1$ and $N_2$, then $N_1$ and $N_2$ commute, since any commutator lies in the intersection $N_1 \cap N_2$, which is trivial since $N_1$ and $N_2$ are minimal.  It follows that $G$ would then have a subgroup isomorphic to $N_1 \times N_2$, which has order divisible by $p^2$.  Since the order of $G$ divides $p!$, which is not divisible by $p^2$, this cannot happen, so $G$ must have a unique minimal normal subgroup, $N$.  Finally, since $N$ is nontrivial, it acts transitively, so by the orbit-stabilizer theorem, we find $[N:H\cap N] = p = [G:H]$.  We conclude that $[H:H \cap N] = [G:H \cap N]/p = [G:N]$, as desired.
\end{proof}

Let $\mathrm{PSL}(n,q)$ and $\mathrm{P}\Gamma\mathrm{L}(n,q)$ denote the projective special linear and projective semilinear groups of rank $n$ over the finite field $\mathbb{F}_q$, respectively, and write $M_{11}$ and $M_{23}$ for the Mathieu groups of rank $11$ and rank $23$.  The transitive subgroups of $S_p$ and their unique minimal normal nontrivial subgroups are classified as follows.

\begin{lemma}\label{lem:transitive-p}
    Let $p$ be a prime, let $G \subseteq S_p$ be transitive, and let $N\unlhd G$ denote its unique minimal normal subgroup.  Then $G$, $N$, and $p$ satisfy one of the following:
    \begin{enumerate}
        \item[(a)] $N = \mathbb{Z}/p\mathbb{Z}$, $G \simeq \mathbb{Z}/p\mathbb{Z} \rtimes H$ for some $H \subseteq \mathrm{Aut}(\mathbb{Z}/p\mathbb{Z})$;
        \item[(b)] $N=G=A_p$ or $N = A_p$ and $G=S_p$;
        \item[(c)] $p=11$, $N=G = \mathrm{PSL}(2,11)$ or $N=G=M_{11}$;
        \item[(d)] $p=23$, $N=G=M_{23}$; or
        \item[(e)] there is some integer $n \geq 2$ and prime power $q$ for which $N=\mathrm{PSL}(n,q)$, $N \subseteq G \subseteq \mathrm{P}\Gamma\mathrm{L}(n,q)$, and $p = (q^n-1)/(q-1)$.
    \end{enumerate}
\end{lemma}
\begin{proof}
	This follows from the classification of finite simple groups.  See \cite[Corollary 4.2]{Feit}.
\end{proof}

As a consequence of Lemma \ref{lem:transitive-p}, we obtain the following.

\begin{lemma}\label{lem:transitive-p-applicable}
    Let $p$ be a prime, let $G$ be a transitive subgroup of $S_p$, and let $N \trianglelefteq G$ be its unique minimal normal subgroup.  Then either $G$ is solvable or the index of $N$ in $G$ is a prime power.  In particular, Hypothesis \hyperlink{hyp:T}{$\mathrm{T}(G,N)$} holds.
\end{lemma}
\begin{proof}
    In the first case of Lemma \ref{lem:transitive-p}, the group $G$ is monomial since every irreducible representation  is either $1$-dimensional or induced from the normal subgroup $\mathbb{Z}/p\mathbb{Z}$.  For cases (b)-(d), the index of $N$ is either $1$ or $2$.
    
    For case (e), write $q = \ell^m$ for some prime $\ell$ and integer $m \geq 1$.  Then $p = (\ell^{mn}-1)/(\ell^m-1)$, an expression that may be factored in terms of cyclotomic polynomials.  An elementary argument then shows that for $(\ell^{mn}-1)/(\ell^m-1)$ to be prime, $n$ must be prime and $m$ must be a power of $n$.  Next, since $(q^n-1)/(q-1)$ is a prime, it must be the case that $\mathrm{gcd}(n,q-1)=1$.  It follows that $\mathrm{PSL}(n,q)$ and $\mathrm{PGL}(n,q)$ coincide, and thus the quotient $\mathrm{P}\Gamma\mathrm{L}(n,q) / \mathrm{PSL}(n,q)$ is isomorphic to $\mathrm{Gal}(\mathbb{F}_q/\mathbb{F}_\ell) \simeq \mathbb{Z}/m\mathbb{Z}$.  Since $m$ is a power of the prime $n$, we conclude that $N$ must have prime power index in $G$.
\end{proof}

We next prove Corollary \ref{cor:approximate-dedekind}.

\begin{proof}[Proof of Corollary \ref{cor:approximate-dedekind}]
Let $G \subseteq S_p$ be transitive.  It follows by Lemma \ref{lem:degree-p} that $G$ admits a unique minimal nontrivial subgroup $N$.  Appealing to Theorem \ref{thm:zero-density-intro}, we find that for all except $O_{|G|,[k:\mathbb{Q},\epsilon}(Q^\epsilon)$ fields $K \in \mathfrak{F}_{k}^{G}(Q)$ that $\zeta_K(s)/\zeta_{K^N}(s)$ is non-vanishing in the region $\Omega_K(\epsilon)$.  For all $K$ outside the exceptional set, by Lemma \ref{lem:transitive-p-applicable} and Theorem \ref{thm:zfr-transfer}, it follows that $L(s,\rho)$ is holomorphic and non-vanishing in $\Omega_K(\epsilon)$ for each nontrivial irreducible Artin representation of $K/k$ whose kernel does not contain $N$.  (Such representations are in fact precisely the faithful representations of $G$, since $N$ is minimal, but we do not need this.)

Let now $F/k$ be a degree $p$ extension whose normal closure is $K$.  Lemma \ref{lem:degree-p} implies that $F \cap K^N = k$, since $F$ corresponds to one of the conjugate stabilizer subgroups of $G$.  This implies that $\zeta_F(s)/\zeta_k(s)$ may be decomposed as a product of irreducible Artin $L$-functions whose kernels do not contain $N$, since for example $\zeta_{K^N}(s)$ may be decomposed exactly as the product over the unfaithful representations of $G$.  Thus, $\zeta_F(s)/\zeta_k(s)$ is holomorphic and non-vanishing in $\Omega_K(\epsilon)$.

We prove the corollary for the family $\mathscr{F}_k^p$ by considering each of the finitely many transitive subgroups $G \subseteq S_p$ in turn.  For the family $\mathscr{F}_k^{n,S_n}$ with $n \geq 2$, it follows as above from Theorem \ref{thm:zero-density-intro} by taking $G = S_n$ and $N$ to be the unique minimal nontrivial normal subgroup of $G$.
\end{proof}

\section{An effective Chebotarev density theorem for fibers} \label{sec:chebotarev-proof}

In this section, we prove Theorem \ref{thm:fibered-chebotarev}.

\begin{proof}[Proof of Theorem \ref{thm:fibered-chebotarev}]
We begin by outlining our strategy in broad terms since the ultimate proof will be almost immediate once set up.  Let $K/k$ be a normal extension with Galois group $G$.  For any $x$, let $\Pi_K(x)$ be the class function on $G$ defined by
\begin{equation*} \label{eqn:class-function}
\Pi_K(x) = \sum_{\mathcal{C}} \frac{1}{|\mathcal{C}|}\pi_\mathcal{C}(x;K/k) \mathbf{1}_\mathcal{C},
\end{equation*}
where the summation runs over the conjugacy classes $\mathcal{C}$ of $G$ and $\mathbf{1}_\mathcal{C}$ denotes the indicator function of the class $\mathcal{C}$.  Because it is a class function, $\Pi_K(x)$ may be decomposed in terms of the irreducible characters of $G$, namely
\begin{equation} \label{eqn:class-function-decomp}
\Pi_K(x) = \sum_{\rho \in \mathrm{Irr}(G)} \langle \Pi_K(x), \chi_\rho\rangle \chi_\rho,
\end{equation}
where $\mathrm{Irr}(G)$ denotes the set of irreducible complex representations of $G$.  Since $\langle \Pi_K(x),\mathbf{1}_G\rangle$ is the average value of $\Pi_K(x)$ across $G$, we find
\[
\langle \Pi_K(x),\mathbf{1}_G\rangle = \frac{1}{|G|} \pi_k^{ur}(x),
\]
where $\pi_k^{ur}(x)$ denotes the number of primes of $k$ with bounded norm that are unramified in $K/k$.  In particular, the usual Chebotarev density theorem follows if $\langle \Pi_K(x),\chi_\rho\rangle$ is small for each nontrivial irreducible $\rho$, since then the difference between $\Pi_K(x)$ and $\langle \Pi_K(x),\mathbf{1}_G\rangle \mathbf{1}_G$ would be small by \eqref{eqn:class-function-decomp}.  

More generally, let $N \trianglelefteq G$ be a normal subgroup.  Then we may regard the analogous class function $\Pi_{K^N}(x)$ of $G/N$ as a class function on $G$, and we find
\[
\mathrm{Proj}_{\mathcal{R}_{\mathbb{C}}(G/N)} \Pi_K(x) = \frac{1}{|N|} \Pi_{K^N}(x),
\]
where $\mathrm{Proj}_{\mathcal{R}_\mathbb{C}(G/N)} \Pi_K(x)$ denotes the orthogonal projection of $\Pi_K(x)$ onto $\mathcal{R}_\mathbb{C}(G/N)$, the space of class functions on $G/N$.  Theorem \ref{thm:fibered-chebotarev} is true if and only if $\Pi_K(x)$ is ``close'' to $\frac{1}{|N|}\Pi_{K^N}(x)$.  Since $\mathcal{R}_\mathbb{C}(G/N)$ is spanned by the irreducible characters that factor through $G/N$, we find
\[
\Pi_K(x) - \mathrm{Proj}_{\mathcal{R}_{\mathbb{C}}(G/N)} \Pi_K(x)
	= \sum_{\substack{ \rho \in \mathrm{Irr}(G) \\ N \not\subseteq \ker \rho}} \langle \Pi_K(x), \chi_\rho\rangle \chi_\rho.
\]
Thus, to prove Theorem \ref{thm:fibered-chebotarev}, our goal is to show that $\langle \Pi_K(x), \chi_\rho\rangle$ is small for each $\rho$ whose kernel does not contain $N$.  In fact, we find
\[
\langle \Pi_K(x), \chi_\rho\rangle 
	= \langle \Pi_K(x), \bar\chi_\rho\rangle 
	= \frac{1}{|G|} \sum_{ \mathrm{N} \mathfrak{p} \leq x} \chi_\rho(\mathfrak{p}).
\]
Thus, 
Theorem \ref{thm:artin-pnt-later} may be used directly to control the inner product $\langle \Pi_K(x), \chi_\rho\rangle$.  Theorem \ref{thm:fibered-chebotarev} follows.
\end{proof}

\section{Applications to class groups} 
\label{sec:class-groups}

In this section, we prove \cref{thm:ell-torsion} on the $\ell$-torsion subgroups of class groups, and we prove Theorem \ref{thm:extremal-class-number} on the extremal order of class numbers.

\subsection{Bounds on $\ell$-torsion subgroups}

We begin by recalling the key lemma of \cite{EV}.

\begin{lemma}[Ellenberg--Venkatesh]\label{lem:ellenberg-venkatesh}
Suppose that $F/k$ is a degree $d$ extension of number fields, and let $\ell \geq 2$.  Suppose that $\mathfrak{P}_1,\dots,\mathfrak{P}_M$ are prime ideals of $F$ with norm at most $D_{F}^{1/(2\ell(d-1))-\delta}$ for some $\delta>0$ and that are not extensions of prime ideals from any proper subfield of $F/k$.  Then
\[
|\mathrm{Cl}(F)[\ell]| \ll_{k,\ell,\epsilon} D_F^{1/2+\epsilon}/M.
\]
\end{lemma}

\begin{proof}[Proof of Theorem \ref{thm:ell-torsion}]
Let $n\geq 2$ be an integer, $Q\geq 1$, and $F\in\mathscr{F}_k^{n,S_n}(Q)$.  For a number field $L$, let $\pi_L(x)$ be the prime ideal counting function for $L$.  Apart from at most $O_{n,[k:\Q],\epsilon}(Q^{\epsilon})$ such fields $F$, it follows from \eqref{eqn:partialsums_n} in \cref{cor:approximate-dedekind} that if $x\geq (\log D_F)^{81(n!)^2/\epsilon}$, we have
\[
\pi_F(x) = \pi_k(x)+\sum_{\N_{k/\Q}\kp\leq x}\chi_{\rho_F}(\kp) = \pi_k(x)+O_{n,[k:\Q],\epsilon}(x\exp(-\Cr{approxDedekind_n}\sqrt{\log x})).
\]
The contribution from the prime ideals of $F$ of degree larger than two is $\ll_n \sqrt{x}$, so this is absorbed by the error term.  Since there exist effectively computable constants $\Cl[abcon]{ck1}=\Cr{ck1}(k)>0$ and $\Cl[abcon]{ck2}=\Cr{ck2}(k)>0$ such that $\pi_k(x)\geq \Cr{ck1}x/\log x$ for all $x\geq\Cr{ck2}$, it follows that for all fixed $\delta>0$, there exists an effectively computable constant $\Cl[abcon]{nke}=\Cr{nke}(n,k,\delta,\epsilon)>0$ such that there are at least $\Cr{nke}D_F^{1/(2\ell(n-1))-\delta}/\log D_F$ degree one prime ideals with norm at most $D_F^{1/(2\ell(n-1))-\delta}$.  The result then follows from \cref{lem:ellenberg-venkatesh}, since degree $1$ primes are necessarily not the extension of a prime ideal in a proper subfield of $F$.  For a prime $p$, the same conclusion holds for all except $O_{p,[k:\Q],\epsilon}(Q^{\epsilon})$ fields $F\in\mathscr{F}_k^p(Q)$ by appealing to \eqref{eqn:partialsums_p} in \cref{cor:approximate-dedekind}.
\end{proof}

\subsection{The extremal order of class numbers}

We turn now to the proof of Theorem \ref{thm:extremal-class-number}.  Thus, let $r_1$ and $r_2$ be non-negative integers with $n:=r_1+2r_2$ at least $2$.  We wish to construct degree $n$ $S_n$-extensions $F/\mathbb{Q}$ of signature $(r_1,r_2)$ with large class number.  Our approach is inspired by the conditional work of Duke \cite{Duke03} in the totally real case, $r_2=0$.  In particular, Duke considered a slight modification of a family of polynomials first considered by Ankeny, Brauer, and Chowla \cite{ABC} for which there is an explicit full-rank subgroup of the units of the resulting fields.  Building on this, let $a_1,\dots,a_{r_1}$ and $b_{r_1+1},\dots,b_{r_1+r_2},c_{r_1+1},\dots,c_{r_1+r_2}$ be integers with each $b_j^2-4c_j < 0$ such that the polynomial
\begin{equation} \label{eqn:g-def}
    g(x) = \prod_{i=1}^{r_1} (x-a_i) \prod_{j=r_1+1}^{r_1+r_2} (x^2+b_j x + c_j)
\end{equation}
has only simple roots, the derivative $g^\prime(x)$ also only has simple roots $\beta_1,\dots,\beta_{n-1}$, and  $g(\beta_i) \neq g(\beta_j)$ for $i \neq j$.  That such integers exist follows, for example, by noting that these conditions are generic for real coefficients, that they may therefore be satisfied for rational coefficients by continuity, and then for integral coefficients by rescaling.  Now, consider a polynomial $f(t,x)$ over $\mathbb{Q}(t)$ defined by
\begin{equation}\label{eqn:f-def}
    f(t,x)
        := t^n g(x/t) - t 
        = \prod_{i=1}^{r_1} (x-a_i t) \prod_{j=r_1+1}^{r_1+r_2} (x^2 + b_j tx + c_j t^2) - t.
\end{equation}
The key properties of the polynomial $f(t,x)$ we shall need are the following.

\begin{lemma}\label{lem:abc-poly}
Let $g(x)$ and $f(t,x)$ be as defined above.  Then:
\begin{itemize}
    \item[(i)] The polynomial $f(t,x)$ is irreducible and has Galois group $S_n$ over $\mathbb{Q}(t)$.
    \item[(ii)] When $\tau \in \mathbb{Z}$ is squarefree and sufficiently large, the polynomial $f(\tau,x)$ is irreducible and the field $F = \mathbb{Q}(x)/f(\tau,x)$ has signature $(r_1,r_2)$, is totally ramified at primes $p \mid \tau$, and has regulator satisfying $\mathrm{Reg}_F \asymp_n (\log D_F)^{r_1+r_2-1}$ if $\mathrm{Gal}(\widetilde{F}/\mathbb{Q}) \simeq S_n$ and $\mathrm{Reg}_F \ll_n (\log D_F)^{r_1+r_2-1}$ in general.
\end{itemize}
\end{lemma}
\begin{proof}
We provide a complete proof since our family is slightly more general than his, but all of the essential ideas are due to Duke \cite{Duke03}.

When $\tau \neq \pm 1$ is squarefree, the polynomial $f(\tau,x)$ is Eisenstein and thus irreducible.  Thus, $f(t,x)$ must be irreducible over $\mathbb{Q}(t)$.  To compute its Galois group, we follow Duke \cite{Duke03} and note that the splitting field of $f(t,x)$ over $\mathbb{Q}(t)$ is the same as that of $t^nf(1/t, x/t) = g(x) - t^{n-1}$.  By the monodromy computation of \cite[Lemma 1]{Duke03}, it follows that the Galois group is generated by transpositions and is thus equal to $S_n$.

For (ii), the claim about irreducibility and ramification follows from the previous observation that $f(\tau,x)$ is Eisenstein.  For the claim about the signature of $F$, we note that as $\tau \to \infty$, the roots of $f(\tau,x\tau)$ approach those of $g(x)$, which was constructed to have signature $(r_1,r_2)$.  The polynomial $g(x)$ has simple roots, and as complex roots come in conjugate pairs, it follows that for sufficiently large $\tau$, $f(\tau,x)$, and therefore $F$, must have signature $(r_1,r_2)$ as well.  It also follows that to each root of $g(x)$ we may associate an embedding of $F$.  Explicitly, if we write $F = \mathbb{Q}(\xi)$ with $\xi$ an arbitrary root of $f(\tau,x)$, then for $j \leq r_1$, we associate the embedding $\sigma_j \colon F \to \mathbb{R}$ given by assigning $\xi$ the value of the root approximating $\tau a_j$, and for $r_1 + 1 \leq j \leq r_1+r_2$, we associate the embedding $\sigma_j \colon F \to \mathbb{C}$ such that $\sigma_j(\xi)/\tau$ approximates a root of $x^2+b_j x + c_j$.  For each $\sigma_j$, let $\gamma_j$ denote the associated root of $g(x)$.  Additionally, to each $\sigma_j$, we attach the usual absolute value $|\cdot|_j\colon F \to \mathbb{R}_{\geq 0}$ given by $|\alpha|_j = |\sigma_j(\alpha)|$ if $\sigma_j$ is real and $|\alpha|_j = |\sigma_j(\alpha)|^2$ if $\sigma_j$ is complex.

We now construct units $\epsilon_1,\dots,\epsilon_{r_1+r_2}$ in $F$ as follows.  For $i \leq r_1$, define $\epsilon_i = \tau(\xi-a_i\tau)^{-n}$, and for $r_1 +1 \leq i \leq r_1 + r_2$, define $\epsilon_i = \tau^2 (\xi^2 + b_i \tau \xi + c_i\tau^2)^{-n}$.  Then
\[
    \prod_{i=1}^{r_1+r_2} \epsilon_i
        = \tau^n (f(\tau,\xi)+\tau)^{-n} = 1.
\]
We also find that, in the ring $\mathcal{O}_F / (\tau)$, we have
\[
    (\xi-a_i \tau)^n 
        = \xi^n 
        = f(\tau,\xi)
        = 0
\]
and, using this, that $(\xi^2+b_i\tau \xi + c_i \tau^2)^n \equiv 0 \pmod{\tau^2}$.  It follows that each $\epsilon_i^{-1}$ is integral and a unit, so each $\epsilon_i \in \mathcal{O}_F^\times$.  We next find if $i \leq r_1$ and $j \neq i$ that as $\tau \to \infty$, then
\[
|\epsilon_i|_j
    \sim
    \begin{cases}
        \tau^{1-n} |\gamma_j-\gamma_i|^{-n}, & \text{if } j \leq r_1 \\
        \tau^{2-2n} |\gamma_j-\gamma_i|^{-2n}, & \text{if } r_1 + 1 \leq j \leq r_1 + r_2,
    \end{cases}
\]
while if $r_1 +1 \leq i \leq r_1 + r_2$ and $j \neq i$, then
\[
|\epsilon_i|_j
    \sim
    \begin{cases}
        \tau^{2-2n}|(\gamma_j-\gamma_i)(\gamma_j - \bar\gamma_i)|^{-n}, & \text{if } j \leq r_1, \\
        \tau^{4-4n}|(\gamma_j-\gamma_i)(\gamma_j - \bar\gamma_i)|^{-2n}, & \text{if } r_1+1 \leq j \leq r_1 + r_2.
    \end{cases}
\]
It then follows, 
using either the product formula or the relation $\epsilon_1 \dots \epsilon_{r_1+r_2} = 1$, that 
\[
    |\sigma_i(\epsilon_i)|
        \sim 
        \begin{cases} \displaystyle
            \tau^{(n-1)^2} \Big|\mathrm{Res}\Big(x-a_i,\frac{g(x)}{x-a_i}\Big)\Big|^{n}, & \text{if } i \leq r_1,\\
            \displaystyle
            \tau^{2(n-1)(n-2)} \Big|\mathrm{Res}\Big(x^2+b_ix+c_i,\frac{g(x)}{x^2+b_ix+c_i}\Big)\Big|^n, & \text{if } r_1 + 1 \leq i \leq r_1+r_2,\\
        \end{cases}
\]
where $\mathrm{Res}(\cdot,\cdot)$ denotes the polynomial resultant.  It follows that for $\tau$ sufficiently large, the $(r_1+r_2)\times (r_1+r_2)$ matrix
\[
    \begin{pmatrix}
        \log|\epsilon_1|_1 & \dots & \log|\epsilon_1|_{r_1+r_2} \\
        \vdots & & \vdots & \\
        \log|\epsilon_{r_1+r_2}|_1 & \dots & \log |\epsilon_{r_1+r_2}|_{r_1+r_2} \\
    \end{pmatrix}
\]
has positive diagonal entries, negative off-diagonal entries, and rows that sum to $0$.  It then follows by a lemma of Minkowski \cite{Minkowski} that any principal $(r_1+r_2-1) \times (r_1+r_2-1)$ minor has positive determinant.  Thus, any $r_1+r_2-1$ of the units $\epsilon_1, \dots, \epsilon_{r_1+r_2}$ are multiplicatively independent, thus forming a full-rank subgroup of $\mathcal{O}_F^\times$, and it follows from the above that $\mathrm{Reg}_F \ll_n (\log \tau)^{r_1+r_2-1} \ll_n (\log D_F)^{r_1+r_2-1}$, the latter inequality holding because $F$ is totally ramified at each prime dividing $\tau$.  The corresponding asymptotic lower bound follows from work of Remak \cite{Remak} (see also \cite{Silverman-Regulator}).
\end{proof}

By comparison with the cited work of Remak \cite{Remak}, Lemma \ref{lem:abc-poly} constructs fields whose regulators are essentially as small as possible.  To show that such fields often have class number as large as possible, we will show that the residue of the associated Dedekind zeta function $\mathrm{Res}_{s=1} \zeta_F(s)$ can be as large as conjecturally possible.  The key is to show that in many cases this residue may be approximated by a short Euler product.  We will do so in more generality.  Recall that the family $\mathfrak{F}_{k}^{G}$ is the set of normal extensions $K/k$ with Galois group isomorphic to $G$.  We show for the Artin representations $\rho$ of $K \in \mathfrak{F}_{k}^{G}$ considered in Theorem \ref{thm:main_result} that the value $L(1,\rho)$ may be approximated by a short Euler product.

\begin{proposition}\label{prop:euler-truncation}
Suppose for some $K \in \mathfrak{F}_{k}^{G}$, with $D_K$ sufficiently large with respect to $|G|$, $[k:\Q]$, and $\epsilon$, that $\zeta_K(s)/\zeta_{K^N}(s)$ is non-vanishing in the region $\Omega_K(\epsilon)$ for some normal subgroup $N \trianglelefteq G$ for which Hypothesis \hyperlink{hyp:T}{$\mathrm{T}(G,N)$} holds.  Let $\epsilon>0$.  For any Artin $L$-function $L(s,\rho)$ that does not factor through $K^N$, there holds for any $x \geq (\log D_K)^{81|G|/\epsilon}$, 
\[
L(1,\rho)
	= \Big(1 + O_{|G|,[k:\mathbb{Q}],\epsilon}\big(\exp(-\Cr{main}\sqrt{\log x})\big)\Big) \prod_{\mathrm{N}\kp \leq x} L_\kp(1,\rho),
\]
where $L_\kp(s,\rho)$ denotes the Euler factor of $L(s,\rho)$ at the prime $\mathfrak{p}$.  In particular, if $A>0$, then
\[
L(1,\rho)
    \asymp_{|G|,[k:\mathbb{Q}],A,\epsilon} \prod_{\mathrm{N}\kp \leq (\log D_K)^{A}} L_\kp(1,\rho).
\]
\end{proposition}
\begin{proof}
Write 
\[
\log L(s,\rho)
	=: \sum_{\mathfrak{n}} \frac{A_\rho(\mathfrak{n})}{\mathrm{N}\mathfrak{n}^s},
\]
where $A_\rho(\mathfrak{n})$ is supported on prime powers and $A_\rho(\mathfrak{p}) = \chi_\rho(\mathfrak{p})$.

We first claim that this series converges at $s=1$.  Indeed, by \eqref{eqn:pnt_a.e.} in \cref{thm:main_result} and partial summation, we have
\[
	\sum_{\mathrm{N} \mathfrak{p} \geq x} \frac{\chi_\rho(\mathfrak{p})}{\mathrm{N}\mathfrak{p}}
		\ll_{|G|,[k:\Q],\epsilon} \exp(-\Cr{main} \sqrt{\log x})
\] 
for any $x \geq (\log D_K)^{81|G|/\epsilon}$, and we trivially bound the contribution from prime powers by
\[
	\sum_{\substack{\mathrm{N}\mathfrak{p}^j \geq x \\ j \geq 2}} \frac{A_{\rho}(\mathfrak{p}^j)}{\mathrm{N}\mathfrak{p}^j}
		\ll_{[K:k]} \sum_{\substack{\mathrm{N}\mathfrak{p}^j \geq x \\ j \geq 2}} \frac{1}{j \mathrm{N}\mathfrak{p}^j}
		\ll_{[K:\mathbb{Q}]} x^{-1/2}.
\]
Thus, the series defining $\log L(s,\rho)$ converges at $s=1$.  It follows that it converges uniformly in sectors to the right of $s=1$, and that it converges to $\log L(1,\rho)$, where $L(1,\rho)$ is defined via the standard Dirichlet series.  Moreover, the above argument also shows that
\[
	\log L(1,\rho)
		= \sum_{\mathrm{N}\mathfrak{p} \leq x} \log L_\mathfrak{p}(1,\rho) + O(\exp(-\Cr{main}\sqrt{\log x})).
\]
The first claim follows upon exponentiating.  For the second claim, we apply Mertens's theorem to trivially bound the contribution from the range $(\log D_K)^A \leq \mathrm{N}\mathfrak{p} \leq (\log D_K)^{81|G|/\epsilon}$.
\end{proof}

We next show for families of fields obtained by specializing extensions of $\mathbb{Q}(t)$ that it is possible to choose the values of $L_\mathfrak{p}(1,\rho)$ at small primes to be as large as possible.  See also \cite[Proposition 3]{Duke03} and \cite{Serre-GaloisTheory}.

\begin{lemma}\label{lem:rigged-primes}
Let $f(t,x) \in \mathbb{Q}[x,t]$ be irreducible over $\mathbb{Q}(t)$.  There exists a constant $\Cl[abcon]{cf}=\Cr{cf}(f)>0$ such that for every prime $p > \Cr{cf}$, there exists $t_p \in \mathbb{Z}$ such that the polynomial $f(t_p,x) \pmod{p}$ splits completely.
\end{lemma}
\begin{proof}
As the polynomial $f(t,x)$ is irreducible over $\mathbb{Q}(t)$, it defines a separable extension of $\mathbb{Q}(t)$.  The normal closure of this extension must also be separable, and is hence cut out by an irreducible polynomial $g(t,x)$.  As $g(t,x)$ is irreducible, the curve $g(t,x) \equiv 0 \pmod{p}$ over $\mathbb{F}_p$ is non-singular for all but finitely many primes $p$.  Moreover, the Weil bound implies that there is a point $(t_p,x_p)$ on this curve for every sufficiently large $p$.  As $g(t,x)$ defines a normal extension of $\mathbb{Q}(t)$, it follows that $g(t_p,x)$ splits completely modulo $p$, and hence that $f(t_p,x)$ splits completely as well.
\end{proof}

We are now ready to prove Theorem \ref{thm:extremal-class-number}.  We do so in the following slightly stronger form.  For any possible signature $(r_1,r_2)$, we set $n = r_1+2r_2$ and let
\[
    \mathscr{F}_{\Q}^{r_1,r_2}(Q)
        := \{ F\in\mathscr{F}_{\Q}^{n,S_n}(Q)\colon  \mathrm{sgn}(F)=(r_1,r_2)\},
\]
where $\widetilde{F}/\mathbb{Q}$ denotes the normal closure of $F/\mathbb{Q}$.

\begin{theorem}\label{thm:extremal-class-number-body}
Let $r_1,r_2\geq 0$ be integers with $n = r_1 + 2 r_2 \geq 2$.  Let $0 < \epsilon < 1/(n!(n^2-n))$ and let $0 < \eta < 1/(n^2-n) $.  There exists a constant $c_{n,\eta,\epsilon}>0$ such that if $Q$ is sufficiently large in terms of $n$, $\eta$, and $\epsilon$, then for at least $c_{n,\eta,\epsilon} Q^{\frac{1}{n^2-n}-\eta}$ fields $F \in \mathscr{F}_{\Q}^{r_1,r_2}(Q)$, we have
\[
    |\mathrm{Cl}(F)|
        \asymp_{n,\epsilon} \frac{ D_F^{1/2} (\log\log D_F)^{n-1}}{(\log D_F)^{r_1+r_2-1}}
\]
and $\zeta_{\widetilde{F}}(s)/\zeta_{\widetilde{F}^{N}}(s)$ is non-vanishing in $\Omega_{\widetilde{F}}(\epsilon)$, where $N$ is the unique minimal nontrivial normal subgroup of $S_n$.
\end{theorem}

\begin{proof}
Let $(r_1,r_2)$ be a possible signature and consider the polynomial $f(t,x)$ defined by \eqref{eqn:f-def}.  Let $T$ be sufficiently large in terms of $n$.  By Lemma \ref{lem:rigged-primes}, there is a constant $c$ such that for every prime $p$ between $c$ and $(\log T)^{1/2}$, there is a congruence class $t_p \pmod{p}$ for which $f(t_p,x)$ splits completely $\pmod{p}$.  Define a congruence class $a \pmod{M}$ by
\[
M = \prod_{c \leq p \leq (\log T)^{1/2}} p, \quad\quad a \equiv t_p \pmod{p} \text{ for all } c \leq p \leq (\log T)^{1/2}.
\]
Note that $M = \exp(O(\sqrt{\log T}))$ by the prime number theorem.  

Consider squarefree $\tau \leq T$ for which $\tau \equiv a \pmod{M}$.  A quantitative version of the Hilbert irreducibility theorem due to Cohen \cite{Cohen} shows that for at most $O(T^{1/2} \log T)$ values $\tau$ the polynomial $f(\tau,x)$ does not cut out an $S_n$-extension of $\mathbb{Q}$.  Additionally, by Lemma \ref{lem:abc-poly}, if $f(\tau,x)$ cuts out a field $F$, then necessarily $\tau|D_F$.  Since $D_F = O_n(T^{n^2-n})$, it follows for any $\eta > 0$ that altogether there are $\gg_{n,\eta} T^{1-\eta}$ distinct degree $n$ $S_n$-extensions with signature $(r_1,r_2)$ produced in this way.

Let $F/\mathbb{Q}$ be such an extension and let $\widetilde{F}$ denote its normal closure over $\mathbb{Q}$.  Since $D_{\widetilde{F}} \leq D_F^{n!} \ll_n T^{n! (n^2-n)}$, it follows from Theorem \ref{thm:zero-density-intro} that for any $\epsilon < 1/(n!(n^2-n))$ all but $O_{n,\epsilon}(T^{\epsilon n!(n^2-n)})$ of the fields $F$ are such that $\zeta_{\widetilde{F}}(s)/\zeta_{\widetilde{F}^{N}}(s)$ is non-vanishing in the region $\Omega_{\widetilde{F}}(\epsilon)$.  For any such $F$, we find by Proposition \ref{prop:euler-truncation} that
\[
    \mathrm{Res}_{s=1}\zeta_F(s)
        \asymp_{n,\epsilon} \prod_{p \leq (\log T)^{1/2}} L_p(1,\rho_F) 
        \asymp_{n,\epsilon} (\log\log T)^{n-1} 
        \asymp_{n,\epsilon} (\log\log D_F)^{n-1}
\]
by our choice of the congruence class $a \pmod{M}$.  Using the estimate for the regulator provided by Lemma \ref{lem:abc-poly}, we conclude by the analytic class number formula that
\[
    |\mathrm{Cl}(F)|
        \asymp_{n,\epsilon} \frac{ D_F^{1/2} (\log\log D_F)^{n-1}}{(\log D_F)^{r_1+r_2-1}}.
\]
The result follows once we choose $Q = c T^{n^2-n}$ for a suitable constant $c$ depending on $n$.
\end{proof}

While the family of polynomials \eqref{eqn:f-def} provides a self-contained exposition, it is possible to obtain the conclusion of Theorem \ref{thm:extremal-class-number-body} for a larger set of fields by using different families while simultaneously imposing certain other constraints.  For example, by working instead with a family of polynomials considered by Bilu and Luca \cite{BiluLuca}, we obtain the following result for the family $\mathscr{F}_{\Q}^{n,0}$.

\begin{theorem}\label{thm:extremal-bilu-luca}
Let $n \geq 2$ and $\ell,\ell_0 \geq 2$ be integers.  Let $0 < \epsilon < 1/(2\ell n!(n-1))$ and $0 < \eta < 1/(2\ell_0(n-1))$.  There are positive constants $\Cl[abcon]{bilu_Luca1}$ and $\Cl[abcon]{bilu_Luca2}$ (depending at most on $n$, $\ell$, $\ell_0$, $\epsilon$, and $\eta$) such that if $Q\geq \Cr{bilu_Luca1}$, then at least $\Cr{bilu_Luca2}Q^{\frac{1}{2\ell(n-1)} - \eta}$ fields $F \in \mathscr{F}_{\Q}^{n,0}(Q)$ have the following properties:
\begin{itemize}
    \item[(i)] $\zeta_{\widetilde{F}}(s)/\zeta_{\widetilde{F}^{N}}(s)\neq 0$ in $\Omega_{\widetilde{F}}(\epsilon)$, where $N$ is the unique minimal normal subgroup of $S_n$;
    \item[(ii)] $|\mathrm{Cl}(F)| \displaystyle
        \asymp_{n,\ell,\epsilon} D_F^{1/2} (\log\log D_F)^{n-1}/(\log D_F)^{n-1}$;
    \item[(iii)] $\mathrm{Cl}(F)$ contains an element of exact order $\ell$; and
    \item[(iv)] $|\mathrm{Cl}(F)[\ell_0]|\ll_{n,\ell,\ell_0,\epsilon,\eta}D_F^{\frac{1}{2}-\frac{1}{2\ell_0(n-1)} + \eta}$.
\end{itemize}
\end{theorem}

\begin{proof}
Bilu and Luca consider the family of polynomials
\[
    f_\ell(t,x) = (x-a_1)\dots(x-a_{n-1})\Big(x - (-1)^{n-1} \frac{t^\ell - 1}{a_1\dots a_{n-1}}\Big) - 1,
\]
and show for any $\ell$ that there exist integers $a_1,\dots,a_{n-1}$ such that all but $O_{n,\ell}(T^{1/2} \log T)$ integers $|t| \leq  T$ subject to a fixed congruence condition cut out an extension $F \in \mathscr{F}_{n,0}( c T^{2n-2})$ with a point of order $\ell$ in the associated class group, and that any given field arises for at most $n(n-1)(n-2)$ values $t$.  Additionally, they also show for such values of $t$ that the regulator $\mathrm{Reg}_F$ of the associated field $F$ satisfies $\mathrm{Reg}_F \ll_{n} (\log D_F)^{n-1}$.  Appealing to Theorem \ref{thm:zero-density-intro}, Proposition \ref{prop:euler-truncation}, and Lemma \ref{lem:rigged-primes} as in the proof of Theorem \ref{thm:extremal-class-number-body}, along with Theorem \ref{thm:ell-torsion} for (iv), the result follows.  
\end{proof}

\section{Applications to subconvexity and periodic torus orbits} \label{sec:subconvexity}

\subsection{Subconvexity}

We begin with an bound for $L(\frac{1}{2}+it,\rho)$ depending only on $C(\rho,t)$ (see \eqref{eqn:analytic_conductor}) and the number $N_{\rho}(\sigma,T)$ of zeros $\beta+i\gamma$ of $L(s,\rho)$ such that $\beta\geq\sigma$ and $|\gamma|\leq T$.

\begin{lemma}
\label{lem:subconvex1}
Let $\rho$ be an $n$-dimensional Artin representation defined over a field $k$.  Suppose that $L(s,\rho)$ has a pole of order $0\leq r\leq n$ at $s=1$ and that $(s-1)^r L(s,\rho)$ is entire.  If $0\leq\Delta<\frac{1}{2}$ and $t\in\R$, then
\begin{align*}
\log |L(\tfrac{1}{2}+it,\rho)|&\leq \Big(\frac{1}{4}-\frac{\Delta}{10^9}\Big)\log(C(\rho)(|t|+1)^{n[k:\Q]})\\
&+\frac{\Delta}{10^7}\#\{\beta+i\gamma\colon \beta\geq 1-\Delta,~|\gamma-t|\leq 6,~L(\beta+i\gamma,\rho)=0\}+O_{n,[k:\Q]}(1)\\
&\leq \Big(\frac{1}{4}-\frac{\Delta}{10^9}\Big)\log(C(\rho)(|t|+1)^{n[k:\Q]})+\frac{\Delta}{10^7}N_{\rho}(1-\Delta,|t|+6)+O_{n,[k:\Q]}(1).
\end{align*}
\end{lemma}
\begin{proof}
Viewing $L(s,\rho)$ as a degree $n[k:\Q]$ $L$-function over $\Q$, our hypotheses imply that $L(s,\rho)$ is an $L$-function in the class $\mathcal{S}(n[k:\Q])$ defined by Soundararajan and Thorner in \cite[Section 1]{ST}.  As such, it follows from \cite[Theorem 1.1]{ST} that
\[
\log|L(\tfrac{1}{2},\rho)|\leq \Big(\frac{1}{4}-\frac{\Delta}{10^9}\Big)\log C(\rho)+\frac{\Delta}{10^7}N_{\rho}(1-\Delta,6)+2\log|L(\tfrac{3}{2},\rho)|+O((n[k:\Q])^2),
\]
where $C(\rho)$ is defined in \eqref{eqn:analytic_conductor}.  Following an observation of Heath-Brown in \cite{HB}, we have for any $t\in\R$ the bound
\begin{align*}
\log|L(\tfrac{1}{2}+it,\rho)|&\leq \Big(\frac{1}{4}-\frac{\Delta}{10^9}\Big)\log C(\rho,t)+2\log|L(\tfrac{3}{2}+it,\rho)|+O((n[k:\Q])^2)\\
&+\frac{\Delta}{10^7}\#\{\beta+i\gamma\colon \beta\geq 1-\Delta,~|\gamma-t|\leq 6,~L(\beta+i\gamma,\rho)=0\},
\end{align*}
where $C(\rho,t)$ is defined in \eqref{eqn:analytic_conductor}.  Since  $|\alpha_{j,\rho}(\kp)|\leq 1$ uniformly, we have that
\[
2\log|L(\tfrac{3}{2}+it,\rho)|+O((n[k:\Q])^2)\ll_{n,[k:\Q]}1.
\]
Moreover, we have the crude bound $\#\{\beta+i\gamma\colon \beta\geq 1-\Delta,~|\gamma-t|\leq 6,~L(\beta+i\gamma,\rho)=0\}\leq N_{\rho}(1-\Delta,|t|+6)$.  The result now follows from \eqref{eqn:t-bound}.
\end{proof}

\begin{proof}[Proof of \cref{thm:subconvexity-nice-families}]
We will prove the result for $F\in\mathscr{F}_k^{p}(Q)$, where $p$ is prime.  For integers $n\geq 2$, the corresponding proof for $F\in\mathscr{F}_k^{n,S_n}(Q)$ is essentially identical.  In what follows, let $G\subseteq S_p$ be a transitive subgroup, and let $N\unlhd G$ be its unique minimal nontrivial normal subgroup.

Let $0<\epsilon<1$, and let $\delta=\epsilon/(20|G|)$.  Let $K\in\mathfrak{F}_{k}^{G}$ with $D_K$ sufficiently large with respect to $|G|$, $[k:\Q]$, and $\epsilon$.  By \cref{lem:cyclic-transfer}, the Artin $L$-function $\zeta_K(s)/\zeta_{K^N}(s)$ is holomorphic and the corresponding Artin representation $\psi_K$ has dimension $d=|G|-|G/N|$. Assume that $\zeta_K(s)/\zeta_{K^N}(s)$ is non-vanishing in the region $\Omega_K(\epsilon)$.  When $D_K$ is sufficiently large with respect to $|G|$, $[k:\Q]$, and $\epsilon$, the region $\Omega_K(\epsilon)$ contains the rectangle
\[
[1-\delta,1]\times[-D_K^{1/(|G|[k:\Q])}-6,D_K^{1/(|G|[k:\Q])}+6].
\]
It follows that $N_{K/K^N}(1-\delta,D_K^{1/(|G|[k:\Q])}+6)=0$.  As the analytic conductor satisfies $C(\psi_K)\ll_{|G|,[k:\Q]}D_K$, \cref{lem:subconvex1} implies that if $|t|\leq D_K^{1/(|G|[k:\Q])}$, then
\[
\log\Big|\frac{\zeta_K(\frac{1}{2}+it)}{\zeta_{K^N}(\frac{1}{2}+it)}\Big|\leq \Big(\frac{1}{4}-\frac{\delta}{10^9}\Big)\log\Big(\frac{D_K}{D_{K^N}}(1+|t|)^{d[k:\Q]}\Big)+O_{|G|,[k:\Q]}(1).
\]
For $|t|>D_K^{1/(|G|[k:\Q])}$, we appeal to the fact that neither $\zeta_K(s)$ nor $\zeta_{K^N}(s)$ vanishes on the line $\re(s)=1$.  The same must hold for  $\zeta_K(s)/\zeta_{K^N}(s)$.  We may therefore apply \cref{lem:subconvex1} with $\Delta = 0$ so that if $|t|>D_K^{1/(|G|[k:\Q])}$, then
\begin{align*}
    \log\Big|\frac{\zeta_K(\frac{1}{2}+it)}{\zeta_{K^N}(\frac{1}{2}+it)}\Big|&\leq \frac{1}{4}\log\Big(\frac{D_K}{D_{K^N}}(1+|t|)^{d[k:\Q]}\Big)+O_{|G|,[k:\Q]}(1)\\
    &\leq \Big(\frac{1}{4}-\frac{\delta}{10^9}\Big)\log \frac{D_K}{D_{K^N}}+d[k:\Q]\Big(\frac{1}{4}+\frac{\delta}{10^9}\Big)\log(1+|t|)+O_{|G|,[k:\Q]}(1).
\end{align*}
Combining bounds for both ranges and noting $d \leq |G|-1$, we conclude for all $t \in \R$ that 
\begin{align*}
\Big|\frac{\zeta_K(\frac{1}{2}+it)}{\zeta_{K^N}(\frac{1}{2}+it)}\Big|&\ll_{|G|,[k:\Q]}\Big(\frac{D_K}{D_{K^N}}\Big)^{\frac{1}{4}-\frac{\epsilon}{2\cdot 10^{10}|G|}}(1+|t|)^{(\frac{|G|}{4}+1)[k:\Q]}.
\end{align*}
By \cref{lem:subconvex1} with $\Delta=0$, we have $|\zeta_{K^N}(\tfrac{1}{2}+it)|\ll_{|G|,[k:\Q]}D_{K^N}^{1/4}(1+|t|)^{|G/N|[k:\Q]/4}$.  If $\zeta_{K^N}(\tfrac{1}{2}+it)=0$, then so does $\zeta_{K}(\tfrac{1}{2}+it)$, so we may assume that $\zeta_{K^N}(\tfrac{1}{2}+it)\neq 0$.  Since $N$ is nontrivial and $D_{K^N} \leq D_K^{1/|N|}$, we conclude that
\begin{align*}
   |\zeta_K(\tfrac{1}{2}+it)|\ll_{|G|,[k:\Q]}D_K^{\frac{1}{4}-\frac{\epsilon}{2\cdot 10^{10}|G|}} D_{K^N}^{\frac{\epsilon}{2\cdot 10^{10}|G|}}(1+|t|)^{2|G|[k:\Q]}\ll_{|G|,[k:\Q]}D_K^{\frac{1}{4}-\frac{\epsilon}{4\cdot 10^{10}|G|}}(1+|t|)^{2|G|[k:\Q]}.
\end{align*}
The desired result now follows.

Let $F\in\mathscr{F}_k^{p}(Q)$ be a subfield of $K/k$, in which case $F\cap K^N=k$. We have the factorization $(s-1)\zeta_F(s) = (\frac{\zeta_F(s)}{\zeta_k(s)})((s-1)\zeta_k(s))$.  The function $(s-1)\zeta_k(s)$ is entire, and the ratio $\zeta_F(s)/\zeta_k(s)$ is holomorphic in the region $\Omega_K(\epsilon)$ by \cref{cor:approximate-dedekind}.  Consequently, the righthand side is a product of two functions which are holomorphic in $\Omega_K(\epsilon)$, and we have the bound
\begin{align*}
    &\#\{\beta+i\gamma\colon \beta\geq 1-\delta,~|\gamma-t|\leq 6,~\zeta_F(\beta+i\gamma)=0\}\\
    &\leq \#\{\beta+i\gamma\colon \beta\geq 1-\delta,~|\gamma-t|\leq 6,~\zeta_F(\beta+i\gamma)/\zeta_k(\beta+i\gamma)=0\}\\
    &+\#\{\beta+i\gamma\colon \beta\geq 1-\delta,~|\gamma-t|\leq 6,~\zeta_k(\beta+i\gamma)=0\}.
\end{align*}
By \cref{lem:subconvex1} with $\Delta = \delta$, we deduce that if $|t|\leq D_K^{1/(|G|[k:\Q])}$, then
\begin{align*}
    \log|\zeta_F(\tfrac{1}{2}+it)|&\leq \Big(\frac{1}{4}-\frac{\delta}{10^9}\Big)\log(D_F(1+|t|)^{[F:\Q]})\\
    &+\frac{\delta}{10^7}\#\{\beta+i\gamma\colon \beta\geq 1-\delta,~|\gamma-t|\leq 6,~\zeta_F(\beta+i\gamma)=0\}+O_{[F:\Q]}(1)\\
    &\leq \Big(\frac{1}{4}-\frac{\delta}{10^9}\Big)\log(D_F(1+|t|)^{[F:\Q]})\\
    &+\frac{\delta}{10^7}\#\{\beta+i\gamma\colon \beta\geq 1-\delta,~|\gamma-t|\leq 6,~\zeta_k(\beta+i\gamma)=0\}+O_{[F:\Q]}(1).
\end{align*}
By \cite[Proposition 5.7]{IK}, we find that
\[
10^{-7}\delta\#\{\beta+i\gamma\colon \beta\geq 1-\delta,~|\gamma-t|\leq 6,~\zeta_k(\beta+i\gamma)=0\}\ll |G|^{-1}\epsilon\log(D_k(1+|t|)^{[k:\Q]}),
\]
hence
\[
|\zeta_F(\tfrac{1}{2}+it)|\ll_{|G|,[k:\Q]}D_k^{O(\frac{\epsilon}{|G|})}D_F^{\frac{1}{4}-\frac{\epsilon}{2\cdot 10^{10}|G|}}(1+|t|)^{O(1+\frac{\epsilon[k:\Q]}{|G|})}.
\]
This proves the desired result for $\zeta_F(\frac{1}{2}+it)$ when $|t|\leq D_K^{1/[k:\Q]}$.  By arguing as we did for $\zeta_K(\frac{1}{2}+it)$ when $|t|>D_K^{1/[k:\Q]}$ using the convexity bound (which follows from \cref{lem:subconvex1} with $\Delta=0$), we arrive at the desired result for all $t$.
\end{proof}

\subsection{Equidistribution of periodic torus orbits}
\label{ssec:torus-proof}

As indicated in \cite[Section 1.6.3]{ELMV-Annals}, the equidistribution statement of Theorem \ref{thm:torus-orbits} follows once a suitable subconvexity bound is known.  We elaborate slightly on their general setup, with notation consistent with this paper rather than \cite{ELMV-Annals}, before specializing to the case of interest to us.  

Thus, let $k$ be a number field and let $S$ be a finite set of places of $k$ containing all archimedean places and such that the finite primes in $S$ generate the class group of $k$.  Let $k_S := \prod_{v \in S} k_v$ and let $\mathcal{O}_{k,S}$ denote the $S$-integers of $k$, i.e. the elements of $k$ that are integral away from primes in $S$.  To show that torus orbits inside $\mathrm{PGL}_p(\mathcal{O}_{k,S}) \backslash \mathrm{PGL}_p(k_S)$ associated to orders $\mathcal{O}$ in  degree $p$ extensions $F/k$ become equidistributed, Einsiedler, Lindenstrauss, Michel, and Venkatesh prove two key lemmas, namely \cite[Lemmas 13.3 and 13.4]{ELMV-Annals}, both relying on a subconvexity hypothesis stated formally as \cite[Equation (71)]{ELMV-Annals}. This subconvexity hypothesis is that for any Hecke character $\chi$ of $k$ ramified only at primes in $S$, there are constants $A$ and $\delta>0$ depending at most on $k$ and $p$ such that for any $t$, 
\begin{equation} \label{eqn:general-subconvexity-elmv}
    L(\tfrac{1}{2} + it, \mathrm{Ind}_{G_F}^{G_k} \mathbf{1} \otimes \chi)
        \ll_{[F:\Q]} (q_\chi \cdot (1+|t|))^A D_F^{\frac{1}{4}-\delta},
\end{equation}
where $q_\chi$ is the conductor of the $L$-function $L(s,\chi)$ and $\mathrm{Ind}_{G_F}^{G_k} \mathbf{1}$ denotes the induction of the trivial character of the absolute Galois group $G_F$ to the absolute Galois group $G_k$ of $k$.  In other words, the $L$-function $L(s, \mathrm{Ind}_{G_F}^{G_k} \mathbf{1} \otimes \chi)$ is the twist of the Dedekind zeta function of $F$ by the character $\chi$.

This general statement appears to be outside the scope of our methods, but in the special case that $k = \mathbb{Q}$ and $S$ consists only of the infinite place, \eqref{eqn:general-subconvexity-elmv} reduces to requiring
\begin{equation}\label{eqn:subconvexity-elmv}
    \zeta_F(\tfrac{1}{2}+it)
        \ll_{[F:\Q]} (1+|t|)^A D_F^{\frac{1}{4}-\delta}.
\end{equation}
This is provided by Theorem \ref{thm:subconvexity-nice-families} and Lemma \ref{lem:degree-p} by the assumptions of Theorem \ref{thm:torus-orbits}.   Thus, for orders $\mathcal{O}$ inside fields for which \eqref{eqn:subconvexity-elmv} holds, both Lemmas 13.3 and 13.4 of \cite{ELMV-Annals} hold.  Lemma 13.4 is used to control the ``escape of mass'' of the measure $\mu_\mathcal{O}$ associated to $\mathcal{O}$, and in particular it follows that any weak-* limit of the measures $\mu_\mathcal{O}$ for orders considered in Theorem \ref{thm:torus-orbits} must be a probability measure.  Lemma 13.3 is used to show that any weak-* limit is such that almost every ergodic component has positive entropy with respect to the action of a regular element in $H_p$.  (See also \cite[Theorem 1.9]{ELMV-Duke} for a weaker but more general statement.)  By a measure rigidity theorem of Einsiedler, Katok, and Lindenstrauss \cite{EKL} (also restated as \cite[Theorem 2.5]{ELMV-Annals}), any ergodic $H_p$-invariant measure on $\mathrm{PGL}_p(\mathbb{Z}) \backslash \mathrm{PGL}_p(\mathbb{R})$ with positive entropy must be Haar measure.  Thus, any weak-* limit of the $\mu_\mathcal{O}$ has almost every ergodic component given by Haar measure on $\mathrm{PGL}_p(\mathbb{Z}) \backslash \mathrm{PGL}_p(\mathbb{R})$, and equidistribution follows.

\section{Heuristics for the intersection multiplicity} \label{ssec:multiplicity} 

While our work is strongest for groups $G$ that have a unique minimal normal subgroup, our results may apply in other situations as well, provided that there is sufficient control over the intersection multiplicity $\mathfrak{m}_k^{G,N}(Q)$.  We therefore find it worthwhile to record a conjecture for how $\mathfrak{m}_k^{G,N}(Q)$ should grow; by taking $N=G$, this will also describe how the intersection multiplicity $\mathfrak{m}_k^G(Q)$ that was relevant in the previous works \cite{PTW,ThornerZaman} (either implicitly or explicitly) should grow.  We work in somewhat more generality.  

Fix a transitive and faithful permutation representation $\pi \colon G \to S_n$ with $n \geq 2$.  Given a field $K \in \mathfrak{F}_k^G$, such permutation representations correspond to subextensions of $K$ whose normal closure over $k$ is $K$ by taking the fixed field of a stabilizer subgroup.  Let $K^\pi$ denote the associated subextension, so for example $K^\pi = K$ when $\pi$ is the regular representation.  Define
	\[
		\mathfrak{F}_k^{G,\pi}(Q)
			:= \{ K \in \mathfrak{F}_k^G : D_{K^\pi} \leq Q\}.
	\]
Malle's conjecture \cite{Malle} predicts that
	\begin{equation} \label{eqn:Malle}
		Q^{\frac{1}{a_\pi(G)}} \ll_{k,G} \#\mathfrak{F}_k^{G,\pi}(Q) \ll_{k,G,\epsilon} Q^{\frac{1}{a_\pi(G)} + \epsilon},
	\end{equation}
where $a_\pi(G) := \min \{ n - \#\mathrm{Orb}_\pi(g) : g \neq \mathrm{id}\}$ and  $\mathrm{Orb}_\pi(g)$ denotes the set of orbits of the action of $\pi(g)$ on the set $\{1,\dots,n\}$.
	
Now, for a normal subgroup $N \unlhd G$, generalizing \eqref{eqn:normalized_intersection_multiplicity}, define
	\[
		\mathfrak{m}_k^{G,N,\pi}(Q):=\max_{K_1\in\mathfrak{F}_k^{G,\pi}(Q)}|\{K_2\in\mathfrak{F}_k^{G,\pi}(Q)\colon K_1\cap K_2\neq K_1^N\cap K_2^N\}|,
	\]
Thus, $\mathfrak{m}_k^{G,N,\pi}(Q)$ measures how often two fields $K_1,K_2 \in \mathfrak{F}_k^{G,\pi}(Q)$ have an intersection outside their associated subfields fixed by $N$.  In order for such an intersection to occur, there must exist normal subgroups $N_1,N_2 \unlhd G$ not containing $N$, possibly equal to each other, for which $K_1^{N_1} = K_2^{N_2}$.  Thus, $\mathfrak{m}_k^{G,N,\pi}(Q)$ will be bounded above by 
	\[
		\sum_{\substack{N^\prime \unlhd G \\ N \not\subseteq N^\prime}} \max_F \#\{K \in \mathfrak{F}_k^{G,\pi}(Q) : K^{N^\prime} = F\},
	\]
where the summation runs over the normal subgroups $N^\prime \unlhd G$ not containing $N$ and the maximum runs over all extensions $F/k$ inside the fixed choice of $\bar k$.  Notice that if $\mathfrak{p}$ is a tamely ramified prime in $K/k$, then $\mathfrak{p}$ is unramified in $K^{N^\prime}$ precisely when the inertia subgroup at $\mathfrak{p}$ is contained in $N^\prime$.  Motivated by the heuristic reasoning behind Malle's conjecture, set
	\[
		a_\pi(G,N^{\prime})
			:= \min \{ n - \#\mathrm{Orb}_\pi(g) : g \in N^\prime, g \neq \mathrm{id}\}.
	\]
If $F = K^{N^\prime}$ for some $K \in \mathfrak{F}_k^{G,\pi}(Q)$, we then expect
	\[
		Q^{\frac{1}{a_\pi(G,N^\prime)}} \ll_{F,G} \#\{K \in \mathfrak{F}_k^{G,\pi}(Q) : K^{N^\prime} = F\} \ll_{F,G,\epsilon} Q^{\frac{1}{a_\pi(G,N^\prime)} + \epsilon}.
	\]
Consequently, define
	\[
		m_\pi(G,N)
			:= \max_{\substack{N^\prime \unlhd G\\ N \not\subseteq N^\prime}} a_\pi(G,N^\prime)^{-1}
	\]
if there is at least one nontrivial such $N^\prime$, and define $m_\pi(G,N) = 0$ if there is no such $N^\prime$ (as is the case if either $N$ is the unique minimal normal subgroup of $G$, or if $N=G$ and $G$ is simple).  We then conjecture:
	\begin{conjecture} \label{conj:mgn-pi}
		With notation as above, as $Q \to \infty$, 
			\[
				Q^{m_\pi(G,N)} \ll_{k,G} \mathfrak{m}_k^{G,N,\pi}(Q) \ll_{k,G,\epsilon} Q^{m_\pi(G,N) + \epsilon}.
			\]
	\end{conjecture}
A few remarks are in order.  First, taking $\pi$ to be the (right) regular reprenentation of $G$, so that $\mathfrak{m}_k^{G,N,\pi}(Q) = \mathfrak{m}_k^{G,N}(Q)$, the orbits of $\pi(g)$ are exactly the (left) cosets of the cyclic subgroup $\langle g \rangle$.  There are $|G|/|\langle g\rangle |$ such orbits, and it follows that 
	\[
		m_{\mathrm{reg}}(G,N)
			= \min_{ \substack{ N^\prime \unlhd G \\ N \not\subseteq N^\prime}} \max_{ 1 \neq g \in N^\prime} \frac{|\langle g \rangle|}{|G|(|\langle g \rangle|-1)}.
	\]
Thus, for $\mathfrak{m}_k^{G,N}(Q)$ defined in \eqref{eqn:normalized_intersection_multiplicity}, Conjecture \ref{conj:mgn-pi} implies
	\begin{equation}
		\mathfrak{m}_k^{G,N}(Q) 
			\gg_{k,G} Q^{m_{\mathrm{reg}}(G,N)}.
	\end{equation}
Specializing further by taking $N=G$, this implies
	\begin{equation} \label{eqn:mg-conjecture}
		\mathfrak{m}_k^G(Q) 
			\gg_{k,G} Q^{m(G)}, \quad \text{where } m(G) = \min_{ \substack{ N \unlhd G \\ 1 \neq N \neq G}} \max_{1 \neq g \in N} \frac{|\langle g \rangle|}{|G|(|\langle g \rangle|-1)}.
	\end{equation}

Second, we now consider when we should expect $\mathfrak{m}_k^G(Q)$ to be $\gg_{k,G,\epsilon} Q^{-\epsilon} \#\mathfrak{F}_k^G(Q)$ for every $\epsilon > 0$, or more generally when $\mathfrak{m}_k^{G,N,\pi}(Q) \gg_{k,G,\epsilon} Q^{-\epsilon} \#\mathfrak{F}_k^{G,\pi}(Q)$.  Comparing Conjecture \ref{conj:mgn-pi} with Malle's conjecture \eqref{eqn:Malle}, it is apparent we should expect this to hold whenever there is a normal subgroup $N^\prime$ not containing $N$ for which $a_\pi(G,N) = a_\pi(G)$.  Focusing again on the regular representation, this will occur whenever there is a non-identity element of minimal order in $G$ contained in such an $N^\prime$, and when $N=G$, whenever there is an element of minimal order contained in any nontrivial proper normal subgroup.  This situation arises, for example, if $G=S_n$ for some $n \geq 4$, as $A_n$ contains elements of order $2$, and this leads to our aforementioned speculation that $\mathfrak{m}_k^{S_n}(Q) \gg_{k,n,\epsilon} Q^{-\epsilon} \#\mathfrak{F}_k^{S_n}(Q)$ for every $\epsilon > 0$.  

Though it does not follow directly from the above conjectures, we also remark that there are groups $G$ for which one should even expect $\mathfrak{m}_k^{G}(Q) \gg_{k,G} \#\mathfrak{F}_k^G(Q)$.  In particular, this is what one should expect if every non-identity element of minimal order in $G$ is contained in the same proper normal subgroup; this is ensured, for example, if $G$ is not a $p$-group and the Sylow $p$-subgroup of $G$ is normal for the smallest prime divisor $p$ of $|G|$.

Lastly, we comment on the role of the ramification restrictions imposed by Pierce, Turnage-Butterbaugh, and Wood.  As described above, for certain groups $G$, they fix conjugacy invariant subsets $\mathscr{R}_G \subseteq G$ and ask that all tamely ramified primes have inertia subgroups in $\mathscr{R}_G$.  An analogous story to the above holds, but with the key quantity $a_\pi(G,N^\prime)$ replaced by 
	\[
		a_\pi(G,N;\mathscr{R}_G)
			:= \min\{ n - \#\mathrm{Orb}_\pi(g) : g \in N^\prime \cap \mathscr{R}_G, g \neq \mathrm{id}\}.
	\]
For the groups they consider ($S_n$ for $n \geq 3$, the alternating group $A_4$, the dihedral group $D_p$, and the the cyclic group $C_n$ for any $n$), the set $\mathscr{R}_G$ is chosen specifically so that the intersection $N^\prime \cap \mathscr{R}_G$ is empty for every proper normal subgroup $N^\prime \unlhd G$.  Thus, with the analogous definition, one should expect $\mathfrak{m}_k^G(Q,\mathscr{R}_G) \ll_{k,G,\epsilon} Q^\epsilon$ for every $\epsilon > 0$.  We note, however, that apart from $G = C_p$, none of these groups are simple, so \eqref{eqn:mg-conjecture} implies that we should nonetheless expect for the full family without restrictions on inertia that $\mathfrak{m}_k^G(Q) \gg_{k,G} Q^{m(G)}$, where $m(G) > 0$.  

As touched upon earlier, a key advantage of our work is that all of the non-abelian groups $G$ considered in \cite{PTW} have a unique minimal normal subgroup $N$, and thus unconditionally we have $\mathfrak{m}_k^{G,N}(Q) = 1$ independently of any of the conjectural analysis above.  This is the case for many other groups as well. 

\bibliographystyle{abbrv}
\bibliography{GeneralizedLinnik}

\end{document}